\documentclass[a4paper,11pt]{article}
\usepackage{color,xcolor,ucs}
\usepackage[top=0.3in, bottom=0.5in, left = 0.65in, right = 0.65in]{geometry}
\usepackage[linkcolor=black,colorlinks=true,urlcolor=blue]{hyperref}
\usepackage{mathtools}

\usepackage{color,xcolor,ucs}
\usepackage{mathtools}   \usepackage{tikz} 
\usepackage{ amssymb }
\usepackage{extarrows} 
\usepackage{pgf,tikz}
\usepackage{float}
\usetikzlibrary{positioning}
\usetikzlibrary{shapes.geometric}
\usetikzlibrary{shapes.misc}
\usetikzlibrary{arrows}
\usepackage{caption}
\usepackage{mathrsfs}
\usetikzlibrary{arrows,shapes,automata,backgrounds,petri,positioning}
\usetikzlibrary{decorations.pathmorphing}
\usetikzlibrary{decorations.shapes}
\usetikzlibrary{decorations.text}
\usetikzlibrary{decorations.fractals}
\usetikzlibrary{decorations.footprints}
\usetikzlibrary{shadows}
\usetikzlibrary{calc}
\usetikzlibrary{spy}
\usepackage{amsmath}
\usepackage{array}
\usepackage{ amssymb }
\usepackage{braket}
\usepackage{qcircuit}
\usepackage{soul}
\usepackage{braket} 
\usepackage{relsize}

\usepackage{amsmath}
\usepackage{ amssymb }
\usepackage{braket}
\usepackage{qcircuit}
\usepackage{soul}
\usepackage{braket}

\title{{\LARGE Scaling limit of the triangular kinetic prudent walk}}
\author{Pete Rigas \footnote{Newport Beach, California 92625, ORCID: 0009-0003-1053-9720}}
\date{}

\begin{document}

\maketitle

\begin{abstract}
 Expressions for scaling limits of random walks, such as those obtained in several areas of the Probability theory literature, are of great significance in characterizing long term, stationary behavior of random processes. Presumably, in the limit of extremely long periods of time random walks and other stochastic processes including percolation are expected to fall into several \textit{universality classes}; such classes are of great significance in being able to disregard local, microscopic, dynamics, in favor of similar stationary dynamics in the bulk over long times. While several probabilistic models of interest are expected to behave similarly with respect to time, differences in fluctuations of limit shapes, and other stationary profiles, are expected to emerge. Given expressions obtained for the scaling limit of the kinetic prudent walk that have been obtained by Beffara, Friedli, and Velenik, in addition to the scaling limit for the kinetic uniform walk that have been obtained by Petrelis, Sun and Torri, we answer several questions pertaining to limit shapes, for the triangular prudent walk. While one would expect that the kinetic prudent walks, either over the square or triangular lattices, would have differing limit shapes, determining whether any differences between the scaling limit normalizations, and related quantities, emerge is of interest. In comparison to previously obtained scaling limits, we incorporate symmetries of the triangular lattice for computing the scale limit, beginning with a suitable basis. Steps for obtaining the scaling limit include: computing the probability that the prudent walker exists a suitably defined region, as well as computing conditionally defined probabilities that the prudent walker travel some fixed horizontal and vertical distances, after some arbitrary, strictly positive, time. Altogether, as the time $t \longrightarrow + \infty$, we argue that a similar scaling limit to those that have been previously obtained holds. \footnote{\textbf{MSC Class}: 60G50; 60G17; 82C41}
\end{abstract}


\section{Introduction}

\subsection{Overview}

\noindent Random walks are amongst the most classically studied objects in Probability theory. As a class of  objects that can be conditioned to avoid any step, ie transition, to a site of the underlying lattice that has been occupied at a previous time, prudent random walks, in comparison to ordinary random walks, are of interest to study as they share relations with polymers and closely related topics. While seminal results within the Probability theory literature have quantified the dependency of random walks on initial starting configurations, particularly through mixing rates and accompanying qualitative behaviors, other aspects of random walks remain of interest to explore. Namely, besides the mixing rate, it is of related interest to determine how the \textit{stationary behavior} of a prudent random walk emerges after having taken an infinite number of steps. Albeit the fact that computations for the scaling limit of such walks have been previously obtained in seminal work, {\color{blue}[1]}, it remains of interest, from a question of Beffara, Friedli, and Velenik, to determine how such limits would differ depending upon the \textit{underlying lattice} over which steps of the prudent random walk are taken. That is, with respect to an infinite number of steps, or close to an infinite number of steps, would there be any parameters of the prudent triangular walker that would degenerate? The computations provided by Beffara, Friedli and Velenik, for computing the scaling limit of the ordinary prudent walk, rely upon: (1) exploiting the symmetries of $\textbf{Z}^2$ for computing the probabilities that horizontal, and vertical, excursions of a random walk are larger than $k^{\frac{4}{3}}$, up to a constant; (2) making use of a suitably defined corners process for truncating horizontal, and vertical, excursions of the prudent walker; (3) computing the $L^1$ norm of a normalization in the scaling limit, which is related to a convergence occurring at a $\frac{3}{7} $ rate; (4) concluding that the desired scaling limit holds by sending the number of steps, and hence the time, of the prudent walker to $+ \infty$.

While Beffara, Friedli and Velenik remark that similar expressions for the scaling limit are expected to hold, several differences emerge. First and foremost, while the underlying lattice over which the prudent walker takes steps may initially not appear to impact the expression for the desired scaling limit, one must introduce several modifications for the triangular kinetic scaling limit, including: (1) a Green's function adapted for $\textbf{T}$; (2) computation of a \textit{system}, in comparison to a single, Radon-Nikodym derivative; (3) relating the behavior of the Radon-Nikodym derivative to a suitable effective walk representation. Over $\textbf{Z}^2$, the effective walk representation has been used for concluding that, under a straightforward computation involving the Radon-Nikodym derivative, that a desired normalization for the scaling limit holds. For the triangular kinetic prudent walk considered in this work, an adaptation of the Radon-Nikodym derivative computation initially presented for $\textbf{Z}^2$ can be adapted for $\textbf{T}$.

While one expects that the Radon-Nikodym derivative, whether corresponding to a random walk over $\textbf{Z}^2$ and $\textbf{T}$, would dictate the total number of steps in a path up to a certain time, the derivative over $\textbf{T}$ presents striking differences. Over $\textbf{T}$, one of the derivatives in the Radon-Nikodym derivative normalizes the \textit{slab} measure of the triangular effective random walk. The \textit{slab} measure, namely the measure supported over a slab over $\textbf{T}$, encodes characteristic information of the triangular scaling limit. While previous results for the uniform kinetic walk over $\textbf{Z}^2$, {\color{blue}[4]}, have compared the difference between a strictly positive scaling limit normalization constant $\mathscr{L}$ with the occupation time of $\mathrm{log}  \mathscr{L}  $ excursions, over $\textbf{T}$ we establish the following sequence of approximations for the scaling limit normalization. First, we introduce Brownian motion over the triangular lattice with the following diffusive scaling.

\bigskip

\noindent \textbf{Definition} \textit{0} (\textit{Brownian motion over the triangular lattice given a mesh size $a>0$ and diffusive scaling}). For the probability measure $\textbf{P}_{\textbf{T}} \big[ \cdot \big]$,

\begin{align*}
\textbf{P}_{\textbf{T}} \big[ \cdot \big]  : = \underset{T \in \textbf{T}}{\bigcup } \textbf{P}_{T} \big[ \cdot \big]  , 
\end{align*}

\noindent supported over the triangular lattice, denote $\big\{ X_k \big\}_{k \geq 0 } $ as the Markov chain,

{\small \begin{align*}
  \textbf{P}_{\textbf{T}} \bigg[         X_{k+1} - X_k = a \mathrm{exp} \big[ \frac{ i \pi j}{3}              \big] \bigg] = \frac{1}{6}  , 
\end{align*} }

\noindent for $0 \leq j \leq 5$. Over the triangular lattice,

{\small \begin{align*}
 \textbf{T} : =     \bigg\{ \forall   a > 0 , \exists m , n \in \textbf{Z} : a \big( m  + n \mathrm{exp} \big[  \frac{i \pi}{3} \big]   \big)                \bigg\}      \subsetneq \textbf{R}^2   , 
\end{align*}}

\noindent with nearest neighbors,

\begin{align*}
  \bigg\{  0 , \pm \frac{\pi}{3} , \pm \frac{2 \pi }{3} , \pi  \bigg\} , 
\end{align*}

\noindent \textit{Brownian motion},

\begin{align*}
    \mathcal{B}^{(a)}_t : = \textit{Brownian motion at time $t>0$ associated with mesh size $a>0$ of $\textbf{T}$} , 
\end{align*}

\noindent is obtained through the following diffusive scaling,

\begin{align*}
   \mathcal{B}^{(a)}_t : =  a X_{\lfloor   \frac{t}{a^2}  \rfloor }  . 
\end{align*}

\noindent For the below result, we compare the magnitude scaling limit normalization $\mathscr{L}$ over $\textbf{T}$ with the time that is spend at $\mathrm{log} \big[ \mathscr{L}^2 \big]$ many excursions. In the second result below, given the time spent by the random walk over $\textbf{T}$ we quantify whether factors of $\mathrm{log} \big[ \mathscr{L}^2 \big]$, $\mathrm{log} \big[ \mathscr{L}^4 \big]$, or $\mathrm{log} \big[ \mathscr{L}^6 \big]$ appear in the following indicator functions. Moreover, in comparison to scaling limits that have previously been identified over $\textbf{Z}^2$, {\color{blue}[4]}, \textit{isotropy} is not obtained through the limit as the scaling limit tends to $+\infty$; rather, for the triangular scaling limit isotropy is \textit{already} present taking the limit, despite the fact that the result still asserts a convergence to Brownian motion (see \textbf{Theorem} in \textit{1.2}).

\begin{itemize}
   \item[$\bullet$] \textit{Comparing the magnitude of the scaling limit normalization,} $\mathscr{L}$, \textit{over the triangular lattice to the total time spent at} $\mathrm{log} \big[ \mathscr{L}^2 \big]$ \textit{many excursions}. Denote,

    \begin{align*}
     \mathscr{T} : = \big[  \textit{Time spent along the first degree of freedom}     ,  \textit{Time spent along the second degree of freedom}     \big] \\ = \big[ \mathcal{T} , \widetilde{\mathcal{T}} \big]  ,
    \end{align*}

   \noindent from which one has that,

    \begin{align*}
    \frac{\mathscr{L} - \underset{1 \leq i \leq \mathrm{log} [ \mathscr{L}^2 ]}{\sum} \mathscr{T}_i }{\big( \mathrm{log} \mathscr{L} \big)^4} \longrightarrow 0     ,
    \end{align*}

\noindent as $\mathscr{L} \longrightarrow + \infty$.

   \item[$\bullet$] \textit{Comparing the time spent at  $\mathrm{log} \big[ \mathscr{L}^4 \big]$ many excursions of the triangular random walk}. Denote the indicator function,

   \begin{align*}
     \textbf{1}_{\big\{  \mathcal{T} + \widetilde{\mathcal{T}} \geq \mathrm{log} [ \mathscr{L}^4  ]     \big\}}      ,
   \end{align*}

   \noindent corresponding to the occurrence of the event $\big\{  \mathcal{T} + \widetilde{\mathcal{T}} \geq \mathrm{log} [ \mathscr{L}^4  ]   \big\}$, and,

    \begin{align*}
     \mathscr{L} : = \big[ \mathscr{L}_1 , \mathscr{L}_2 \big]   ,
    \end{align*}

   \noindent from which one has that,

{\tiny \begin{align*}
        \bigg\{    \frac{\mathscr{L}_1 - \underset{1 \leq i \leq \mathrm{log} [ \mathscr{L}^2 ]}{\sum} \mathcal{T}_i }{\big( \mathrm{log} \mathscr{L} \big)^2} \longrightarrow 0 ,     \frac{\mathscr{L}_2 - \underset{1 \leq i \leq \mathrm{log} [ \mathscr{L}^2 ]}{\sum} \widetilde{\mathcal{T}_i}  }{\big( \mathrm{log} \mathscr{L} \big)^2} \longrightarrow 0    \bigg\} \Longleftrightarrow     \bigg\{          \frac{\mathscr{L} - \underset{1 \leq i \leq \mathrm{log} [ \mathscr{L}^2 ]}{\sum} \mathscr{T}_i }{\big( \mathrm{log} \mathscr{L} \big)^4} \longrightarrow 0     \bigg\}           , 
\end{align*}}

\noindent as $\mathscr{L} \longrightarrow + \infty$.

\item[$\bullet$] \textit{Comparing the time spend at $\mathrm{log} \big[ \mathscr{L}^6 \big]$ many excursions of the effective triangular random walk}. Denote,

\begin{align*}
    \widetilde{\widetilde{\mathcal{T}}} > 0 , 
\end{align*}

\noindent corresponding to the time spent about the third degree of freedom of the effective triangular random walk, and the indicator function,

\begin{align*}
  \textbf{1}_{\big\{   \mathcal{T} + \widetilde{\mathcal{T}} +    \widetilde{\widetilde{\mathcal{T}}} \geq \mathrm{log} [ \mathscr{L}^6 ]    \big\}}  ,
\end{align*}

\noindent corresponding to the occurrence of the event $\big\{   \mathcal{T} + \widetilde{\mathcal{T}} +    \widetilde{\widetilde{\mathcal{T}}} \geq \mathrm{log} [ \mathscr{L}^6 ]    \big\}$, from which one has that,

   {\small  \begin{align*}
    \frac{\mathscr{L}^3  - \underset{1 \leq i \leq \mathrm{log} [ \mathscr{L}^2 ]}{\sum} \big[  \mathcal{T}_i + \widetilde{\mathcal{T}_i} +    \widetilde{\widetilde{\mathcal{T}_i}} \big]  }{\big( \mathrm{log} \mathscr{L} \big)^8} \longrightarrow 0    .
    \end{align*} }

\end{itemize}

\noindent The system of Radon-Nikodym derivatives appearing in expressions provided in \textit{2.2} for slab measures determines the normalization appearing in the scaling limit. Furthermore, the Radon-Nikodym derivatives, in comparison to \textit{the} Radon-Nikodym derivative for random walks over $\textbf{Z}^2$, satisfies:

\begin{itemize}
    \item[$\bullet$] \textit{Normalization for the derivative about the first degree of freedom}. The Radon-Nikodym derivative,

\begin{align*}
\mathrm{d}  \textbf{P}^{*}_{E_1}    := \textit{First Radon-Nikodym derivative taken over probability measures $\textbf{P}^{*} \big[ \cdot \big]$ of the triangular ef-} \\ \textit{fective random walk},
\end{align*}

    \noindent about the first degree of freedom $E_1$ is normalization by the summation of Radon-Nikodym derivatives,

\begin{align*}
{\mathrm{d}   \textbf{P}_{E_1} + \mathrm{d} \textbf{P}_{E_2}  }  : = \textit{Summation of Radon-Nikodym derivatives over probability measures $\textbf{P} \big[ \cdot \big]$ of tri-} \\ \textit{angular random walks} . 
\end{align*}

\bigskip

     \item[$\bullet$] \textit{Normalization for the derivative about the second degree of freedom}. The Radon-Nikodym derivative,

\begin{align*}
\mathrm{d}  \textbf{P}^{*}_{E_2} := \textit{Second Radon-Nikodym derivative taken over probability measures $\textbf{P}^{*} \big[ \cdot \big]$ of the triangular ef-} \\ \textit{fective random walk}  ,
\end{align*}

    \noindent about the first degree of freedom is normalization by the summation of Radon-Nikodym derivatives,

\begin{align*}  
     {\mathrm{d}   \textbf{P}_{E_1} + \mathrm{d} \textbf{P}_{E_2}  }   .
\end{align*}

\end{itemize}

\noindent To quantify how the above system of Radon-Nikodym derivatives is used for analyzing the triangular Green's function,

\begin{align*}
G_{\textbf{T}} \big( \lambda \big) , 
\end{align*}

\noindent for some real $\lambda$, we extend arguments provided in {\color{blue}[4]} for obtaining the scaling limit of the uniform kinetic prudent walk to the scaling limit of the triangular kinetic prudent walk. Over $\textbf{T}$, the stochastic process associated with the triangular random walk, along with its effective representation, satisfy:

\begin{itemize}
    \item[$\bullet$] \textit{Decomposition about each degree of freedom}. The representation of the triangular random walk over $\textbf{T}$ implies that the desired Green's function satisfies,

    \begin{align*}
      G_{\textbf{T}} \big( \lambda \big)  \propto \textbf{1}_{ \{ R = \textit{first degree of freedom} \}} \# \big\{ \textit{steps} : \textit{steps take along $R$}    \}    + \textbf{1}_{ \{ R = \textit{second degree of freedom} \}}  \# \big\{ \textit{steps} \\ : \textit{steps take along $R$}    \}   ,
    \end{align*}

\item[$\bullet$] \textit{Restriction of probability measures to slab measures appearing in the approximations for the scaling limit constant}. $\mathscr{L}$ is obtained through estimates related to the quantity,

\begin{align*}
    \textit{Slab measures over $\textbf{T}$} \propto \mathrm{exp} \big[ \textit{First Radon-Nikodym derivative} \big]  + \mathrm{exp} \big[ \textit{Second Radon-Nikodym} \\ \textit{ derivative} \big]          . 
\end{align*}

\item[$\bullet$] \textit{Sending the scaling limit normalization to $+\infty$}. As $\mathscr{L} \longrightarrow + \infty$,

{\small \begin{align*}
    \textit{Triangular effective expectation of }    \bigg\{  \begin{bmatrix}
\textit{Martingale along the first degree of freedom at time $t>0 $}  \\ 
    \textit{Martingale along the second degree of freedom at time $t^{\prime} > 0 $} \end{bmatrix} \bigg\}_{\{ t \neq t^{\prime} > 0  \} } \\ = \begin{bmatrix}
    1 \\ 1 \end{bmatrix} ,   
\end{align*} }

\noindent for the expectation under the triangular effective random walk. This expectation is introduced in the next section with $\mathscr{E} \big[ \cdot \big]$.

\bigskip

\item[$\bullet$] \textit{Formulating ratios of the restricted slab measures}. One manipulates objects of the form,

\begin{align*}
   \bigg\{ \frac{\textit{Slab measure $i$}}{\textit{Slab measure $i + 1$}}     \bigg\}_{i \in \textbf{N}}  , 
\end{align*}

\noindent to obtain the desired scaling limit.

\end{itemize}

\bigskip

\noindent While the degree of every vertex of $\textbf{Z}^2$ is $4$ and the degree of every vertex of $\textbf{T}$, the triangular lattice,

\begin{align*}
   \textbf{T} : = \big( V \big( \textbf{T} \big) , E \big( \textbf{T} \big) \big)  , 
\end{align*}

\noindent is $3$, nontrivial differences are expected to emerge between the scaling limits of the ordinary prudent walk, versus that of the triangular prudent walk, when applying computations similar to those provided by Beffara, Frieldi, and Velenik. To begin, denote the spanning set,

\begin{align*}
 \textit{Set of all linear combinations of basis vectors of the triangular lattice} := \mathrm{span}_{\textbf{T}} \big\{ \vec{\textbf{e}_1} , \vec{\textbf{e}_2} , \vec{\textbf{e}_3} \big\}  \\ = \mathrm{span}_{\textbf{T}} \big\{ {\textbf{e}_1} , {\textbf{e}_2} , {\textbf{e}_3} \big\}  , 
\end{align*}

\noindent of basis vectors,

\begin{align*}
       \textbf{e}_1 :=           \big[ \frac{\sqrt{3}}{2} , \frac{1}{2} , 0 \big]   , \\ \\  \textbf{e}_2 :=            \big[ \frac{\sqrt{3}}{2} , - \frac{1}{2} , 0 \big] , \\  \\    \textbf{e}_3 :=          \big[ 0 , 0 , 1 \big]            , 
\end{align*}

\noindent of the triangular lattice, $\textbf{T}$. As will be described, a random walk over $\textbf{T}$ can take steps in any of the following possible directions,

\[ \left\{\!\begin{array}{ll@{}>{{}}l} 
 \textit{Direction along } \textbf{e}_1 : =  \big[ \frac{\sqrt{3}}{2} , \frac{1}{2} , 0 \big]   , \\ \\   \textit{Direction along - degree of freedom of first coordinate of } \textbf{e}_1 : =  \big[ - \frac{\sqrt{3}}{2} , \frac{1}{2} , 0 \big]  ,  \\ \\  \textit{Direction along - degree of freedom of second coordinate of } \textbf{e}_1 : =  \big[  \frac{\sqrt{3}}{2} , - \frac{1}{2} , 0 \big]  ,  \\ \\  \textit{Direction along } - \textbf{e}_1 : =  \big[ -  \frac{\sqrt{3}}{2} , - \frac{1}{2} , 0 \big]  ,  \\  \\ \textit{Direction along } \textbf{e}_2 : =  \big[   \frac{\sqrt{3}}{2} , - \frac{1}{2} , 0 \big] , \\ \\ \textit{Direction along - degree of freedom of first coordinate of } \textbf{e}_2 : =  \big[  -  \frac{\sqrt{3}}{2} , - \frac{1}{2} , 0 \big] , \\    \\ \textit{Direction along - degree of freedom of second coordinate of } \textbf{e}_2 : =  \big[    \frac{\sqrt{3}}{2} ,  \frac{1}{2} , 0 \big] ,  \\   \\   \textit{Direction along } - \textbf{e}_2 : =  \big[   - \frac{\sqrt{3}}{2} ,  \frac{1}{2} , 0 \big] ,   \\  \\  \textit{Direction along } \textbf{e}_3 : =   \big[ 0 , 0 , 1 \big]  , \\ \\  \textit{Direction along } - \textbf{e}_3 : =  \big[ 0 , 0 , -1 \big] .
\end{array}\right. 
\]

\noindent Besides the definition of the basis and spanning set of $\textbf{T}$, to discuss expected differences in scaling limits, for the origin,

\begin{align*}
 \textit{Origin of the triangular lattice} : =  \big( 0 , 0, 0 \big) \in V \big( \textbf{T} \big)      ,
\end{align*}

\noindent introduce the sample space,

\begin{align*}
    \Omega_{\mathrm{Triangular}} := \bigg\{  \textit{paths} \text{ } \big( \pi_j \big)_{0 \leq j \leq L}  : \bigg\{   \pi_0 : = \textit{location of path at time 0}  = \big( 0 , 0, 0 \big) \bigg\}  ,  \bigg\{     \pi_{j+1} - \pi_j \in \big\{ e_1,  e_2 \\ , e_3 \big\} ,    \forall j \in \big[ 0 , L-1 \big] \bigg\} , \bigg\{   \big(     \pi_j + \textbf{N} \big( \pi_{j+1} - \pi_j  \big)      \big)          \cap   \big(   \pi_{[0,j]}          \big)  = \emptyset    ,    \forall j \in \big[ 0 , L-1 \big]     \bigg\}                     \bigg\} , 
\end{align*}

\noindent satisfying,

\begin{align*}
    \big\{ \Omega_{\mathrm{Triangular}}  \subsetneq \textbf{T} \big\}  , 
\end{align*}

\noindent of the triangular prudent walk. The definition of the above sample space resembles that provided for the uniform prudent walk, {\color{blue}[4]}; in words, the three conditions of the sample space for the triangular prudent walk above stipulate: (1) the restriction of a path of the triangular prudent walk, at time $0$, begins at the origin; (2) increments of the triangular prudent walk, $\pi_{j+1}-\pi_j$, can belong to the direction aligned with any one of the three possible basis vectors $e_1, e_2, e_3$; (3) the path of the triangular prudent walk, when restricted to $[0,j]$, and the path $\pi_j + \textbf{N}\big( \pi_{j+1} - \pi_j \big)$, are disjoint. Over such a collection of possible prudent paths, one can impose a probability measure,

\begin{align*}
     \textit{Uniform probability measure over triangular prudent paths } \pi    = \textbf{P}_{\mathrm{Unif},\textbf{T}} \big[ \pi \in \Omega_{\mathrm{Triangular}} \big]  \\   := \textbf{P}_{\mathscr{U},\textbf{T}} \big[ \pi \in \Omega_{\mathrm{Triangular}} \big]    = \big| \Omega_{\mathrm{Triangular}} \big|^{-1}      , 
\end{align*}

\noindent corresponding to each possible path $\pi$ being sampled with probability $\big| \Omega_{\mathrm{Triangular}} \big|^{-1} >0$.

\bigskip

\noindent In comparison to the \textit{kinetic} probability measure that will be introduced in \textit{1.3}, the \textit{uniform} probability measure over $ \Omega_{\mathrm{Triangular}} $ is used for concluding that:

\begin{itemize}
    \item[$\bullet$] \textit{As the triangular scaling limit normalization approaches positive infinity, the conditional probability that the excursion, after $\mathrm{log} \big[ \mathscr{L}^2 \big]$ many excursions, crosses the maximally displaced distance along the first degree of freedom occurs with probability zero}. One has that,

\begin{align*}
      \frac{1}{\sqrt{3} }    \underset{\mathscr{L} \longrightarrow + \infty}{\mathrm{lim}}                  \textbf{P}_{\mathscr{U},\textbf{T}} \big[   \exists i \in \big\{ \mathrm{log} \big[ \mathscr{L}^2 \big], \cdots ,  \gamma_{e_1} \big\}    :  \textit{The triangular walk crosses the maximally displaced} \\ \textit{ distance about $e_1$, the first degree of freedom, for every i}  \big]      ,
\end{align*}

     \item[$\bullet$] \textit{As the triangular scaling limit normalization approaches positive infinity, the conditional probability that the excursion, after $\mathrm{log} \big[ \mathscr{L}^2 \big]$ many excursions, crosses the maximally displaced distance along the second degree of freedom occurs with probability zero}. One has that,

     \begin{align*}
      \frac{1}{\sqrt{3} }    \underset{\mathscr{L} \longrightarrow + \infty}{\mathrm{lim}}         \textbf{P}_{\mathscr{U},\textbf{T}} \big[   \exists i \in \big\{ \mathrm{log} \big[ \mathscr{L}^2 \big] , \cdots ,  \gamma_{e_2} \big\}     :    \textit{The triangular walk crosses the maximally displaced} \\ \textit{ distance about $e_2$, the second degree of freedom, for every i}  \big]     ,
\end{align*}

      \item[$\bullet$] \textit{As the triangular scaling limit normalization approaches positive infinity, the conditional probability that the excursion, after $\mathrm{log} \big[ \mathscr{L}^2 \big]$ many excursions, crosses the maximally displaced distance along the third degree of freedom occurs with probability zero}. One has that,

      \begin{align*}
      \frac{1}{\sqrt{3} }    \underset{\mathscr{L} \longrightarrow + \infty}{\mathrm{lim}}          \textbf{P}_{\mathscr{U},\textbf{T}} \big[    \exists i \in \big\{ \mathrm{log} \big[ \mathscr{L}^2 \big] , \cdots ,  \gamma_{e_3} \big\}     :   \textit{The triangular walk crosses the maximally displaced} \\ \textit{ distance about $e_3$, the third degree of freedom, for every i}  \big]   ,
\end{align*}

\noindent for,

  \begin{align*}
        \gamma_{e_1} : = \textit{Path along the first degree of freedom of the triangular effective random walk} , \\ \\   \gamma_{e_2} : = \textit{Path along the second degree of freedom of the triangular effective random walk} , \\ \\   \gamma_{e_3} : = \textit{Path along the third degree of freedom of the triangular effective random walk} .
    \end{align*}

\end{itemize}

\noindent The above collection of results differs from that for the uniform kinetic prudent walk, whose stationary behavior is captured through the following three results:

\bigskip

\noindent Fix $\pi \in \Omega_{\textbf{Z}^2}$, where,

{\small \begin{align*}
   \Omega_{\textbf{Z}^2} : = \textit{sample space of the uniform, prudent, effective random walker over } \textbf{Z}^2 \\ =     \underset{\textit{path}}{\bigcup}  \bigg\{ \textit{path over } \textbf{Z}^2 :  \big\{ \textit{path lies strictly above the y, and x, axes of }  \textbf{Z}^2 \big\}  ,    \big\{   \textit{path } \\ \textit{  avoids any step that it has previously taken}          \big\}         \bigg\}  .           
\end{align*} } 

\noindent Denote,

{\small \begin{align*}
   \mathcal{E}_i \big( \pi \big) : =     \textbf{1} \bigg[ \textit{the uniform, prudent, effective random path } \pi \textit{ crosses a distance } \lambda > 0 \textit{ along either } R_1 \textit{, or } R_2    \bigg]        . 
\end{align*} }

\begin{itemize}

    \item[$\bullet$] \textbf{Lemma} (\textbf{Lemma} \textit{5.1}, {\color{blue}[4]}, \textit{under the square uniform probability measure, an effective path deviates more than the range} $\lambda$ \textit{with probability zero}). There exists $\delta > 0 $ such that,
    
    \begin{align*}
     \underset{L \longrightarrow + \infty}{\mathrm{lim}} \textbf{P}_{\textit{Unif}, L} \big[ \big\{  \exists i \in \big\{ \delta \mathrm{log} L , \cdots , \gamma_L \big( \pi \big) \big\} \big\}  : \mathcal{E}_i \big( \pi \big) = 1 \big]  = 0  .
    \end{align*}

     \item[$\bullet$]             \textbf{Lemma} (\textbf{Lemma} \textit{5.2}, {\color{blue}[4]}, \textit{under the square uniform probability measure, an effective path does not spend more time on order of} $\big(\mathrm{log} L \big)^2$, \textit{for} $\mathrm{log} L$ \textit{many excursions}). For every $\delta > 0$, there exists $\kappa > 0$ such that,

    \begin{align*}
     \underset{L \longrightarrow + \infty}{\mathrm{lim}} \textbf{P}_{\textit{Unif}, L} \big[ T_1 \big( \pi \big) + \cdots + T_{\delta \mathrm{log} L} \big( \pi \big) \geq \kappa \big( \mathrm{log} L \big)^2 \big] = 0    .
    \end{align*}

      \item[$\bullet$] \textbf{Lemma} (\textbf{Lemma} \textit{5.3}, {\color{blue}[4]}, \textit{under the square uniform probability measure, the length of an effective path, subtracted from the time spent at a constant of L many excursions, does not exceed } $\mathrm{log} L$). There exists $\alpha >0$ such that,

      \begin{align*}
       \underset{L \longrightarrow + \infty}{\mathrm{lim}} \textbf{P}_{\textit{Unif}, L} \big[ L - \big( T_1 + \cdots + T_{\gamma L} \big) \geq \alpha \mathrm {log} L    \big] = 0   .
      \end{align*}

\end{itemize}

\noindent In previous computations for the scaling limit of the square, uniform, kinetic prudent walk, a partition, {\color{blue}[4]},

\begin{align*}
  \Omega^+_L  := \textit{Restriction of the sample space of the kinetic, uniform prudent walk to paths} \\ \textit{which move in the northeast direction} ,
\end{align*}

\noindent was introduced over the sample space of uniform prudent walks, $\Omega_{\text{Uniform Prudent}}$ for some $L > 0$. The scaling limit for this variant of the prudent walk supported over $\textbf{Z}^2$ is obtained with the interpolation,

\begin{align*}
  \widetilde{\pi}^L_t : =  \frac{1}{L} \big[ \pi_{\lfloor t L \rfloor} +  \big( t L - \lfloor t L \rfloor \big) \big( \pi_{\lfloor t L \rfloor + 1} -    \pi_{\lfloor t L \rfloor }  \big) \big]      .
\end{align*}

\noindent The normalization, through $L$, demonstrates how one must consider, when computing a scaling limit:

\begin{itemize}

            \item[$\bullet$] The global maxima of Brownian motion paths over the underlying lattice,

             \item[$\bullet$] Truncations of random walks,

             \item[$\bullet$] Whether it is appropriate to rescale collections of path by space, by time, or by space and time simultaneously,

                     \item[$\bullet$] Conditionally defined probabilities involving time, or length, rescaling,

    \item[$\bullet$] Effective random walk representations,

        \item[$\bullet$] Scalings, with respect to length, or with respect to time,

                \item[$\bullet$] Restrictions to smaller sample spaces of prudent walks in the northeast direction.
\end{itemize}

\subsection{The main result}

For obtaining scaling limits of triangular prudent walkers, one expects a different normalization from scaling limits that have previously been obtained. For a path $\pi$ distributed according to the law of the triangular prudent walk, denote,

\begin{align*}
    \widetilde{\pi} : = \frac{1}{\mathscr{L}} \bigg\{                \big[    \pi_{\lfloor R t \rfloor }     +    \big[      t R - \lfloor t R \rfloor  \big] \big[   \pi_{\lfloor t R \rfloor + 1} - \pi_{\lfloor t R \rfloor}    \big]    \big]                 \textbf{1}_{\{ R = R_1 \} }      + \big[    \pi_{\lfloor R t \rfloor }      +      \big[      t R - \lfloor t R \rfloor  \big] \big[   \pi_{\lfloor t R \rfloor + 1} \\ - \pi_{\lfloor t R \rfloor}    \big]    \big] \textbf{1}_{\{ R = R_2 \} }                   \bigg\}    ,
\end{align*}

\noindent corresponding to the rescaled path.

\bigskip

\noindent Fix $c_{ \{ 1,e \} } \neq c_{ \{ 2,e\}  }>0$, and $k>0$. Given to triangular kinetic prudent paths $\gamma_s$ and $\hat{\gamma_s}$, while it is expected that there can exist paths $\widetilde{\pi_s} \equiv \widetilde{\pi_s} \big( \textbf{T} \big) $ and $\hat{\pi_s} \equiv \hat{\pi_s \big( \textbf{T} \big)}$, for which,

{\small \begin{align*}
   \bigg\{ \textbf{P}_{\mathscr{U},\textbf{T}}  \bigg[ \bigg\{   \text{ }     \underset{0 \leq s \leq t}{\mathrm{sup}} \bigg| \bigg|         \frac{1}{t} \big[    \hat{\pi_s} - \pi_s    \big]        \bigg|\bigg|_2 \geq \epsilon \bigg\} \text{ } \bigg| \text{ } \mathcal{Q}_1 \text{ }   \bigg]   \approx 1             \bigg\} \Longrightarrow      \bigg\{           \textbf{P}_{\mathscr{U},\textbf{T}}  \big[   \big|       \widetilde{\pi} t     -     \big[   c_{\{ 1,e \} } \textbf{1}_{    \{   e : e \in \mathrm{span} \{ \textit{first degree of freedom}        \}   \}   }  \\   +  c_{\{ 2,e \} } \textbf{1}_{    \{   e : e \in \mathrm{span} \{ \textit{second degree of freedom}        \}   \}   }     \big] t  \big|        \big]                                                                         \bigg\}        \\ \Downarrow \\     \bigg\{               \frac{1}{\sqrt{\mathscr{L}}}    \bigg\{         X_{\lfloor \mathscr{L} t \rfloor }                                       \bigg\}_{\{ t \geq 0\}  }                         \overset{t \longrightarrow + \infty}{\longrightarrow}   \big( \textit{Covariance matrix} \big)^{\frac{1}{2}}  \mathcal{B}^{(a)}  \bigg\}              ,  \tag{*}
\end{align*} }

\noindent where,

{\small \begin{align*}
 \textit{Covariance matrix} : =  \underset{k \longrightarrow + \infty}{\mathrm{lim}} \bigg\{ \textit{Expected value of $\bigg\{ \big\{    X_k   \big\}_{\{ k \geq 0\}} \big\{    X_{k^{\prime}}   \big\}_{\{ k^{\prime} \geq 0\}} \bigg\} $ }, \textit{conditionally upon the filtration } \\ \textit{up to step $k-1$} \bigg\}  , 
\end{align*} }

\noindent given $\epsilon$ sufficiently small, and,

\begin{align*}
    \mathcal{Q}_1 : =     \underset{\textit{paths } \gamma}{\bigcup} \big\{ \gamma : \gamma \text{ } \textit{is in the first quadrant of the triangular lattice at step n} \big\}        , \\ \\ \textit{filtration up to step $k-1$} : = \sigma \big( X_0 , \cdots , X_k \big) , 
\end{align*}

\noindent is upper bounded by some constant approaching $0$, one must introduce the random variables,

\begin{align*}
      \mathscr{R}\mathscr{V}_1 : =    \big\{ \textit{excursions of Brownian motion over the triangular lattice along the first degree of freedom} \big\}        , \\ \\       \mathscr{R}\mathscr{V}_2 : =  \big\{ \textit{excursions of Brownian motion over the triangular lattice along the second degree of freedom} \big\}      , \\ \\       \mathscr{R}\mathscr{V}_3 : =      \big\{ \textit{excursions of Brownian motion over the triangular lattice along the third degree of freedom} \big\}         , 
\end{align*}

\noindent associated with Brownian motion, $\mathscr{B}$, over $\textbf{T}$.

\bigskip

\noindent For the $\big| \big| \cdot \big| \big|_2$ norm,

\begin{align*}
  \big| \big|  \cdot \big| \big|_2 : =    \sqrt{\underset{1 \leq i \leq n}{\sum} \big| \cdot_i \big|^2  }   , 
\end{align*}

\noindent given a path $\gamma$ over $\textbf{T}$ that is used for defining the above difference,

\begin{align*}
   \bigg| \bigg| \frac{1}{t} \big( \hat{\pi_s} - \pi_s \big)  \bigg| \bigg|_2  ,
\end{align*}

\noindent with respect to the $l-2$ norm,if the supremum of the $\mathrm{l}$-2 norm between $\hat{\gamma_s}$ and $\gamma_s$ exceeds $\epsilon$, then there exists a coupling between $\hat{\gamma_s}$ and $\gamma_s$ for which:

\begin{itemize}
\item[$\bullet$] \textit{Property 1}: There exists a mapping between $\hat{\gamma_s}$ and $\gamma_s$ for which,

\begin{align*}
 \phi : \hat{\pi_s} \longrightarrow \pi_s :  \text{ } \text{there exists countably many } k \text{ } \text{for which }   \big\{    \hat{\pi_s} \neq \pi_s       \big\}    . 
\end{align*}

\item[$\bullet$] \textit{Property 2}: Dependent upon the occurrence of the event $\mathcal{Q}_1$,

\begin{align*}
  \frac{\textbf{P} \bigg[ \big\{  \underset{0 \leq s \leq t}{\mathrm{sup}} \big| \big|         \frac{1}{t} \big(    \hat{\pi_s} - \pi_s    \big)       \big|\big|_2 \geq \epsilon        \big\} \cap \big\{ \mathcal{Q}_1  \big\}        \bigg] }{\textbf{P} \big(   \mathcal{Q}_1       \big)} \leq \epsilon   . 
\end{align*}

\item[$\bullet$] \textit{Property 3}. Almost surely, there exists a coupling between paths for which,

\begin{align*}
 \underset{t \geq 0}{\mathrm{sup}}   \big| \big| \hat{\gamma_t} - \gamma_t \big| \big|_1 < + \infty. 
\end{align*}

\end{itemize}

\noindent Each $\mathscr{R}\mathscr{V}$ above is used to compute several probabilistic quantities, particularly the ones given below,

\begin{align*}
 \textbf{P}_{\mathscr{U},\textbf{T}} \big[  \mathscr{A}^{\mathrm{v}}_k                   \big|    \mathcal{H}_{\mathscr{T}_k} \geq c_1 k^{\frac{4}{3}}     \big]  , \\  \\  \textbf{P}_{\mathscr{U},\textbf{T}} \big[       \mathcal{W}_{\mathscr{T}_k}  < c_2 k^{\frac{4}{3}}  \big]  , \\ \\ \textbf{P}_{\mathscr{U},\textbf{T}} \big[       \mathcal{H}_{\mathscr{T}_k} < c_2 k^{\frac{4}{3}}  \big] , 
 \end{align*}

\noindent corresponding to the probability that a \textit{vertical} excursion, $\mathscr{A}^{\mathrm{v}}$ occurs, conditionally upon $\big\{   \mathcal{H}_{\mathscr{T}_k} \geq c_1 k^{\frac{4}{3}}  \big\}$, the probability that, at time $\mathscr{T}_k$, the vertical increment of the prudent walk, at time $\mathscr{T}_k$, is strictly smaller than $k^{\frac{4}{3}}$, up to a constant, and lastly, the probability that the horizontal increment of the prudent walk, at time $\mathscr{T}_k$, is strictly smaller than $k^{\frac{4}{3}}$, up to a constant. Related to the conditional probability provided in $(*)$, the main result that will be shown to hold for computing the desired scaling limit claims:

\bigskip

\noindent \textbf{Theorem} (\textit{the probability of the supremum exceeding some sufficiently small parameter, conditional upon the prudent walker settling in the first quadrant of the triangular lattice, is zero}). Denote $\mathcal{Q}_1$ as the event that the prudent walker's last step of an excursion is in the first quadrant of the triangular lattice. For $\epsilon$ sufficiently small, one has that $\textit{(*)}$ holds.

\bigskip

\noindent While one can expect that probabilistic arguments, as given in the conditions provided above which parallel those initially provided by Beffara, Friedli, and Velenik, {\color{blue}[1]}, would apply to the triangular prudent walk, one must introduce different scaling factors rather than $k^{\frac{4}{3}}$, up to a constant. For now, we denote the power appearing in the same constant $k$, with some strictly positive constant, $\beta$. As the prudent walker over $\textbf{T}$ takes an infinite number of steps, $\beta$ is expected to dictate the scale along which the probabilities,

\begin{tabular}{|l|l|}
\hline\parbox[t]{0.25\textwidth}{
\begin{itemize}
\item \textbf{Definition} \textit{1}
\item \textbf{Definition} \textit{2}
\item \textbf{Definition} \textit{3}
\item \textbf{Definition} \textit{4}
\item \textbf{Definition} \textit{5}
\item \textbf{Definition} \textit{6}
\item \textbf{Definition} \textit{7}
\item \textbf{Definition} \textit{8}
\item \textbf{Definition} \textit{9}
\item \textbf{Definition} \textit{10}
\item \textbf{Definition} \textit{11}
\item \textbf{Definition} \textit{12}
\item \textbf{Definition} \textit{13}
\item \textbf{Definition} \textit{14}
\item \textbf{Definition} \textit{15}
\item \textbf{Definition} \textit{16}
\item \textbf{Definition} \textit{17}
\item \textbf{Definition} \textit{18}
\item \textbf{Definition} \textit{19}
\item \textbf{Definition} \textit{20}
\item \textbf{Definition} \textit{21}
\item \textbf{Definition} \textit{22}
\item \textbf{Definition} \textit{23}
\item \textbf{Definition} \textit{24}
\item \textbf{Definition} \textit{25}
\end{itemize}}& 
\parbox[t]{0.68\textwidth}{
\begin{itemize}
\item \textit{the probability measure over square uniform prudent paths}
\item \textit{the uniform probability measure over prudent walks}
\item \textit{nonnegative excursions of the square effective prudent walk}
\item \textit{Radon-Nikodym derivative over} $\textbf{Z}^2$
\item  \textit{the Green's function, in terms of the expected number of visits}
\item \textit{the probability measure over triangular uniform prudent paths}
\item \textit{the uniform probability measure over square effective prudent walks}
\item \textit{nonnegative excursions of the square effective prudent walk}
\item \textit{Radon-Nikodym derivative over} $\textbf{T}$
\item \textit{the triangular Green's function}
\item \textit{sample space of the effective triangular random walk}
\item \textit{excursions sampled over the collection of random paths of length L}
\item \textit{the total length, and number of increments, of random paths}
\item \textit{the total length, and number of increments, of triangular effective random paths of length L}
\item \textit{truncation of the effective random walk over } $\textbf{T}$
\item \textit{ordinary, and truncated, sample spaces of the triangular effective walk over} $\textbf{T}$
\item \textit{truncation of the effective random walk over } $\textbf{T}$
\item \textit{inductive definition of the length, and number of excursion}
\item \textit{Slabs over } $\textbf{Z}^2$
\item \textit{stopping times}
\item \textit{Radon-Nikodym derivative of the triangular effective random walk}
\item \textit{excursions of the triangular effective random walk about the second degree of freedom}
\item \textit{slab measures over} $\textbf{N} \times \textbf{N} \times\textbf{N} \times \big[ 0 , 1 \big]$
\item \textit{slab measures over} $\textbf{N} \times \textbf{N} \times\textbf{N} \times \big[ 0 , 1 \big]$ \textit{for incomplete excursions}
\item \textit{the set of triangular effective paths within a slab that end at 0, or at R}
\end{itemize}}\\
\hline
\end{tabular}
\noindent \textit{Table 1}. An overview of each \textbf{Definition} provided in the next section.

\bigskip

\begin{align*}
   \textbf{P}_{\mathscr{U},\textbf{T}} \big[ \textit{k th step of an excursion along the direction of the basis vector } e_1 \big| \textit{k th} \\ \textit{    step of an excursion along the direction of the basis vector } e_2 \gtrsim k^{\beta} \big]    , \\ \\   \textbf{P}_{\mathscr{U},\textbf{T}} \big[ \textit{k th step of an excursion along the direction of the basis vector } e_2 \big| \textit{k th} \end{align*}

\begin{tabular}{|l|l|}
\hline\parbox[t]{0.25\textwidth}{
\begin{itemize}
\item \textbf{Definition} \textit{26}
\item \textbf{Definition} \textit{27}
\item \textbf{Definition} \textit{28} 
\item \textbf{Definition} \textit{29}
\item \textbf{Definition} \textit{30}
\item \textbf{Definition} \textit{31}
\item \textbf{Definition} \textit{32}
\item \textbf{Definition} \textit{33}
\item \textbf{Definition} \textit{34}
\item \textbf{Definition} \textit{35}
\item
\textbf{Definition} \textit{36}
\end{itemize}}& 
\parbox[t]{0.68\textwidth}{
\begin{itemize}
\item \textit{excursions about each degree of freedom}
\item \textit{displacement of the triangular effective walk}
\item \textit{computation of the displacement of the triangular effective walk}
\item \textit{normalizing constant over the triangular lattice}
\item \textit{properties of the normalizing constant}
\item \textit{Indicator function for the triangular effective random walk}
\item \textit{Uniform, three-sided representation}
\item \textit{expected value of the slab, and renormalized, slab measures}
\item \textit{renormalized paths}
\item \textit{truncation procedure over} $\textbf{Z}^2$
\item \textit{stopping times for the square effective random walk}
\end{itemize}}\\
\hline
\end{tabular}
\noindent \textit{Table 2}. An overview of each remaining \textbf{Definition} provided in the next section.

\bigskip

   \begin{align*} \textit{    step of an excursion along the direction of the basis vector } e_3 \gtrsim k^{\beta} \big]   , \\ \\   \textbf{P}_{\mathscr{U},\textbf{T}} \big[ \textit{k th step of the excursion along the direction of } e_1  \lesssim k^{\beta} \big]   , \\ \\   \textbf{P}_{\mathscr{U},\textbf{T}} \big[ \textit{k th step of the excursion along the direction of } e_2  \lesssim  k^{\beta} \big]   , \\ \\   \textbf{P}_{\mathscr{U},\textbf{T}} \big[  \textit{k th step of the excursion along the direction of } e_3 \lesssim k^{\beta}  \big]   ,   
\end{align*}

\noindent occur, under the identification,

\begin{align*}
     \textit{Prudent path at time step 1} = \textit{Prudent path at time 1} : =    \pi_{[0,1]} = \pi_0  \overset{\cdot}{\cup} \pi_1   , \\ \vdots \\  \textit{Prudent path at time step L} = \textit{Prudent path at time L} : =    \pi_{[0,L]} = \pi_0  \overset{\cdot}{\cup} \pi_1   \overset{\cdot}{\cup} \cdots  \overset{\cdot}{\cup} \pi_L ,   \\ \\      \textit{Restriction of the Prudent path at time step 1 to the direction of } e_1  = \textit{Prudent path at} \end{align*}

     \begin{align*} \textit{ time 1} \bigg|_{e_1}  : =    \pi_{[0,1]} \bigg|_{e_1}   = \big[  \pi_0  \overset{\cdot}{\cup} \pi_1 \big] \bigg|_{e_1}   , \\   \vdots \\ \textit{Restriction of the Prudent path at time step L to the direction of } e_1 = \textit{Prudent path at time L} \bigg|_{e_1} \\  : =    \pi_{[0,L]} \bigg|_{e_1}  =  \big[ \pi_0  \overset{\cdot}{\cup} \pi_1   \overset{\cdot}{\cup} \cdots  \overset{\cdot}{\cup} \pi_L \big] \bigg|_{e_1}  ,   \\  \vdots \\  \textit{Restriction of the Prudent path at time step L to the direction of } e_3 = \textit{Prudent path at time L} \bigg|_{e_3} \end{align*}
     
    \begin{align*} : =    \pi_{[0,L]} \bigg|_{e_3}  =  \big[ \pi_0  \overset{\cdot}{\cup} \pi_1   \overset{\cdot}{\cup} \cdots  \overset{\cdot}{\cup} \pi_L \big] \bigg|_{e_3}  . 
\end{align*}

\noindent Before introducing objects associated with paths $\pi \sim \Omega_{\mathrm{Triangular}}$, denote,

\begin{align*}
 \textit{Vertices of the prudent path } \pi :=    V \big( \pi \big)  \subsetneq V  \big( \textbf{T} \big)    ,  \\ \\ \textit{Edges of the prudent path } \pi :=    E \big( \pi \big)  \subsetneq E \big( \textbf{T} \big)    ,
\end{align*}

\noindent corresponding to the edges of the triangular lattice that some path $\pi$ occupies,

\begin{align*}
  \textit{Number of steps taken by a triangular prudent walker} : =  \underset{1 \leq t \leq L}{\bigcup}   \textit{Number of steps taken by a triangular} \\ \textit{ prudent walker at time t} =         \underset{1 \leq t \leq L}{\bigcup}    \pi_t                    , 
\end{align*}

\noindent corresponding to the number of walks that the path $\pi_t$ takes, as a function of the time $t$,

\begin{align*}
    \textit{Sites of } \pi_t : =  \mathrm{Sites} \big( \pi_t  \big) = \underset{v. \in V ( \textbf{T})}{\bigcup}   \big\{ v :       v \cap V \big( \pi_t \big) \neq \emptyset \big\}     , 
\end{align*}

\noindent corresponding to the sites of $\pi$ that have been occupied up to some time $t$,  in addition to,

\begin{align*}
 \textit{Prudent walkers} : = \underset{\omega  \in \Omega_{\mathrm{Triangular}}}{\bigcup} \textit{Prudent walker} \big( \omega \big)           , 
\end{align*}

\noindent corresponding to set of all triangular prudent walker, with each respective configuration $\omega$ drawn uniformly at random from $\Omega_{\mathrm{Triangular}}$.

\bigskip

\noindent Under the natural decomposition,

\begin{align*}
  E \big( \pi \big)  := \underset{0 \leq j \leq L}{\bigcup} E \big( \pi_j \big)  ,
\end{align*}

\noindent it is convenient to introduce,

\begin{align*}
  \textit{Restriction of the prudent path } \pi \textit{ along the basis vector } e_1 :=  \pi \bigg|_{e_1} = \underset{0 \leq j \leq L}{\bigcup } \big\{   e \in E \big( \textbf{T} \big) :  e \cap E \big( \pi_j \big) \\  \neq 0  , < e , e_1 >   = 0      \big\}  ,  \\ \\ 
  \textit{Restriction of the prudent path } \pi \textit{ along the basis vector } e_2 :=  \pi \bigg|_{e_2} =   \underset{0 \leq j \leq L}{\bigcup } \big\{    e \in E \big( \textbf{T} \big) :  e \cap E \big( \pi_j \big) \\ \neq 0  , < e , e_2 >  = 0           \big\}   , \\ \\ 
  \textit{Restriction of the prudent path } \pi \textit{ along the basis vector } e_3 := \pi \bigg|_{e_3} =    \underset{0 \leq j \leq L}{\bigcup } \big\{    e \in E \big( \textbf{T} \big) :  e \cap E \big( \pi_j \big) \\  \neq 0  , < e , e_3 >  = 0           \big\}   , 
\end{align*}

\noindent from which one can introduce the corresponding times,

\begin{align*}
  \textit{Time that the triangular prudent walker is at a distance of 1 from the maximum achieved in the k th} \\ \textit{step of the excursion in the direction of } e_1 := \mathcal{T}_{1}   , \\ \\   \textit{Time that the triangular prudent walker is at a distance of 1 from the maximum achieved in the k th} \\ \textit{step of the excursion in the direction of } e_2 :=  \mathcal{T}_{2}  , \\ \\   \textit{Time that the triangular prudent walker is at a distance of 1 from the maximum achieved in the k th} \\ \textit{step of the excursion in the direction of } e_3  :=   \mathcal{T}_{3}  , 
\end{align*}

\noindent which respectively correspond to the steps of the excursion for which the a prudent path is at distance $1$ from the maximum along the basis vectors $e_1, e_2$, and $e_3$, respectively. 

\bigskip

\noindent As the prudent walker takes infinitely many steps, to quantitatively determine the \textit{asymptotic order of fluctuation} over $\textbf{T}$, one would compute the probabilities,

{\small \begin{align*}
  \underset{t \longrightarrow + \infty}{\mathrm{lim}}    \textbf{P}_{\mathscr{U},\textbf{T}} \bigg[ \underset{0 \leq s \leq t}{\mathrm{max}} \bigg\{   \mathrm{min} \bigg\{   n_1 : \big\{ < n_1, e_1 > = 0 \big\} , \big\{ \mathrm{diam}  \big[ \mathscr{E}_{s,1}  \big]  \geq \big( n_1 \big)^s       \big\}                 \bigg\}      \bigg\}  \bigg]   , \\ \\    \underset{t \longrightarrow + \infty}{\mathrm{lim}} \textbf{P}_{\mathscr{U},\textbf{T}} \bigg[ \underset{0 \leq s \leq t}{\mathrm{max}} \bigg\{   \mathrm{min} \bigg\{   n_1 : \big\{ < n_1, e_2 > = 0 \big\} , \big\{ \mathrm{diam}  \big[ \mathscr{E}_{s,2}  \big]  \geq \big( n_1 \big)^s       \big\}                 \bigg\}      \bigg\}  \bigg]      , \\ \\     \underset{t \longrightarrow + \infty}{\mathrm{lim}}  \textbf{P}_{\mathscr{U},\textbf{T}} \bigg[ \underset{0 \leq s \leq t}{\mathrm{max}} \bigg\{   \mathrm{min} \bigg\{   n_1 : \big\{ < n_1, e_3 > = 0 \big\} , \big\{ \mathrm{diam}  \big[ \mathscr{E}_{s,3}  \big]  \geq \big( n_1 \big)^s       \big\}                 \bigg\}      \bigg\}  \bigg]    , 
\end{align*} }

\noindent for the excursions,

\begin{align*}
 \mathscr{E}_{s,1}   \equiv \textit{excursion along first degree of freedom of the triangular lattice}  , \\ \\ \mathscr{E}_{s,2}  \equiv \textit{excursion along second degree of freedom of the triangular lattice}  , \\ \\  \mathscr{E}_{s,3}  \equiv \textit{excursion along third degree of freedom of the triangular lattice}   . 
\end{align*}

\noindent Given the times $\mathcal{T}_1$, $\mathcal{T}_2$, and $\mathcal{T}_3$, another component of the arguments due to Beffara, Friedli, and Velenik, as has been adopted for computing the scaling limit of the uniform prudent walk, {\color{blue}[4]}, remain applicable for the triangular prudent walk.

\noindent Fix times $\mathcal{T}_1, \mathcal{T}_2, \mathcal{T}_3$, satisfying,

\begin{align*}
0 <   \mathcal{T}_1 < \mathcal{T}_2 < \mathcal{T}_3   ,
\end{align*}

\noindent That is, one can, along the lines of the quantities,

\begin{align*}
\textit{Restriction of a triangular lattice path along the first degree of freedom} : =  \pi \bigg|_{e_1}  ,    \\ \\ \textit{Restriction of a triangular lattice path along the second degree of freedom} : =   \pi \bigg|_{e_2}  , \\ \\ \textit{Restriction of a triangular lattice path along the third degree of freedom} : =   \pi \bigg|_{e_3}  ,
\end{align*}

\noindent one can also introduce,

\begin{align*}
  \pi \bigg|_{e_1} \textit{between times 0 and } \mathcal{T}_1 - 1   : =     \pi \bigg|_{e_1, [0, \mathcal{T}_1 - 1 ]}  ,    \\ \\  \pi \bigg|_{e_2} \textit{between times } \mathcal{T}_1 + 1 \textit{ and } \mathcal{T}_2 - 1  : = \pi \bigg|_{e_2, [ \mathcal{T}_1 +1 , \mathcal{T}_2 - 1 ]} , \\ \\   \pi \bigg|_{e_3} \textit{between times } \mathcal{T}_2 +1 \textit{ and } \mathcal{T}_3 - 1 : =     \pi \bigg|_{e_3, [ \mathcal{T}_2 +1 , \mathcal{T}_3 - 1 ]}  ,
\end{align*}

\noindent corresponding to \textit{restrictions} of $\pi$, respectively along $e_1, e_2,$ and $e_3$, at times \textit{before} the walk achieves the largest distance from $\big( 0 , 0 , 0 \big)$. From such a collection of prudent paths, the triangular \textit{effective random walk} over $\textbf{T}$, takes the form,

{\small \begin{align*}
  \mathscr{E}\mathscr{W}_{\textbf{T}} : = \bigg\{ \omega \in \Omega_{\mathrm{Triangular} } : \bigg\{ \bigg[ \bigg|  \pi \bigg|_{e_1, [0, \mathcal{T}_1 - 1 ]}  \cap \pi \bigg|_{e_1} \bigg| \bigg]  =  \big| \# \big\{ i : \mathcal{T}_0 \leq \mathcal{T}_i \leq \mathcal{T}_1 - 1 \big\} \big|  \bigg\} ,   \bigg\{ \bigg[ \bigg|    \pi \bigg|_{e_2, [ \mathcal{T}_1 +1 , \mathcal{T}_2 - 1 ]}  \cap  \pi \bigg|_{e_2} \bigg| \bigg] \\  =        \big| \# \big\{ i : \mathcal{T}_1 + 1 \leq \mathcal{T}_i \leq \mathcal{T}_2 - 1  \big\} \big|    \bigg\} , \bigg\{ \bigg[ \bigg|   \pi \bigg|_{e_3, [ \mathcal{T}_2 +1 , \mathcal{T}_3 - 1 ]}   \cap     \pi \bigg|_{e_3} \bigg| \bigg]   =         \big| \# \big\{ i : \mathcal{T}_2 + 1 \leq \mathcal{T}_i \leq \mathcal{T}_3 - 1  \big\} \big|           \bigg\}          \bigg\}         . 
\end{align*} }

\noindent In this work, we not only seek to determine how key components of the argument due to Beffara, Friedli, and Velenik, can be applied from the ordinary prudent walk to the triangular prudent walk, but also how the \textit{scaling limit} of a random walk depends upon the underlying lattice.

\bigskip

In order to study the scaling limit of the triangular prudent walk, one must introduced the following probabilistic objects:

\begin{itemize}
    \item[$\bullet$] \textit{Times at which a triangular prudent path exits an interval of strictly positive Lebesgue measure}.

    \item[$\bullet$] \textit{Displacements along each basis element of the triangular lattice}.

    \item[$\bullet$] \textit{The normalization constant, which can depend upon time or space, for the scaling limit of the triangular kinetic prudent walk}.

        \item[$\bullet$] \textit{The rate of convergence of the normalizing constant}.

        \item[$\bullet$] \textit{The rate of convergence of the numerator of the scaling limit before normalization}.
\end{itemize}

Besides other closely related types of random walks to the prudent walk, including those discussed in {\color{blue}[2,3,4]}, we answer one question raised by Beffara, Friedli, and Velenik, raised in {\color{blue}[1]}, as to whether the computations for obtaining the scaling limit of the prudent walk over the square lattice can be extended to obtain the scaling limit of the prudent walk over the triangular lattice. In the case of the square lattice, the probability that some prudent walker intersection with some point $m$, conditionally upon the path explored up to some time $\mathscr{T}_k$, is equal to the probability that the prudent walk arrives to $m$ at $\mathscr{T}_k$.

\subsection{Paper organization}

\noindent From general characteristics of random walks discussed in \textit{1.1}, in \textit{2} below we describe how excursions of the triangular kinetic prudent walk can be partitioned. Associated with such partitions are measures over not only over excursion of the effective walk representation, but also over the number of horizontal excursion paths. In previous work on the square kinetic prudent walk, {\color{blue}[4]}, a computation of the Radon-Nikodym derivative was used to determine a $\frac{4}{3}$ factor appearing in the probability measure of effective random walks.

\bigskip

\noindent The representation of the probability measure  $  \textbf{P}^{*}_{\textit{Kinetic}, \textbf{T}} \big[ \cdot \big] $ depends upon: (1) the Radon-Nikodym computation over the shared sample space of kinetic triangular prudent paths, $\mathscr{V}_{\textbf{T}}$; (2) two measures, $\mathscr{L}$ and $\mathscr{L}^{\prime}$, respectively, supported over $\textbf{N} \times \textbf{N} \times \textbf{N} \times \big[ 0 , 1 \big]$ and a volume $\mathcal{V} \subsetneq \textbf{T}$; (3) a \textit{truncation} of the effective random walk; (4) a suitable counterpart for the Green's function, originally introduced through $G$ over $\textbf{Z}^2$. As a function that depends upon the number of vertices, and hence, edges, that the random walk traverses, one can make use of the Green's function for:

\begin{itemize}
    \item[$\bullet$] \textit{Obtaining estimates on a series of the Green's function after a strictly positive number of steps}. Fix $r>0$ as the number of steps of a random walk over $\textbf{T}$. For the Green's function supported over $\textbf{T}$, $G_{\textbf{T}} \big( \lambda \big)$, the series,

    \begin{align*}
      \underset{r > 0}{\sum} G^r_{\textbf{T}} \big( \lambda \big)     ,
    \end{align*}

    \noindent is strictly upper bounded by a suitable, strictly positive, constant, where the decomposition of the Green's function takes the form,

\begin{align*}
 G^r_{\textbf{T}} \big( \lambda \big)  =      G^r_{T \subsetneq \textbf{T}} \big( \lambda \big)     \textbf{1}_{\{ R = R_1 : R_1 \text{ } \textit{spans}\text{ }  T \} }           +  G^r_{T^{\prime} \subsetneq      \textbf{T}} \big( \lambda \big)     \textbf{1}_{\{ R = R_2 : R_2 \text{ } \textit{spans} \text{ } T^{\prime} \} }                      , 
\end{align*}

    \noindent for,

    \begin{align*}
   R_1 : =   \textit{First degree of freedom of the triangular effective walk over the slab}          ,  \\ \\ R_2 : =          \textit{Second degree of freedom of the triangular effective walk over the slab}   . 
\end{align*}

      \item[$\bullet$] \textit{Leveraging behavior of the Green's function}. The behavior of the Green's function provided in the previous item above depends upon strictly positive $\alpha^{\prime}$, for which,

        \begin{align*}
         \alpha :  = \mathrm{log} \bigg[ \mathscr{E} \big[ \mathrm{exp} \big[ - \lambda \big| \mathscr{U}_1 \big| \big] \big]  \bigg]     . 
        \end{align*}

        \item[$\bullet$] \textit{Relating the Radon-Nikodym computation to the parameter of the Green's function}. The parameter introduced in the previous item above in the definition of the Green's function is equal to the Radon-Nikodym derivative,

\begin{align*}
 \frac{\mathrm{d}\textbf{P}^{*}_{\textit{Kinetic}, \textbf{T}} \big[  \textit{Paths}_{\infty, \textbf{T}} \big]}{\mathrm{d} \textbf{P}_{\textit{Kinetic}, \textbf{T}} \big[  \textit{Paths}_{\infty, \textbf{T}} \big] }  : =  \frac{\mathrm{d}\textbf{P}^{*}_{\textit{Kinetic}, \textbf{T}} }{\mathrm{d} \textbf{P}_{\textit{Kinetic}, \textbf{T}}  }          \big[    \textit{Paths}_{\infty, \textbf{T}}    \big] = \mathrm{d} \frac{\textbf{P}^*}{\textbf{P}}  , 
\end{align*}

\noindent for the collection of paths,

\begin{align*}
    \textit{Paths}_{\infty, \textbf{T}} :=  \underset{N \geq 1}{\bigcup} \bigg\{  \forall i \leq N, \exists \big(  \textit{Paths} \big)_{0 \leq i \leq N} \subsetneq \textbf{T}  : \big\{ \textit{Paths}_0 =  0 \big\}  ,  \big\{ \textit{Paths}_i           \geq 0 \big\} , \big\{ \textit{Paths}_N = 0   \big\}      \bigg\} . 
\end{align*}

        \item[$\bullet$] \textit{Obtaining the triangular kinetic scaling limit}. Sending the number of steps of the triangular prudent walk to $+\infty$ implies that the desired scaling limit mentioned in \textit{1.1} holds, with respect to the probability measure,

   \begin{align*}
        \mathscr{P}^{*}_{\textbf{T}} \big[ \cdot \big] : =   \mathscr{P}^{*}    \big[ \cdot \big]  =   \textbf{P}_{\mathscr{U},\textbf{T}}  \big[   \textit{effective paths } P : \big\{ P  \textit{ lies strictly above } \textbf{T} \big\} , \big\{  P_0 = P_L = 0 \big\}     \big]            , 
        \end{align*}

\noindent described further in \textit{2.2}.

\end{itemize}

\section{Probabilistic setting}

\subsection{Square Kinetic Prudent walker objects}

\noindent We provide an overview of the objects associated with the square uniform kinetic walk, {\color{blue}[4]}.

\bigskip

\noindent \textbf{Definition} \textit{1} (\textit{the probability measure over the sample space of square uniform prudent paths}). Given a random walk $V \subsetneq \textbf{Z}^2$, $t \in \textbf{N}$, and the time,

\begin{align*}
    \eta_t =   \mathrm{min} \bigg\{ i \geq 1: i + \underset{1 \leq j \leq i}{\sum} \big| U_j \big| \geq t \bigg\}       , 
\end{align*}

\noindent corresponding to the first probability measure, $K$, define,

\begin{align*}
  K \big( t \big)  : =  K = \textbf{E} \bigg\{   \mathrm{exp} \bigg[ \mathrm{log} \bigg[ \frac{3}{2 } \bigg] \eta_t  \bigg]                            \textbf{1}_{ \{ V_i \geq 0, \forall i \leq \eta_t, V_{\eta_t} = 0 , \eta_t + \underset{1 \leq j \leq \eta_t}{\sum} | U_j | = t \}}         \bigg\}         . 
\end{align*}

\noindent A factorization of the first probability measure $K \big( t \big)$ above,

\begin{align*}
   \underset{1 \leq i \leq r}{\prod} K \big( t_i \big)  ,
\end{align*}

\noindent for some $r \leq L$ can be related to the sample space of square uniform prudent paths, $\Omega_L$, with,

{\small \begin{align*}
  2^{-L} \big| \Omega_L \big| = \bigg\{  \underset{1 \leq i \leq r}{\prod} 2^{-t_i} \bigg\} \big| \Omega_L \big|           =  \underset{r \geq 1}{\sum}  \bigg\{ \underset{t_1 + \cdots + t_r = L}{\sum} \bigg\{  \underset{1 \leq i \leq r}{\prod} 2^{-t_i} \big| \mathcal{I}_{t_i} \big| \bigg\} \bigg\}   = \underset{r \geq 1}{\sum}  \bigg\{ \underset{t_1 + \cdots + t_r = L}{\sum} \bigg\{  \underset{1 \leq i \leq r}{\prod} K \big( t_i \big)  \bigg\} \bigg\}                                                        , 
\end{align*} }

\noindent up to a $2^L$ normalization, where,

{\small \begin{align*}
      \mathcal{I}_{T} : = \underset{1 \leq i \leq t}{\bigcup} \mathcal{I}_i =   \bigg\{\forall i \in \big\{ 0 , \cdots , t \big\} = T ,  \exists \textit{paths } \pi \equiv \big[ \pi_0 , \cdots , \pi_t \big]  \subsetneq \textbf{Z}^2 : \big\{ \pi_0 = \big( 0 , 0 \big) , \pi_1 = \big( 1 , 0 \big) \big\} \\ , \big\{ \pi_{i,2 } \geq 0   , \pi_{t,2} = 0 \big\} \bigg\}           ,
\end{align*} }

\noindent and,

\begin{align*}
   \mathcal{I}_t \bigg|_{\textit{time }i} = \mathcal{I}_{t_i}  .
\end{align*}

\noindent \textbf{Definition} \textit{2} (\textit{the uniform probability measure over prudent walks from the probability measure over square effective prudent walks}). Denote the times $T_i$, where $1 \leq i \leq L$, for which the square uniform prudent walk occupies $i \in V \big( \textbf{Z}^2 \big)$. Corresponding to the second probability measure, $\mathscr{P}_{\textit{Uniform}}    \big[ \cdot \big] $, denote the probability measure over square effective random walks,

\begin{align*}
   \mathscr{P}^{*}_{\textbf{Z}^2} \big[ \cdot \big] : =   \mathscr{P}^{*}    \big[ \cdot \big]  =   \textbf{P}_{\mathscr{U},\textbf{Z}^2} \big[   \textit{effective paths } P : \big\{  P  \textit{ lies strictly above } \textbf{Z}^2 \big\} ,  \big\{ P_0 = P_L = 0 \big\}     \big]          , 
\end{align*}

\noindent from which one can define,

{\small \begin{align*}
\mathscr{P}_{\textit{Uniform}}    \big[ \cdot \big]  : = \mathscr{P}_{\textit{Uniform},\textbf{Z}^2}    \big[ \cdot \big]     =   \frac{\textbf{E} \bigg[  \textbf{1}_{\cdot}   \times  \textbf{1}_{\underset{1 \leq i \leq L}{\sum} \widetilde{T}_i = L  }                \bigg]}{\mathscr{P}^{*} \bigg[  \underset{1 \leq i \leq L}{\sum} \widetilde{T}_i = L       \bigg] }      , 
\end{align*} }

\noindent from the measure,

\begin{align*}
     \textit{Uniform probability measure over square prudent paths } \pi    = \textbf{P}_{\mathrm{Unif},\textbf{Z}^2} \big[ \pi \in \Omega_{\mathrm{Triangular}} \big]   \\  := \textbf{P}_{\mathscr{U},\textbf{Z}^2} \big[ \pi \in \Omega_{\mathrm{Square}} \big]    = \big| \Omega_{\mathrm{Square}} \big|^{-1}      , 
\end{align*}

\noindent and sample space over $\textbf{Z}^2$,

\begin{align*}
    \Omega_{\mathrm{Square}}  . 
\end{align*}

\noindent \textbf{Definition} \textit{3} (\textit{nonnegative excursions of the square effective prudent walk}). The sample space associated with the probability measure $\mathscr{P}^{*} \big[ \cdot \big] $ is given by,

\begin{align*}
    \textit{Paths}_{\infty, \textbf{Z}^2} :=  \underset{N \geq 1}{\bigcup} \bigg\{ \forall i \leq N, \exists \big(  \textit{Paths} \big)_{0 \leq i \leq N} \subsetneq \textbf{Z}^2 : \big\{ \textit{Paths}_0 =  0  \big\} ,  \big\{ \textit{Paths}_i           \geq 0 \big\} , \big\{  \textit{Paths}_N = 0  \big\}       \bigg\} . 
\end{align*}

\noindent \textbf{Definition} \textit{4} (\textit{Radon-Nikodym derivative over} $\textbf{Z}^2$). The factor of $\frac{4}{3}$ in the definition of the first, and second, probability measures above comes from the observation that,

\begin{align*}
     \frac{\mathrm{d}   \mathscr{P}^{*}        }{\mathrm{d} \mathscr{P}}  \big[ \textit{Paths}_{\infty, \textbf{Z}^2}  \big]             \propto \frac{4}{3} . 
\end{align*}

\noindent \textbf{Definition} \textit{5} (\textit{the Green's function, in terms of the expected number of visits to each} $U_i$). Fix, 

\begin{align*}
 \alpha  = \mathrm{log} \bigg\{  \textbf{E} \big[ \mathrm{exp} \big[ - \lambda \big| U_1 \big| \big] \big]    \bigg\} , 
\end{align*}

\noindent and some real $\lambda$. For the expectation,

\begin{align*}
  \mathscr{E} \big[ \cdot \big]   : =  \mathscr{E}_{\textbf{Z}^2} \big[ \cdot \big] \propto \underset{\textbf{Z}^2}{\sum} \cdot  \mathrm{d} \mathscr{P}_{\textbf{Z}^2} \big[ \cdot \big]
\end{align*}

\noindent denote the Green's function $G$ with,

\begin{align*}
   G \big( \lambda \big) =   \frac{\mathrm{exp} \big[ \alpha \big]}{3} \bigg[ \frac{1}{2} + \frac{\mathrm{exp} \big[ - \lambda ]}{4} \bigg]   + \frac{1}{2} \bigg[1 - \frac{\mathrm{exp} \big[ - \lambda ]}{ 2 } \bigg]   \mathscr{E}  \bigg\{     \mathrm{exp} \big[ \alpha \tau - \lambda \underset{1 \leq i \leq \tau}{\sum} \big| U_i \big| \big]     \bigg\}  .
\end{align*}

\noindent For the triangular kinetic prudent walk, we make use of the above five objects. In particular, we demonstrate how one can make use of measures over the sample space $\Omega_{\mathrm{Triangular}}$. Besides measures over the sample space $\Omega$, a straightforward representation of the effective random walk over $\textbf{T}$ can be inferred from that over $\textbf{Z}^2$. However, albeit the fact that one can represent probability measures supported over $\textbf{P}_{\textbf{T}} \big[ \cdot \big]$ with $\textbf{E}_{\textbf{T}} \big[ \cdot \big]$, a different factor than $\frac{4}{3}$ appears. To this end, we perform computations for determining a suitable constant of proportionality for,

\begin{align*}
 \frac{\mathrm{d}\textbf{P}^{*}_{\textit{Kinetic}, \textbf{T}} \big[  \textit{Paths}_{\infty, \textbf{T}} \big]}{\mathrm{d} \textbf{P}_{\textit{Kinetic}, \textbf{T}} \big[  \textit{Paths}_{\infty, \textbf{T}} \big] }  : =  \frac{\mathrm{d}\textbf{P}^{*}_{\textit{Kinetic}, \textbf{T}} }{\mathrm{d} \textbf{P}_{\textit{Kinetic}, \textbf{T}}  }          \big[    \textit{Paths}_{\infty, \textbf{T}}    \big] = \mathrm{d} \bigg[ \frac{\textbf{P}^*}{\textbf{P}} \bigg]  , 
\end{align*}

\noindent where,

\begin{align*}
    \textit{Paths}_{\infty, \textbf{T}} :=  \underset{N \geq 1}{\bigcup} \bigg\{ \forall i \leq N, \exists \big(  \textit{Paths} \big)_{0 \leq i \leq N} \subsetneq \textbf{T}  : \big\{ \textit{Paths}_0 =  0  \big\} , \big\{ \textit{Paths}_i           \geq 0  \big\} , \big\{  \textit{Paths}_N = 0    \big\}     \bigg\} . 
\end{align*}

\subsection{Triangular Kinetic Prudent walker objects}

\noindent Over $\textbf{T}$, the prudent walk is defined as the random walk which avoids any of the sites that it has visited in any previous steps of its excursion up to some strictly positive $t$, which we denote with $\gamma_{[0,t]}$. That is, if the prudent walk over the triangular lattice begins at the origin, given spanning set,

\begin{align*}
 \mathrm{span}_{\textbf{T}} \big\{ \vec{\textbf{e}_1} , \vec{\textbf{e}_2} , \vec{\textbf{e}_3} \big\}   , 
\end{align*}

\noindent of basis vectors,

\begin{align*}
       \textbf{e}_1 :=           \big[ \frac{\sqrt{3}}{2} , \frac{1}{2} , 0 \big]   , \\  \\ \textbf{e}_2 :=            \big[ \frac{\sqrt{3}}{2} , - \frac{1}{2} , 0 \big] , \\  \\    \textbf{e}_3 :=          \big[ 0 , 0 , 1 \big]            , 
\end{align*}

\noindent of $\textbf{T}$, one can extend the arguments from {\color{blue}[1]} to the triangular lattice by defining the times that the prudent walker occupies at the beginning, and at later times, of an excursion, which are respectively given by,

\begin{align*}
  \mathscr{T}_0 :=0  , \\ \\   \mathscr{U} :=\big\{     t > 0 :         \mathcal{H}_t > 1    \big\}  - 1 , \\ \\  \mathscr{V} : = \big\{     t > 0 :         \mathcal{W}_t > 1    \big\}  - 1  .
\end{align*}

\noindent Similar to the objects introduced in \textbf{Definition} \textit{1}-\textit{5} in \textit{2.1}, introduce:

\bigskip

\noindent \textbf{Definition} \textit{6} (\textit{the probability measure over the sample space of triangular uniform prudent paths}). Given a random walk $\mathscr{V} \subsetneq \textbf{T}$, $t \in \textbf{N}$, and the time,

\begin{align*}
    \eta^{\prime}_t =   \mathrm{min} \bigg\{ i \geq 1: i + \underset{1 \leq j \leq i}{\sum} \big| \mathscr{U}_j \big| \geq t \bigg\}       , 
\end{align*}

\noindent corresponding to the first probability measure, $\mathscr{K}$, define,

\begin{align*}
  \mathscr{K} \big( \lambda \big)  : =  \mathscr{K} = \mathscr{E} \bigg\{  \mathrm{exp} \bigg[ \mathrm{log} \bigg[ \mathrm{d} \frac{\textbf{P}^*}{\textbf{P}} \bigg] \eta_t  \bigg]                            \textbf{1}_{ \{ \mathscr{V}_i \geq 0, \forall i \leq \eta^{\prime}_t, \mathscr{V}_{\eta_t} = 0 , \eta^{\prime}_t + \underset{1 \leq j \leq \eta_t}{\sum} | \mathscr{U}_j | = t \}}         \bigg\}          . 
\end{align*}

\bigskip

\noindent \textbf{Definition} \textit{7} (\textit{the uniform probability measure over prudent walks from the probability measure over square effective prudent walks}). Denote the times $T_i$, where $1 \leq i \leq L$, for which the square uniform prudent walk occupies $i \in V \big( \textbf{Z}^2 \big)$. Corresponding to the second probability measure, $\mathscr{P}_{\textit{Uniform}}    \big[ \cdot \big] $, denote the probability measure over square effective random walks,

\begin{align*}
   \mathscr{P}^{*}_{\textbf{T}} \big[ \cdot \big] : =   \mathscr{P}^{*}    \big[ \cdot \big]  =   \textbf{P}_{\mathscr{U},\textbf{T}}  \big[   \textit{effective paths } P : \big\{ P  \textit{ lies strictly above } \textbf{T} \big\} , \big\{  P_0 = P_L = 0  \big\}    \big]          , 
\end{align*}

\noindent from which one can define,

{\small \begin{align*}
\mathscr{P}_{\textit{Uniform}}    \big[ \cdot \big]  : = \mathscr{P}_{\textit{Uniform},\textbf{T}}    \big[ \cdot \big]     =   \frac{\mathscr{E} \bigg[  \textbf{1}_{\cdot}   \times  \textbf{1}_{\underset{1 \leq i \leq L}{\sum} \widetilde{T}_i = L  }                \bigg]}{\mathscr{P}^{*} \bigg[  \underset{1 \leq i \leq L}{\sum} \widetilde{T}_i = L       \bigg] }      . 
\end{align*} } 

\noindent for the expectation,

\begin{align*}
  \mathscr{E} \big[ \cdot \big]   : =  \mathscr{E}_{\textbf{T}} \big[ \cdot \big] \propto \underset{\textbf{T}}{\sum} \cdot  \mathrm{d} \mathscr{P}_{\textbf{T}} \big[ \cdot \big] .
\end{align*}

\bigskip

\noindent \textbf{Definition} \textit{8} (\textit{nonnegative excursions of the square effective prudent walk}). The definition for $\textit{Paths}_{\infty,\textbf{T}}$ is provided at the end of the previous section, in the definition of $\mathrm{d} \frac{\textbf{P}^{*}}{\textbf{P}}$.

\bigskip

\noindent \textbf{Definition} \textit{9} (\textit{Radon-Nikodym derivative over} $\textbf{T}$). The strictly positive factor, $\mathscr{C}$, below in the definition of the first, and second, probability measures above comes from the observation that,

{\small \begin{align*}
  \bigg|  \mathrm{d} \frac{\mathscr{P}^{*}}{\mathscr{P}}    \bigg|      \propto \mathscr{C} : = \bigg|    \mathrm{d} \bigg[  \frac{\textbf{P}^{*}_{E_1}}{\textbf{P}_{E_1} + \textbf{P}_{E_2}} \bigg]   +   \mathrm{d} \bigg[ \frac{\textbf{P}^{*}_{E_2}}{\textbf{P}_{E_1} + \textbf{P}_{E_2}} \bigg]      \bigg|                    , 
\end{align*} } 

where,

{\small \begin{align*}
 \mathrm{d} \bigg[  \frac{\textbf{P}^{*}_{E_1}}{\textbf{P}_{E_1} + \textbf{P}_{E_2}} \bigg]  : =    \frac{   \mathrm{d} \textbf{P}^{*}_{E_1}}{  \mathrm{d} \textbf{P}_{E_1} +   \mathrm{d}  \textbf{P}_{E_2}}      , \\ \\           \mathrm{d} \bigg[ \frac{\textbf{P}^{*}_{E_2}}{\textbf{P}_{E_1} + \textbf{P}_{E_2}} \bigg]  : =     \frac{ \mathrm{d}  \textbf{P}^{*}_{E_2}}{ \mathrm{d}  \textbf{P}_{E_1} +  \mathrm{d}  \textbf{P}_{E_2}}       . 
\end{align*}    }

\noindent \textbf{Definition} \textit{10} (\textit{the triangular Green's function for each} $\mathscr{U}_i$). Fix, 

\begin{align*}
 \alpha^{\prime}  = \mathrm{log} \mathscr{E} \big[ \mathrm{exp} \big[ - \lambda \big| \mathscr{U}_1 \big| \big] \big]    , 
\end{align*}

\noindent and some real $\lambda$. For the expectation,

\begin{align*}
  \mathscr{E} \big[ \cdot \big]   : =  \mathscr{E}_{\textbf{T}} \big[ \cdot \big] \propto \underset{\textbf{T}}{\sum} \cdot  \mathrm{d} \mathscr{P}_{\textbf{T}} \big[ \cdot \big]
\end{align*}

\noindent denote the Green's function $G$ with,

\begin{align*}
   G_{\textbf{T}} \big( \lambda \big)    .
\end{align*}

\noindent Each excursion introduced above, independently of the initial starting point of the path, satisfies the proportionalities,

\begin{align*}
  \mathscr{T} \propto \mathrm{span}_{\textbf{T}} \big\{ \vec{e_1} , 0 , 0 \big\}  , \\ \\ \mathscr{U} \propto \mathrm{span}_{\textbf{T}} \big\{ 0 , \vec{e_2}  , 0 \big\}  ,\\ \\ \mathscr{W} \propto \mathrm{span}_{\textbf{T}} \big\{0 ,  0 , \vec{e_3} \big\}   . 
\end{align*}

\noindent respectively. In the two times defined above, the latter is dependent upon the height, $\mathcal{H}_t$, of the triangular prudent walker at $t$. Associated with some rectangular box over $\textbf{T}$, $\mathscr{R}_t \subsetneq \textbf{T}$, the stopping times $\eta_L$ introduced below will quantify the probability of whether the triangular walker exits some region of $\textbf{T}$. Relatedly, one can also define the following two other times associated with an excursion,

{\small \begin{align*}
  \mathscr{T}_{k+1} :=\big\{ t > \mathscr{U}_k :  ( \textit{Set of edges spanned by $e_1$ explored by the prudent path} )_t > (\textit{Set of edges spanned} \\ \textit{by $e_1$ explored by the prudent path} )_{t-1} \big\}  - 1 ,\\  \\ \mathscr{U}_{k+1} :=\big\{ t > \mathscr{T}_{k+1} :  ( \textit{Set of edges spanned by $e_2$ explored by the prudent path} )_t > (\textit{Set of edges spanned} \\ \textit{by $e_2$ explored by the prudent path} )_{t-1} \big\} - 1, \\  \\ \mathscr{V}_{k+1} : = \big\{ t >      \mathscr{W}_{t+2}           :    ( \textit{Set of edges spanned by $e_3$ explored by the prudent path} )_t > (\textit{Set of edges spanned} \\ \textit{by $e_3$ explored by the prudent path} )_{t-1}     \big\} - 1. 
\end{align*} }

\noindent corresponding to the times strictly greater than $\mathscr{U}_k$ for which the width of the prudent walker at time $t$, $\mathscr{W}_t$, is strictly greater than the width of the excursion of the prudent walker at time $t-1$, $\mathscr{W}_{t-1}$, and to the times strictly greater than $\mathscr{T}_{k+1}$ for which the height of the prudent walker at time $t$, $\mathcal{H}_t$, strictly exceeds the height of the prudent walker at time $t-1$, $\mathcal{H}_{t-1}$. When $\mathscr{T}_0 :=0$, denote $\mathscr{W}_0 :=\mathcal{H}_0 :=1$.

Introduce the following additional objects:

\bigskip

\noindent \textbf{Definition} \textit{11} (\textit{sample space of the effective triangular random walk}). Introduce,

\begin{align*}
   \Omega_{\textit{Eff}} : =  \big\{ \omega \in \Omega_{\mathrm{Triangular}}  : \omega \in \mathscr{E}\mathscr{W}_{\textbf{T}} \big\}  , 
\end{align*}

\noindent corresponding to the sample space of effective triangular kinetic prudent walks.

\bigskip

\noindent \textbf{Definition} \textit{12} (\textit{excursions sampled against the collection of random paths of length L}). Introduce,

\begin{align*}
  \Omega_{\textit{Eff},L}  : =   \Omega_{\textit{Eff}} \bigg|_L =   \Omega_{\textit{Eff}} \cap \big\{ \textit{effective, prudent random paths of strictly positive length L} \big\}  ,
\end{align*}

\noindent corresponding to the \textit{restriction} of paths for some $L$ taken sufficiently large.

\bigskip

\noindent Often times, we will fix an instance of the triangular effective random walk with $\widetilde{\omega} \in \Omega_{\textit{Eff}}$, and, as expected, $\widetilde{\omega_L} \in \Omega_{\textit{Eff},L}$. Clearly, $\widetilde{\omega}, \widetilde{\omega_L} \sim 
\textbf{P}^{*}$.

\bigskip

\noindent \textbf{Definition} \textit{13} (\textit{the total length, and number of increments, of triangular effective random paths}). Introduce,

\begin{align*}
  \textit{Length} \big[ \omega  \big] : =   \# \big\{\textit{vertex in } V \big( \textbf{T} \big) :  \textit{the triangular effective random path } \omega \textit{ occupies a vertex for some} \\ \textit{ strictly positive time} \big\}     , 
\end{align*}

\noindent corresponding to the total length of a triangular effective random walk, and,

\begin{align*}
  \textit{Increments} \big[ \omega \big] : = \# \big\{ \textit{times t} : \textit{the triangular effective random path } \omega \textit{ changes direction along } e_1,e_2, \\ \textit{ or }  e_3 \textit{ at}  \textit{ time t} \big\}  , 
\end{align*}

\noindent corresponding to the number of increments of a triangular effective random walk.

\bigskip

\noindent \textbf{Definition} \textit{14} (\textit{the total length, and number of increments, of triangular effective random paths of length L}). Introduce,

\begin{align*}
   \textit{Length} \big[ \omega_L   \big] : =   \# \big\{ \textit{vertex in } V \big( \textbf{T} \big)  :  \textit{the triangular effective random path } \omega \textit{ occupies a vertex} \\ \textit{for some strictly positive time, of length L} \big\}     , 
\end{align*}

\noindent corresponding to the total length of a triangular effective random walk, of length $L$, and,

\begin{align*}
     \textit{Increments} \big[ \omega_L  \big] : = \# \big\{ \textit{times t} : \textit{the triangular effective random path } \omega \textit{ changes direction along } e_1,e_2, \\ \textit{ or }    e_3 \textit{ at time t, of length L} \big\}  , 
\end{align*}

\noindent corresponding to the number of increments of a triangular effective random walk, also of length $L$.

\bigskip

\noindent \textbf{Definition} \textit{15} (\textit{truncation of the effective random walk over } $\textbf{T}$). Introduce,

\begin{align*}
    \mathcal{T} :   \Omega_{\textit{Eff}} \longrightarrow \Omega_{\textit{Eff},L}        ,
\end{align*}

\noindent corresponding to the \textit{truncation} mapping between the sample space of the triangular effective random paths and the sample space of the triangular effective random paths of length $L$.

\bigskip

\noindent \textbf{Definition} \textit{16} (\textit{ordinary, and truncated, sample spaces of the triangular effective walk over} $\textbf{T}$ under $\textbf{P}^*$, instead of under $\textbf{P}$). Introduce,

\begin{align*}
\widetilde{\Omega_{\textit{Eff}}} , \\ \\  \widetilde{\Omega_{\textit{Eff},L}}        ,
\end{align*}

\noindent corresponding to the sample spaces of the effective triangular random walk without, and with truncation over length $L$, respectively.

\bigskip

\noindent \textbf{Definition} \textit{17} (\textit{truncation of the effective random walk over } $\textbf{T}$ under $\textbf{P}^*$, instead of under $\textbf{P}$). Introduce,

\begin{align*}
    \widetilde{\mathcal{T}} :   \widetilde{\Omega_{\textit{Eff}}} \longrightarrow \widetilde{\Omega_{\textit{Eff},L}}        ,
\end{align*}

\noindent corresponding to the \textit{truncation} mapping between the sample space of the triangular effective random paths and the sample space of the triangular effective random paths of length $L$.

\bigskip

\noindent \textbf{Definition} \textit{18} (\textit{inductive definition of the length, and number of excursions, of an effective triangular random path}). Introduce,

{\small \begin{align*}
 \bigg[  \widetilde{\mathcal{T}}  , \big\{  \widetilde{U} \textit{ is the effective triangular random path along } e_i : < e_i , e_2 > = 0 \big\}   , \big\{   \widetilde{V} \textit{is the effective} \\ \textit{ triangular random path along } e^{\prime}_i : < e^{\prime}_i , e_3 > = 0    \big\}      \bigg]_{i \geq 1} \sim  \bigg[   \underset{i \geq 1}{\bigcup} \widetilde{\mathcal{T}_i} , \widetilde{\mathscr{U}} , \widetilde{\mathscr{V}}     \bigg]  ,
\end{align*} } 

\noindent where,

{\small \begin{align*}
  \widetilde{\mathscr{U}} : = \underset{\textit{paths } \widetilde{U}  }{\bigcup} \big\{ \widetilde{U} : \widetilde{U} \textit{is the effective triangular random path along } e_2 \big\}   , \\ \\ \widetilde{\mathscr{V}} : = \underset{\textit{paths } \widetilde{V} }{\bigcup} \big\{ \widetilde{V} : \widetilde{V} \textit{is the effective triangular random path along } e_3 \big\}  , 
\end{align*} } 

\noindent and,

{\small \begin{align*}
    \widetilde{\mathscr{U}_0} : =   \underset{\textit{path } }{\bigcup}      \big\{    \textit{Position of the triangular effective random path along } e_2 \textit{ at time 0}   \big\}    , \\ \\    \widetilde{\mathscr{V}_0} : = \underset{\textit{path } }{\bigcup}       \big\{ \textit{Position of the triangular effective random path along } e_3 \textit{ at time 0} \big\}   .
\end{align*} }

\noindent corresponding to an inductive definition of an effective triangular random path.

\bigskip

\noindent To define the next object, recall, {\color{blue}[4]}:

\bigskip

\noindent \textbf{Definition} \textit{19} (\textit{Slabs over } $\textbf{Z}^2$ - \textit{square slabs}). Fix $r,x \in \textbf{N}$. Introduce,

\begin{align*}
    \mathscr{S} \mathscr{S}_{\textbf{Z}^2,r} : =    \big[ 0 , x \big]     \times   \big[ 0 , r \big]                              ,
\end{align*}

\noindent corresponding to slabs over $\textbf{Z}^2$.

\bigskip

\noindent \textbf{Definition} \textit{20} (\textit{stopping times before which the effective random walk exists a triangular slab, along increments of either } $\mathscr{U}$, or $\mathscr{V}$). Fix some $r^{\prime} \in \textbf{N}$. Denote a triangular slab of length $r$ with,

\begin{align*}
    \mathscr{T}\mathscr{S}_{\textbf{T}, r^{\prime}} : =            \big\{ \textbf{Z}^2 \textit{ slab} \big\}              \times \big[ 0 , r^{\prime} \big]   ,
\end{align*}

\noindent Introduce,

\begin{align*}
    \tau : = \mathrm{inf} \bigg\{             i \geq 1 : \big\{ \mathscr{V}_i \geq R + 1 \big\} \textit{ or } \big\{ \mathscr{U}_i \geq R+1 \big\}            \bigg\}      ,
\end{align*}

\noindent corresponding to the time at which the triangular effective random path exists a slab over $\textbf{T}$.

\bigskip

\noindent The computation of the Radon-Nikodym derivative stated in \textbf{Definition} \textit{9} can be alternatively expressed with the following objects. As a matter of shorthand, we denote,

{\small \begin{align*}
   \mathrm{d} \bigg[  \frac{\textbf{P}^{*}}{\textbf{P}} \bigg]  :=  \mathrm{d} \bigg[ \frac{\mathscr{P}^{*}}{\mathscr{P}} \bigg]  =  \mathrm{d} \bigg[ \frac{\textbf{P}^{*}_{\mathscr{U}, \textbf{T}}}{\textbf{P}_{\mathscr{U}, \textbf{T}}} \bigg] .
\end{align*} } 

\noindent \textbf{Definition} \textit{21} (\textit{Radon-Nikodym derivative about each degree of freedom of the triangular effective random walk}). Fix strictly positive $\mathscr{C}_{E_1}, \mathscr{C}_{E_2}$, with $\mathscr{C}_{E_1}, \mathscr{C}_{E_2} < \mathscr{C}$, and,

{\small \begin{align*}
   E_1 : =  \underset{e \in E_1}{\bigcup}  \big\{ \textit{e : e is an edge of the triangular effective random walk along the first degree} \\ \textit{of freedom} \big\} , \\ \\ E_2 : =  \underset{e \in E_2}{\bigcup}  \big\{ \textit{e : e is an edge of the triangular effective random walk along the second degree} \\ \textit{of freedom} \big\}   , 
\end{align*} } 

\noindent Introduce,

{\small \begin{align*}
   \mathrm{d} \bigg[  \frac{\textbf{P}^{*}}{\textbf{P}} \bigg]  : =  \mathrm{d} \bigg[  \frac{\textbf{P}^{*}_{E_1} + \textbf{P}^{*}_{E_2}}{\textbf{P}_{E_1} + \textbf{P}_{E_2}} \bigg]  =    \mathrm{d} \bigg[  \frac{\textbf{P}^{*}_{E_1}}{\textbf{P}_{E_1} + \textbf{P}_{E_2}} \bigg]  +  \mathrm{d} \bigg[ \frac{\textbf{P}^{*}_{E_2}}{\textbf{P}_{E_1} + \textbf{P}_{E_2}} \bigg]       \propto \mathscr{C} = \mathscr{C}_{E_1} + \mathscr{C}_{E_2}                     . 
\end{align*} }

\noindent  corresponding to a decomposition of the Radon-Nikodym derivative,

\begin{align*}
   \mathrm{d} \bigg[ \frac{\textbf{P}^{*}}{\textbf{P}}  \bigg]  , 
\end{align*}

\noindent and probability measures,

\begin{align*}
     \textbf{P}_{E_1} : =  \textit{Triangular effective probability measures over the degree of freedom spanned by } E_1 \in E \big( \textbf{T} \big)    , \\ \\ \textbf{P}_{E_2} : =  \textit{Triangular effective probability measures over the degree of freedom spanned by } E_2 \in E \big( \textbf{T} \big) , 
\end{align*}

\noindent over paths of the triangular effective random walk, and where,

\begin{align*}
  \mathrm{d} \bigg[  \frac{\textbf{P}^{*}_{E_1} + \textbf{P}^{*}_{E_2}}{\textbf{P}_{E_1} + \textbf{P}_{E_2}} \bigg] : =      \frac{\mathrm{d}  \big[ \textbf{P}^{*}_{E_1} + \textbf{P}^{*}_{E_2}\big] }{\mathrm{d}  \big[ \textbf{P}_{E_1} + \textbf{P}_{E_2} \big] }  =     \frac{\mathrm{d}  \textbf{P}^{*}_{E_1} + \mathrm{d} \textbf{P}^{*}_{E_2} }{\mathrm{d}   \textbf{P}_{E_1} + \mathrm{d} \textbf{P}_{E_2}  }                ,   \\ \\ \mathrm{d} \bigg[  \frac{\textbf{P}^{*}_{E_1}}{\textbf{P}_{E_1} + \textbf{P}_{E_2}} \bigg]  : =    \frac{   \mathrm{d} \textbf{P}^{*}_{E_1}}{  \mathrm{d} \textbf{P}_{E_1} +   \mathrm{d}  \textbf{P}_{E_2}}      , \\ \\           \mathrm{d} \bigg[ \frac{\textbf{P}^{*}_{E_2}}{\textbf{P}_{E_1} + \textbf{P}_{E_2}} \bigg]  : =     \frac{ \mathrm{d}  \textbf{P}^{*}_{E_2}}{ \mathrm{d}  \textbf{P}_{E_1} +  \mathrm{d}  \textbf{P}_{E_2}}       . 
\end{align*}

\noindent \textbf{Definition} \textit{22} (\textit{excursions of the triangular effective random walk about the second degree of freedom}). Besides a previous expression for $\textit{Paths}_{\infty, \textbf{T}}$, introduce, 

\begin{align*}
  \textit{Paths}^{\prime}_{\infty, \textbf{T}} : =      \underset{N \geq 1}{\bigcup} \bigg\{  \forall i \leq N, \exists \big(  \textit{Paths}^{\prime} \big)_{0 \leq i \leq N} \subsetneq \textit{second degree of freedom of the triangular ef-} \\ \textit{ fective random walk over } \textbf{T}  : \big\{ \textit{Paths}^{\prime}_0 =  0 \big\}  ,  \big\{ \textit{Paths}^{\prime}_i           \geq 0 \big\} , \big\{ \textit{Paths}^{\prime}_N = 0   \big\}      \bigg\}                                              ,
\end{align*}

\noindent corresponding to the set of excursions of the triangular prudent path about each $e \in E_2$.

\bigskip

\noindent \textbf{Definition} \textit{23} (\textit{slab measures over} $\textbf{N} \times \textbf{N} \times\textbf{N} \times \big[ 0 , 1 \big]$). Fix,

\begin{align*}
    R_1 : =   \textit{First degree of freedom of the triangular effective walk over the slab}          ,  \\ \\ R_2 : =          \textit{Second degree of freedom of the triangular effective walk over the slab}   ,
\end{align*}

\noindent and also,

\begin{align*}
 n : =  \textit{displacement of a triangular effective path along the first degree of freedom}  , \\ \\ n^{\prime} : = \textit{displacement of a triangular effective path along the second degree of fre-} \\ \textit{edom}  . 
\end{align*}

\noindent Introduce,

{\small \begin{align*}
   \mathscr{L}_0  : =     \mathscr{E}  \bigg\{      \textbf{1}_{\{ R = R_1 \} }                \mathrm{exp}     \bigg[   \mathrm{d} \bigg[  \frac{\textbf{P}^{*}_{E_1}}{\textbf{P}_{E_1} + \textbf{P}_{E_2}} \bigg]      n       \textbf{1}_{\{    \widetilde{\mathscr{U}} \in \{ 0 , \cdots, R \}  : \widetilde{\mathscr{U}}_n = 0 ,  \underset{1 \leq i \leq n}{\sum} | \widetilde{\mathscr{U}}_i |         = t - n                                                  \}}                      \bigg]                               +   \textbf{1}_{\{ R = R_2 \} }       \mathrm{exp} \bigg[    \mathrm{d} \bigg[ \frac{\textbf{P}^{*}_{E_2}}{\textbf{P}_{E_1} + \textbf{P}_{E_2}}       \bigg]    n^{\prime}       \\ \times    \textbf{1}_{\{  \widetilde{\mathscr{V}} \in \{ 0 , \cdots, R \}  :   \widetilde{U}_n = 0 ,  \underset{1 \leq i \leq n}{\sum} | \widetilde{\mathscr{V}}_i  |                 =    t - n^{\prime}          \}}                      \bigg]         \bigg\}   ,
\end{align*}

\begin{align*}
   \mathscr{L}_R  : =     \mathscr{E}  \bigg\{      \textbf{1}_{\{ R = R_1 \} }                \mathrm{exp}     \bigg[   \mathrm{d} \bigg[  \frac{\textbf{P}^{*}_{E_1}}{\textbf{P}_{E_1} + \textbf{P}_{E_2}} \bigg]       n       \textbf{1}_{\{    \widetilde{\mathscr{U}} \in \{ 0 , \cdots, R \}  : \widetilde{U}_n = R ,  \underset{1 \leq i \leq n}{\sum} | \widetilde{\mathscr{U}}_i |         = t - n                                                  \}}                      \bigg]                               +   \textbf{1}_{\{ R = R_2 \} }       \mathrm{exp} \bigg[    \mathrm{d} \bigg[ \frac{\textbf{P}^{*}_{E_2}}{\textbf{P}_{E_1} + \textbf{P}_{E_2}}       \bigg]  n^{\prime}         \\ \times    \textbf{1}_{\{  \widetilde{\mathscr{V}} \in \{ 0 , \cdots, R \}  :   \widetilde{U}_n = R ,  \underset{1 \leq i \leq n}{\sum} | \widetilde{\mathscr{V}}_i  |                 =    t - n^{\prime}          \}}                      \bigg]         \bigg\}   ,
\end{align*} }

\noindent corresponding to the expectation that excursions of the triangular effective random walk about $R_1$, and $R_2$, over slabs of length $0$ and $R$, respectively. For shorthand, we denote $\mathscr{L}_0 : = \mathscr{L}_0 \big[ n , n^{\prime} , 0 , R_1 , R_2 \big]$ and $\mathscr{L}_R : = \mathscr{L}_R \big[ n , n^{\prime} , R, R_1 , R_2 \big]$.

\bigskip

\noindent \textbf{Definition} \textit{24} (\textit{slab measures over} $\textbf{N} \times \textbf{N} \times\textbf{N} \times \big[ 0 , 1 \big]$ \textit{for incomplete excursions}). Fix,

\begin{align*}
 n : =  \textit{displacement of a triangular effective path along the first degree of freedom}  , \\ \\ n^{\prime} : = \textit{displacement of a triangular effective path along the second degree of freedom}  . 
\end{align*}

\noindent Introduce,

{\small \begin{align*}
  \hat{\mathscr{L}_0} : =        \mathscr{E}  \bigg\{      \textbf{1}_{\{ R = R_1 \} }                \mathrm{exp}     \bigg[   \mathrm{d} \bigg[  \frac{\textbf{P}^{*}_{E_1}}{\textbf{P}_{E_1} + \textbf{P}_{E_2}} \bigg]       n      \textbf{1}_{\{    \widetilde{\mathscr{U}} \in \{ 0 , \cdots, R \}  :  \widetilde{U}_n = 0  ,  \underset{1 \leq i \leq n}{\sum} | \widetilde{\mathscr{U}}_i |         = t - n                                                  \}}                      \bigg]                               +   \textbf{1}_{\{ R = R_2 \} }       \mathrm{exp} \bigg[    \mathrm{d} \bigg[ \frac{\textbf{P}^{*}_{E_2}}{\textbf{P}_{E_1} + \textbf{P}_{E_2}}       \bigg]    n^{\prime}       \\ \times    \textbf{1}_{\{  \widetilde{\mathscr{V}} \in \{ 0 , \cdots, R \}  : \widetilde{U}_n = 0 ,  \underset{1 \leq i \leq n}{\sum} | \widetilde{\mathscr{V}}_i  |                 =    t - n^{\prime}          \}}                      \bigg]         \bigg\}    , \\   \\  \hat{\mathscr{L}_R} : =        \mathscr{E}  \bigg\{      \textbf{1}_{\{ R = R_1 \} }                \mathrm{exp}     \bigg[   \mathrm{d} \bigg[  \frac{\textbf{P}^{*}_{E_1}}{\textbf{P}_{E_1} + \textbf{P}_{E_2}} \bigg]       n      \textbf{1}_{\{    \widetilde{\mathscr{U}} \in \{ 0 , \cdots, R \}  : 0 < \widetilde{U}_n <  R ,  \underset{1 \leq i \leq n}{\sum} | \widetilde{\mathscr{U}}_i |         = t - n                                                  \}}                      \bigg]                               +   \textbf{1}_{\{ R = R_2 \} }       \mathrm{exp} \bigg[    \mathrm{d} \bigg[ \frac{\textbf{P}^{*}_{E_2}}{\textbf{P}_{E_1} + \textbf{P}_{E_2}}       \bigg]    n^{\prime}       \\ \times    \textbf{1}_{\{  \widetilde{\mathscr{V}} \in \{ 0 , \cdots, R \}  :   0 < \widetilde{U}_n <  R ,  \underset{1 \leq i \leq n}{\sum} | \widetilde{\mathscr{V}}_i  |                 =    t - n^{\prime}          \}}                      \bigg]         \bigg\}    ,
\end{align*} }

\noindent corresponding to a slab measure over $\textbf{N} \times \textbf{N} \times\textbf{N} \times \big[ 0 , 1 \big]$ for \textit{incomplete} excursions, over slabs of length $0$ and $R$, respectively. For shorthand, we denote $\hat{\mathscr{L}_0} : = \hat{\mathscr{L}_0}  \big[ n , n^{\prime} , 0 , R_1 , R_2 \big]$ and $\hat{\mathscr{L}_R} : = \big[ n , n^{\prime} , R , R_1 , R_2 \big] $.

\bigskip

\noindent \textbf{Definition} \textit{25} (\textit{the set of triangular effective paths within a slab that end at 0, or at R}). Introduce,

\begin{align*}
   \big\{  \Omega_{\textit{Eff}} \big\}  \cap       \mathscr{T}\mathscr{S}_{\textbf{T}, R}             ,
\end{align*}

\noindent corresponding to the set of triangular effective paths, contained within a slab, which end at $R$,

\begin{align*}
   \big\{  \Omega_{\textit{Eff}} \big\}  \cap       \mathscr{T}\mathscr{S}_{\textbf{T}, 0}             ,
\end{align*}

\noindent corresponding to the set of triangular effective paths, contained within a slab, which end at $0$, and,

\begin{align*}
   \big\{  \Omega_{\textit{Eff}} \big\}  \cap       \mathscr{T}\mathscr{S}_{\textbf{T}, R}      \cap       \mathscr{T}\mathscr{S}_{\textbf{T}, 0}                 ,
\end{align*}

\noindent corresponding to the set of triangular effective paths, contained within a slab, which end at $0$, or at $R$, for,

\begin{align*}
    \mathscr{T}\mathscr{S}_{\textbf{T}, 0}  : =        \big\{ \textbf{Z}^2 \textit{ slab} \big\}              \times \big\{ 0 \big\}   . 
\end{align*}

\bigskip

\noindent In light of the expressions provided in \textbf{Definition} \textit{20}, and \textit{21}, the measure,

\begin{align*}
  \mathscr{K} \big( \lambda \big)  : =  \mathscr{K} = \mathscr{E} \bigg\{  \mathrm{exp} \bigg[ \mathrm{log} \bigg[ \mathrm{d} \frac{\textbf{P}^*}{\textbf{P}} \bigg] \eta_t  \bigg]                            \textbf{1}_{ \{ \mathscr{V}_i \geq 0, \forall i \leq \eta^{\prime}_t, \mathscr{V}_{\eta_t} = 0 , \eta^{\prime}_t + \underset{1 \leq j \leq \eta_t}{\sum} | \mathscr{U}_j | = t \}}         \bigg\}       , 
\end{align*}

\noindent for the triangular effective walk takes the form,

\begin{align*}
  \hat{\mathscr{K}} \big( \lambda \big)  : =  \hat{\mathscr{K}} = \mathscr{E}  \bigg\{      \textbf{1}_{\{ R = R_1 \} }                \mathrm{exp} \bigg[      \mathrm{d} \bigg[  \frac{\textbf{P}^{*}_{E_1}}{\textbf{P}_{E_1} + \textbf{P}_{E_2}} \bigg]         n      \textbf{1}_{\{    \widetilde{\mathscr{U}} \in \{ 0 , \cdots, R \}  : \widetilde{\mathscr{U}}_n = 0 ,  \underset{1 \leq i \leq n}{\sum} | \widetilde{\mathscr{U}}_i |         = t - n                                                  \}}                      \bigg]                               +   \textbf{1}_{\{ R = R_2 \} }    \\ \times     \mathrm{exp} \bigg[    \mathrm{d} \bigg[ \frac{\textbf{P}^{*}_{E_2}}{\textbf{P}_{E_1} + \textbf{P}_{E_2}}       \bigg]      n^{\prime}   \textbf{1}_{\{  \widetilde{\mathscr{V}} \in \{ 0 , \cdots, R \}  :   \widetilde{U}_n = 0 ,  \underset{1 \leq i \leq n}{\sum} | \widetilde{\mathscr{V}}_i  |                 =    t - n^{\prime}          \}}                      \bigg]  \bigg]        \bigg\}          . 
\end{align*}

\noindent \textbf{Definition} \textit{26} (\textit{excursions of the triangular effective random walk about each degree oof freedom}). Fix the number of steps $k \in \textbf{N}$ of the triangular effective random walk. Introduce,

\begin{align*}
 {E}^{\textit{first degree of freedom over } \textbf{T}}_k  , \\ \\  {E}^{\textit{second degree of freedom over } \textbf{T}}_k , \\ \\ {E}^{\textit{third degree of freedom over } \textbf{T}}_k , 
\end{align*}

\noindent corresponding to the displacements $E$ of the triangular effective walk about each degree of freedom over $\textbf{T}$, where,

\begin{align*}
   {E}^{\textit{first degree of freedom over } \textbf{T}}_k \neq \emptyset \Longleftrightarrow \textit{times } t \in \big( \mathscr{T}_k , \mathscr{U}_k \big] , \\ \\ {E}^{\textit{second degree of freedom over } \textbf{T}}_k   \neq \emptyset \Longleftrightarrow \textit{times } t \in \big( \mathscr{U}_k , \mathscr{T}_{k+1} \big]     , \\ \\  {E}^{\textit{third degree of freedom over } \textbf{T}}_k    \neq \emptyset \Longleftrightarrow     \textit{times } t \in \big( \mathscr{T}_k , \mathscr{W}_{k+1} \big]           . 
\end{align*}

\noindent \textbf{Definition} \textit{27} (\textit{displacement of the triangular effective walk}). Denote the displacements of the walker over $\textbf{T}$ about each degree of freedom with,

\begin{align*}\mathscr{X}_k :=\mathscr{W}_{\mathscr{T}_{k+1}} - \mathscr{W}_{\mathscr{T}_k} , \\ \\  \mathscr{Y}_k :=\mathcal{H}_{\mathscr{T}_{k+1}} - \mathcal{H}_{\mathscr{T}_k}, \\ \\ \mathscr{Z}_k : =   \mathscr{V}_{\mathscr{T}_{k+1}}  - \mathscr{V}_{\mathscr{T}_k}       ,
\end{align*}

\noindent corresponding to the displacement about the first, second, and third, degrees of freedom, respectively. 

\bigskip

\noindent \textbf{Definition} \textit{28} (\textit{computing the displacement of paths of the triangular effective walk}). With the displacements defined in the previous definition above, introduce, for the effective random walk, $\mathcal{S}_n$, over $\textbf{T}$, with,

\begin{align*}
 \mathcal{S}_0 :=0  ,  \\   \\  \mathcal{S}_n   :=\underset{1 \leq i \leq n}{\sum}  \xi_i, 
\end{align*}

\noindent the quantities,

{\small \begin{align*}
\textit{Exit times from the interval with Lebesgue measure $L-1$}: \eta_L,  \\ \\ 
\textit{The first exit time for which the triangular prudent walker exists the interval} \\ \textit{with Lebesgue measure $L-1$}: \big( \eta_L \big)_1 :=\mathrm{inf} \big\{ n>0 : \mathcal{S}_n \not\in \big[ 0 , L - 1 \big] \big\}, \\ \\ \textit{The probability density of $\big\{ \xi_i :=k \big\}$}: \textbf{P} \big[ \xi_i = k \big] = \frac{1}{3} \bigg\{  \frac{1}{2} \bigg\}^{|k|}. 
\end{align*} }

\noindent \textbf{Definition} \textit{29} (\textit{normalizing constant over the triangular lattice}). Denote the triangular normalizing constant as,

{\small \begin{align*}
 Z \bigg[ \textbf{T} , \sigma_1 , \sigma_2 , \sigma_3 , u , \alpha \bigg]   :=Z^{\sigma_1 , \sigma_2 , \sigma_3}_{u , \alpha} :=  \int_0^{\alpha u }      \bigg[       \sigma_1 \textbf{1}_{\{W(s) < 0 \}}      \vec{\textbf{e}_1}        + \sigma_2  \textbf{1}_{\{ W(s) = 0 \}}    \vec{\textbf{e}_2}            + \sigma_3   \textbf{1}_{\{W(s) > 0\} }    \vec{\textbf{e}_3}           \bigg]  \mathrm{d} s , 
\end{align*}} 

\noindent for,

{\small \begin{align*}
 \alpha \geq \bigg| \bigg|   \int_0^{\alpha u }   \bigg[         \sigma_1 \textbf{1}_{\{W(s) < 0 \}}      \vec{\textbf{e}_1}        + \sigma_2  \textbf{1}_{\{ W(s) = 0 \}}    \vec{\textbf{e}_2}            + \sigma_3   \textbf{1}_{\{W(s) > 0\} }    \vec{\textbf{e}_3}           \bigg]   \mathrm{d} s        \bigg| \bigg|_1   . 
\end{align*} }

\noindent In comparison to the $\frac{3}{7}$ convergence rate obtained for the normalizing constant in the scaling limit of the square kinetic prudent walk, {\color{blue}[1]}, that obtained for the triangular kinetic prudent walk is dependent upon the following properties:

\bigskip

\noindent \textbf{Definition} \textit{30} (\textit{properties of the normalizing constant}). The above normalization, over $\textbf{T}$, satisfies,

\begin{itemize}
    \item[$\bullet$] \textit{Strict positivity}. The times at which paths of Brownian motion over $\textbf{T}$ are strictly larger than, equal to, or strictly less than, $0$,

    \item[$\bullet$] \textit{The L-1 norm bound is dependenty upon the below three L-1 norms}. Bounds on the $L$-1 norms of,

       {\small  \begin{align*}
          \int_0^{\alpha u }    \bigg| \bigg|         \sigma_1 \textbf{1}_{\{W(s) < 0 \}}      \vec{\textbf{e}_1}    \bigg| \bigg|_1   \mathrm{d} s         =      \int_0^{\alpha u }   \textbf{1}_{\{W(s) < 0 \}}     \bigg| \bigg|         \sigma_1    \vec{\textbf{e}_1}    \bigg| \bigg|_1   \mathrm{d} s         , \end{align*}

          \begin{align*} \int_0^{\alpha u }    \bigg| \bigg|   \sigma_2  \textbf{1}_{\{ W(s) = 0 \}}    \vec{\textbf{e}_2}     \bigg| \bigg|_1       \mathrm{d} s        =       \int_0^{\alpha u }   \textbf{1}_{\{ W(s) = 0 \}}     \bigg| \bigg|   \sigma_2   \vec{\textbf{e}_2}     \bigg| \bigg|_1       \mathrm{d} s          , \\ \\    \int_0^{\alpha u }    \bigg| \bigg|    \sigma_3   \textbf{1}_{\{W(s) > 0\} }    \vec{\textbf{e}_3}     \bigg| \bigg|_1   \mathrm{d} s            =  \int_0^{\alpha u }  \textbf{1}_{\{W(s) > 0\} }     \bigg| \bigg|    \sigma_3     \vec{\textbf{e}_3}     \bigg| \bigg|_1   \mathrm{d} s     . 
        \end{align*} }

\end{itemize}

\noindent \textbf{Definition} \textit{31} (\textit{indicator function for determining whether the triangular effective random walk crosses the maximum displaced distance about each degree of freedom}). Introduce,

{\small \begin{align*}
\textbf{1}^{\prime} : =   \textbf{1}  \bigg\{  \mathscr{W}  \textit{ crosses the maximum displaced distance about the first degree of freedom over }    \textbf{T}        \bigg\}         , \\ \\ \textbf{2}^{\prime} : =   \textbf{1}  \bigg\{  \mathscr{W}  \textit{ crosses the maximum displaced distance about the second degree of freedom over }    \textbf{T}        \bigg\}   , \\ \\ \textbf{3}^{\prime} : =   \textbf{1}  \bigg\{  \mathscr{W}  \textit{ crosses the maximum displaced distance about the third degree of freedom over }    \textbf{T}        \bigg\}    . 
\end{align*} }

\noindent corresponding to the indicator functions that a triangular effective random path, $\mathscr{W}$, crosses the maximum displaced distance about each degree of freedom over $\textbf{T}$.

\bigskip

\noindent \textbf{Definition} \textit{32} (\textit{uniform, three-sided, representation of the triangular effective random path}). For a subset of the sample space of the triangular prudent walk, 

\begin{align*}
    \Omega_{\{ \textit{U,3-S} \} } \subsetneq  \Omega_{\mathrm{Triangular}}    ,
\end{align*}

\noindent introduce,

{\small \[   \Omega_{ \{ \textit{U,3-S} \} } : =           \left\{\!\begin{array}{ll@{}>{{}}l} 
  \textit{(1). Triangular, effective paths } \pi, \textit{ where } \pi      \textit{begins with an east step at } \big( 1, 1, 0 \big)      ,  \\ \\ \textit{(2). The endpoint } \pi_L, \textit{of } \pi,  \textit{is located at the top left corner of the smallest volume over } \textbf{T} \\ \textit{containing } \pi         ,  \\  \\ \textit{(3). }  \pi \textit{ takes no steps in the quadrant spanned by the set of linear combinations } \underset{v \in \textbf{R}}{\mathrm{span}} \big\{   \big( - v , 0 , 0 \big)  \big\}    \\ \subsetneq \mathrm{span} \textbf{T}    ,
\end{array}\right. 
\]  }

\noindent corresponding to sample space of all uniform, three-sided, paths with support over $\textbf{T}$.

\bigskip

\noindent \textbf{Definition} \textit{33} (\textit{expected value of the slab, and renormalized, slab measures}). Introduce,

{\small \begin{align*}
  \hat{\mathscr{K}^*} \big( \lambda  \big)  : =  \hat{\mathscr{K}^*} = \mathscr{E}  \bigg\{      \textbf{1}_{\{ R = R_1 \} }              \bigg\{   \mathrm{exp} \bigg[      \mathrm{d} \bigg[  \frac{\textbf{P}^{*}_{E_1}}{\textbf{P}_{E_1} + \textbf{P}_{E_2}} \bigg]          n   \textbf{1}_{\{    \widetilde{\mathscr{U}} \in \{ 0 , \cdots, R \}  : \widetilde{\mathscr{U}}_n = 0 ,  \underset{1 \leq i \leq n}{\sum} | \widetilde{\mathscr{U}}_i |         = t - n                                                  \}}      \bigg\} \\ \times  \frac{P_1}{ \big( \textbf{P}_{E_1} + \textbf{P}_{E_2} \big) \bigg[   T > \mathscr{L} - \underset{1 \leq i \leq \mathrm{log} [ \mathscr{L}^2 ]}{\sum} \mathscr{T}_i  \bigg] }                                          \\      +   \textbf{1}_{\{ R = R_2 \} }       \bigg\{ \mathrm{exp} \bigg[    \mathrm{d} \bigg[ \frac{\textbf{P}^{*}_{E_2}}{\textbf{P}_{E_1} + \textbf{P}_{E_2}}       \bigg] n^{\prime}  \textbf{1}_{\{  \widetilde{\mathscr{V}} \in \{ 0 , \cdots, R \}  :   \widetilde{U}_n = 0 ,  \underset{1 \leq i \leq n}{\sum} | \widetilde{\mathscr{V}}_i  |                 =    t - n^{\prime}          \}}                      \bigg]   \bigg]        \bigg\}       \\  \times  \frac{P_2}{ \big( \textbf{P}_{E_1} + \textbf{P}_{E_2} \big) \bigg[   T > \mathscr{L} - \underset{1 \leq i \leq \mathrm{log} [ \mathscr{L}^2 ]}{\sum} \mathscr{T}_i \bigg] }                   \bigg\}                    , 
\end{align*} }

%


\noindent corresponding to the expected value of the slab, and renormalized, slab measures, where,

\begin{align*}
   P_1 : =       \mathrm{exp}     \bigg\{    \mathrm{d} \bigg[  \frac{\textbf{P}^{*}_{E_1}}{\textbf{P}_{E_1} + \textbf{P}_{E_2}} \bigg]        n     \textbf{1}_{\{    \widetilde{\mathscr{U}} \in \{ 0 , \cdots, R \}  : 0 < \widetilde{U}_n <  R ,  \underset{1 \leq i \leq n}{\sum} | \widetilde{\mathscr{U}}_i |         = t - n                                                  \}}                      \bigg\}            , \\ \\ P_2 : = \mathrm{exp} \bigg\{      \mathrm{d} \bigg[ \frac{\textbf{P}^{*}_{E_2}}{\textbf{P}_{E_1} + \textbf{P}_{E_2}}       \bigg]  n^{\prime}       \textbf{1}_{\{  \widetilde{\mathscr{V}} \in \{ 0 , \cdots, R \}  :   0 < \widetilde{U}_n <  R ,  \underset{1 \leq i \leq n}{\sum} | \widetilde{\mathscr{V}}_i  |                 =    t - n^{\prime}          \}}                      \bigg\}       ,
 \end{align*}

 \noindent and,

 {\small \begin{align*}
          \big( \textbf{P}_{E_1} + \textbf{P}_{E_2} \big) \bigg[   T > \mathscr{L} - \underset{1 \leq i \leq \mathrm{log} [ \mathscr{L}^2 ]}{\sum} \mathscr{T}_i  \bigg] =   \textbf{P}_{E_1} \bigg[   T > \mathscr{L} - \underset{1 \leq i \leq \mathrm{log} [ \mathscr{L}^2 ]}{\sum} \mathscr{T}_i  \bigg]   + \textbf{P}_{E_2}  \bigg[   T > \mathscr{L} - \underset{1 \leq i \leq \mathrm{log} [ \mathscr{L}^2 ]}{\sum} \mathscr{T}_i  \bigg]      , \\ \\     \big( \textbf{P}_{E_1} + \textbf{P}_{E_2} \big) \bigg[   T > \mathscr{L} - \underset{1 \leq i \leq \mathrm{log} [ \mathscr{L}^2 ]}{\sum} \mathscr{T}_i \bigg] =   \textbf{P}_{E_1}  \bigg[   T > \mathscr{L} - \underset{1 \leq i \leq \mathrm{log} [ \mathscr{L}^2 ]}{\sum} \mathscr{T}_i \bigg]   + \textbf{P}_{E_2} \bigg[   T > \mathscr{L} - \underset{1 \leq i \leq \mathrm{log} [ \mathscr{L}^2 ]}{\sum} \mathscr{T}_i \bigg]   . 
 \end{align*}
} 

\noindent \textbf{Definition} \textit{34} (\textit{renormalization paths over } $\textbf{T}$ \textit{with the scaling limit normalization factor} $\mathscr{L}$). Introduce,

\begin{align*}
    \widetilde{\pi} : = \frac{1}{\mathscr{L}} \bigg\{                \big[    \pi_{\lfloor R t \rfloor }     +    \big[      t R - \lfloor t R \rfloor  \big] \big[   \pi_{\lfloor t R \rfloor + 1} - \pi_{\lfloor t R \rfloor}    \big]    \big]                 \textbf{1}_{\{ R = R_1 \} }      + \big[    \pi_{\lfloor R t \rfloor }      +      \big[      t R - \lfloor t R \rfloor  \big] \big[   \pi_{\lfloor t R \rfloor + 1} \\ - \pi_{\lfloor t R \rfloor}    \big]    \big] \textbf{1}_{\{ R = R_2 \} }                   \bigg\}    ,
\end{align*}

\noindent corresponding to the renormalization of the path $\pi \in \Omega_{\mathrm{Triangular}}$, $\widetilde{\pi}$, and a $t \in \textbf{N}$.

\bigskip

\noindent As objects of closely related interest, the slab measures over $\textbf{N} \times \textbf{N} \times \textbf{N} \times \big[ 0 , 1 \big]$ that have been introduced above with \textbf{Definition} \textit{23} are related to \textit{truncation} procedures for random walks over $\textbf{Z}^2$. From objects manipulated in {\color{blue}[1]}, we recall:

\bigskip

\noindent \textbf{Definition} \textit{35} (\textit{truncation procedures for horizontal, and vertical, excursions of random walks over} $\textbf{Z}^2$). For,

\begin{align*}
  H_0 :=0 , 
\end{align*}

\noindent denote the \textit{vertical} truncation of a random walk over $\textbf{Z}^2$, with,

\[
\mathcal{E}^{\mathrm{v}}_1 := \text{ } 
\left\{\!\begin{array}{ll@{}>{{}}l}         \hat{\mathcal{E}^{\mathrm{v}}_1} 
 \text{ } \text{if} \text{ } m^{\mathrm{v}}_1 \leq H_0  \text{ , } \\  
\mathrm{Trunc}_{H_0} \big( \hat{\mathcal{E}^{\mathrm{v}}_1}  \big)   \text{ }  \text{if} \text{ } m^{\mathrm{v}}_1 > H_0  
 \text{ . }  \\
\end{array}\right.
\]

\noindent After applying $\mathcal{E}^{\mathrm{v}}_1$ to the random walk, one obtains the truncation,

\begin{align*}
  \mathrm{Trunc}_{H_0} \big( \hat{\mathcal{E}^{\mathrm{v}}_1}  \big)   = \hat{\mathcal{E}^{\mathrm{v}}_1}  \cap    H_0     .
\end{align*}

\noindent In place of $H_0$, one can straightforwardly take the point of truncation for the \textit{vertical} excursions of $\mathcal{E}$ with $H_k$, from which the truncation takes the form,

\[
\mathcal{E}^{\mathrm{v}}_{k+1} := \text{ } 
\left\{\!\begin{array}{ll@{}>{{}}l}         \hat{\mathcal{E}^{\mathrm{v}}_k} 
 \text{ } \text{if} \text{ } m^{\mathrm{v}}_k \leq H_k  \text{ , } \\  
\mathrm{Trunc}_{H_k} \big( \hat{\mathcal{E}^{\mathrm{v}}_{k+1}}  \big)   \text{ }  \text{if} \text{ } m^{\mathrm{v}}_k > H_0  
 \text{ . }  \\
\end{array}\right.
\]

\noindent After applying $\mathcal{E}^{\mathrm{v}}_{k+1}$ to the random walk, one obtains the truncation,

\begin{align*}
   \mathrm{Trunc}_{H_k} \big( \hat{\mathcal{E}^{\mathrm{v}}_{k+1}}  \big)  = \hat{\mathcal{E}^{\mathrm{v}}_{k+1}} \cap H_k .  
\end{align*}

\noindent Straightforwardly, one can obtain define similar expressions corresponding to \textit{horizontal} truncations, in addition to the path one obtains after applying the truncation.

\bigskip

\noindent \textbf{Definition} \textit{36} (\textit{Stopping times of the effective random walk for previously obtained scaling limits over } $\textbf{Z}^2$). Introduce the stopping time of an effective random walk over $\textbf{Z}^2$, before the first step is taken,

\begin{align*}
 \tau_0 :=0   . 
\end{align*}

\noindent Denote the following two stopping times for the square effective walk, $S_n$, as time step $n$ with,

\begin{align*}
 \tau_{2k+1} :=            \big\{ n > \tau_{2k}  : S_n < S_{\tau_{2k}} \big\}   , \\  \\   \tau_{2k+2} :=   \big\{      n > \tau_{2k+1} : S_n > S_{\tau_{2k+1}}      \big\}      . 
\end{align*}

\noindent At any instance of time, $S_n$ and $S_{\tau_{2k_1}}$, or $S_n$ and $S_{\tau_{2k+2}}$ can have overshoots with each other depending upon whether the horizontal/vertical displacement of one random walk exceeds that of the other, which can be computed by taking,

\begin{align*}
  \Delta_{2k+1}   =    - 1 - \big( S_{\tau_{2k+1}} - S_{\tau_{2k}} \big)   , \\ \\     \Delta_{2k+2}      =  1 - \big( S_{\tau_{2k+2}} - S_{\tau_{2k+1}} \big)  . 
\end{align*}

\noindent Denote the event that the triangular random walk at $\mathscr{T}_k$ is located in the bottom (resp. top) corner of the box, and it then at the top (resp. bottom) corner of the box at time $\mathscr{U}_k$. It is also possible that the random walk is located at the bottom (resp. top) corner of the box at $\mathscr{U}_k$, and is then located at the top (resp. bottom) corner of the box at time $\mathscr{U}_k$. We denote any one of these possible events with the A-type crossing event, $\mathscr{A}_k$.

\bigskip

\noindent In previous work for obtaining the scaling limit of interest,  {\color{blue}[4]}, the above representation for $S_n$ is used in the following several steps:

\begin{itemize}

\item[$\bullet$] \textit{(1). Determining the excursions of random walks over $\textbf{Z}^2$ which exceed a $\frac{4}{3}$ threshold}). One can upper bound $\textbf{P} \big( \mathscr{A}_k \big)$, from the superposition,

\begin{align*}
   \textbf{P} \big[ \mathcal{H}_{\mathscr{T}_k} <     c_1 k^{\frac{4}{3}}  \big] + \textbf{P} \big[  \mathscr{A}^{\mathrm{v}}_k                   \big|    \mathcal{H}_{\mathscr{T}_k} \geq c_1 k^{\frac{4}{3}}     \big]  , 
\end{align*}

\noindent of probabilities corresponding to horizontal, and to vertical, crossings, respectively, in which the first term is the chance of $\big\{ \mathcal{H}_{\mathscr{T}_k} <     c_1 k^{\frac{4}{3}} 
 \big\}$ occurring, while the second term is conditioned on the occurrence of $\big\{  \mathcal{H}_{\mathscr{T}_k} \geq     c_1 k^{\frac{4}{3}}       \big\}$ occurring. From this superposition, to bound the conditionally defined probability from above, observe,

 \begin{align*}
  \textbf{P} \big[  \mathscr{A}^{\mathrm{v}}_k                   \big|    \mathcal{H}_{\mathscr{T}_k} \geq c_1 k^{\frac{4}{3}}     \big]  \leq \frac{\big( c_1 \big)^{\prime}}{k^{\frac{4}{3}}}  , 
 \end{align*}

\bigskip

\noindent for a suitable $\big( c_1 \big)^{\prime}$.

\item[$\bullet$] \textit{(2). Relating the above probability for walks over $\textbf{Z}^2$ to the probability of an A type crossing event occurring is bound above with a strictly positive function of the box length}). With the quantities defined above, there exists suitable $c_1$ for which,

\begin{align*}
  \textbf{P} \big( \mathscr{A}_k \big) \leq \frac{c_1}{k^{\frac{4}{3}}} , 
\end{align*}

\noindent given $k$ large enough.

\item[$\bullet$]\textit{(3). There exists an upper bound for the horizontal, and vertical, excursions of the random walk over $\textbf{Z}^2$ which is exponentially small upper bounds}. There exists two constants for which,

\begin{align*}
  \textbf{P} \big(      \mathcal{W}_{\mathscr{T}_k}  < c_2 k^{\frac{4}{3}}  \big) \leq \mathrm{exp} \big[  - c_3 k^{\frac{4}{3}} \big]  , \\ \\ \textbf{P} \big(       \mathcal{H}_{\mathscr{T}_k} < c_2 k^{\frac{4}{3}}  \big) \leq \mathrm{exp} \big[  - c_3 k^{\frac{4}{3}} \big] , 
\end{align*}

\noindent in which $\big\{ \mathcal{W}_{\mathscr{T}_k}  < c_2 k^{\frac{4}{3}} \big\}$, and $\big\{ \mathcal{H}_{\mathscr{T}_k} < c_2 k^{\frac{4}{3}}  \big\}$, occurring are each bounded above by an exponentially decaying constant.

\item[$\bullet$] \textit{(4). The probability of $\big| S_k - \hat{S_k} \big|$ exceeding a threshold vanishes as n approaches infinity}. For any $\delta >0$, the probability,

\begin{align*}
\underset{ n \longrightarrow + \infty}{\mathrm{lim}}  \textbf{P} \bigg[      \text{ }             \underset{0 \leq k \leq n}{\mathrm{sup}}        \big| S_k - \hat{S_k} \big| \geq n^{\frac{1}{2} + \delta}                   \text{ }                \bigg]   , 
\end{align*}

\noindent vanishes.

\noindent To prove the claim, it suffices to demonstrate the existence of some threshold, $N$, for which,

\begin{align*}
 \textbf{P} \bigg[ \text{ }   \underset{0 \leq k \leq n}{\mathrm{sup}}   \big| S_k - \hat{S_k} \big| \geq N  \bigg]  = 0   , 
\end{align*}

\noindent beyond which it will also hold that,

\begin{align*}
   \textbf{P} \bigg[ \text{ }   \underset{0 \leq k \leq n}{\mathrm{sup}}   \big| S_k - \hat{S_k} \big| \geq      n^{\frac{1}{2} + \delta }    \bigg]  = 0   , 
\end{align*}

\noindent given some $m$, with,

\begin{align*}
  n \approx \frac{1}{N^{\frac{1}{2} + \delta}} . 
\end{align*}

\noindent The first statement above can be phrased in terms of the observation that,

\begin{align*}
\textbf{P} \bigg[   \text{ }     \underset{0 \leq  k \leq n}{\mathrm{sup}} \big|   S_k   -    \hat{S_k}  \big| \geq    N^{\prime}   \text{ }     \bigg]   = 0   , 
\end{align*}

\noindent for some strictly positive $N^{\prime}$, then necessarily,

\begin{align*}
\textbf{P} \bigg[   \text{ }     \underset{0 \leq  k \leq n}{\mathrm{sup}} \big|   S_k   -    \hat{S_k}  \big| \geq    N^{\prime\prime}   \text{ }     \bigg]   = 0   , 
\end{align*}

\noindent for some $N^{\prime\prime} > N^{\prime}$.

{\small \item[$\bullet$] \textit{(5). As n approaches infinity, probabilities dependent upon Brownian motion over $\textbf{Z}^2$ also vanish}. Fix some $\delta >0$ sufficiently small, and another constant $\epsilon$ which can also be taken sufficiently small. One has that,

\begin{align*}
 \underset{n \longrightarrow + \infty}{\mathrm{lim}}   \textbf{P} \bigg[         \text{ }   \underset{1 \leq k \leq n}{\mathrm{sup} }          \bigg| \text{ }  \frac{1}{n} \bigg[       \underset{1 \leq i \leq k }{\sum}   \textbf{1}_{\{\widetilde{S_i} \geq 0\} } -  \underset{1 \leq i \leq k }{\sum}   \textbf{1}_{\{S_i \geq n^{\frac{1}{2}} + \delta \} }   \bigg]  \bigg|  > \epsilon  \text{ }           \bigg]     , 
\end{align*}

\noindent and that,

\begin{align*}
  \underset{n \longrightarrow + \infty}{\mathrm{lim}}   \textbf{P} \bigg[         \text{ }   \underset{1 \leq k \leq n}{\mathrm{sup} }          \bigg| \text{ }  \frac{1}{n} \bigg[       \underset{1 \leq i \leq k }{\sum}   \textbf{1}_{\{\widetilde{S_i} < 0\} } -  \underset{1 \leq i \leq k }{\sum}   \textbf{1}_{\{S_i < - n^{\frac{1}{2}} + \delta \} }   \bigg]  \bigg|  > \epsilon  \text{ }     \bigg]       , 
\end{align*}

\noindent each vanish.

\item[$\bullet$] \textit{(6). The probability of the supremum of Brownian motions over $\textbf{Z}^2$ exceeds epsilon vanishes as n approaches infinity}). Fix some integer $k$. One has that,

\begin{align*}
  \underset{n \longrightarrow + \infty}{\mathrm{lim}}  \textbf{P} \bigg[ \text{ }      \underset{k \leq n}{\sup}    \text{ } \bigg| \text{ }     \frac{1}{n} \bigg[ \text{ }     \int_0^k \textbf{1}_{\{ B_s \geq 0 \}}   \mathrm{d} s - \underset{1 \leq i \leq k}{\sum}    \textbf{1}_{\{ S_i \geq n^{\frac{1}{2} 
 + \delta }\}} \text{ } \bigg]     \text{ } \bigg| > \epsilon  \text{ } \bigg]   , 
\end{align*}

\noindent and,

\begin{align*}
   \underset{n \longrightarrow + \infty}{\mathrm{lim}}  \textbf{P} \bigg[ \text{ }   \underset{k \leq n}{\sup}    \text{ } \bigg| \text{ }     \frac{1}{n} \bigg[ \text{ }       \int_0^k \textbf{1}_{\{ B_s < 0 \}} \mathrm{d} s - \underset{1 \leq i \leq k}{\sum}  \textbf{1}_{\{ S_i < n^{-\frac{1}{2} + \delta}\}}  \text{ } \bigg]          \text{ }     \bigg| > \epsilon \text{ } \bigg], 
\end{align*}

\noindent vanish. }

\item[$\bullet$] \textit{(8).} Introduce the stopping time,

\begin{align*}
   \vec{\eta_L} :=    \mathrm{inf} \big\{    n > 0 : S_n \geq L - 1            \big\} , 
\end{align*}

\noindent for the effective random walk, from which one can obtain the desired lower bound for $\textbf{P} \big( \eta_L \geq n \big)$ by observing that,

\begin{align*}
   \textbf{P} \big( \eta_L \geq n \big)   , 
\end{align*}

\noindent can be related to the joint probability measure,

  \begin{align*}  \textbf{P} \bigg[    \bigg\{ \big\{    \eta_{\infty} \geq n         \big\}      ,  \big\{    \vec{\eta_L} \geq n    \big\} \bigg\}                 \bigg] , \end{align*}

\noindent over the intersection of events,

{\small \begin{align*}  \bigg\{ \big\{    \eta_{\infty} \geq n         \big\}      \cap   \big\{    \vec{\eta_L} \geq n    \big\} \bigg\}  = \bigg\{ \big\{    \eta_{\infty} \geq n         \big\}      ,  \big\{    \vec{\eta_L} \geq n    \big\} \bigg\} . \end{align*} }

\item[$\bullet$] \textit{(9). Deducing that the supremum of the l-2 norm between the time re-weighted path and the normalizing constant exceeding an arbitrarily small parameter vanishes}. With previously defined quantities, one has that,

\begin{align*}
  \underset{t \longrightarrow + \infty}{\mathrm{lim}}     \textbf{P} \bigg[    \text{ }   \underset{ 0 \leq s \leq t}{\mathrm{sup}}     \bigg| \bigg| \frac{1}{t} \hat{\gamma_s} -  Z^{1,1}_{u , \alpha } \bigg| \bigg|_2 \geq \epsilon \text{ }  \bigg]   , 
\end{align*}

\noindent vanishes, for,

\begin{align*}
   Z^{1,1}_{u , \alpha } : = \textit{Scaling limit normalization for the random walk over } \textbf{Z}^2  .
\end{align*}

\item[$\bullet$] \textit{(10). Representation a transformation of the effective random walk over $\textbf{Z}^2$}. The expansion,

\begin{align*}
      \hat{S_n} :=S_n + \underset{j \geq 0}{\sum} \Delta_j \textbf{1}_{\{ \tau_j \leq n \}} , 
\end{align*}

\noindent corresponds to a transformation of the effective random walk representation, $S_n$, over $\textbf{Z}^2$, that is obtained with $\hat{\cdot}$. This transformation is used to conclude that the normalization of the scaling limit satisfies the previous item stated in \textit{(9)}. Explicitly,

{\small \begin{align*}
   \textbf{P} \bigg[  \text{ }  \underset{0 \leq k \leq n}{ \mathrm{sup} }      \big| S_k - \hat{S_k} \big| \geq n^{\frac{1}{3} + \delta}         \text{ }            \bigg]    
 \leq 0 \Longleftrightarrow     \textbf{P} \bigg[  \text{ }  \underset{0 \leq k \leq n}{ \mathrm{sup} }      \big| S_k - \hat{S_k} \big| \geq n^{\frac{1}{3} + \delta}         \text{ }            \bigg]  = 0  , 
\end{align*} }

\noindent implies,

{\small \begin{align*}
\underset{n \longrightarrow + \infty}{\mathrm{lim}}  \textbf{P} \bigg[ \text{ }      \underset{k \leq n}{\mathrm{sup} }  \big| \hat{S_k} - \sigma B_k \big| > n^{\frac{1}{3} + \delta }    \text{ }  \bigg]      , 
\end{align*} }

\noindent namely that there exists a probability space for which the supremum of the difference between $S_k$ and $\sigma B_k$ exceeds $n^{\frac{1}{3}}$ is exponentially small.

\end{itemize}

\subsection{Obtaining the closed form representation of the triangular Green's function from slab measures}

\noindent From remarks in  \textit{1.3}, the assumption that the triangular Green's function,

\begin{align*}
 G^r_{\textbf{T}} \big( \lambda \big)   , 
\end{align*}

\noindent has a decomposition of the form,

\begin{align*}
 G^r_{\textbf{T}} \big( \lambda \big)  =      G^r_{T \subsetneq \textbf{T}} \big( \lambda \big)     \textbf{1}_{\{ R = R_1 : R_1 \text{ } \textit{spans}\text{ }  T \} }           +  G^r_{T^{\prime} \subsetneq      \textbf{T}} \big( \lambda \big)     \textbf{1}_{\{ R = R_2 : R_2 \text{ } \textit{spans} \text{ } T^{\prime} \} } ,  \\   \tag{*}                  
\end{align*}

\noindent can be related to $\mathscr{K} \big( \lambda \big)$ through the following result.

\bigskip

\noindent \textbf{Lemma} \textit{1} (\textit{expressing the triangular Green's function in terms of the measure over the triangular effective walk, from arguments in {\color{blue}[4]}} for the Green's function for the square effective walk). Fix $n,n^{\prime}>0$. The Green's function, $ G^r_{\textbf{T}} \big( \lambda \big)$, with the decomposition $(*)$ provide above, takes the form,

\noindent under under the identification,

\begin{align*}
      G^r_{T \subsetneq \textbf{T}} \big( \lambda \big)     \textbf{1}_{\{ R = R_1 : R_1 \text{ } \textit{spans}\text{ }  T \} }           +  G^r_{T^{\prime} \subsetneq      \textbf{T}} \big( \lambda \big)     \textbf{1}_{\{ R = R_2 : R_2 \text{ } \textit{spans} \text{ } T^{\prime} \} }  ,
\end{align*}

\noindent coincides with $\mathscr{K} \big( \lambda \big)$ provided in \textbf{Definition} \textit{25}, under the closed form representation,

{\small \begin{align*}
    G^r_{\textbf{T}} \big( \lambda \big)   =            \frac{1}{18} \bigg[ 1 +   \mathrm{exp} \big[       - 3^{-1} \lambda            \big]           \bigg]  \bigg\{  \mathrm{exp} \bigg[        \bigg\{ \mathrm{log} \bigg[   \frac{ \mathrm{d}  \textbf{P}^{*}_{E_1}}{ \mathrm{d}  \textbf{P}_{E_1} +  \mathrm{d}  \textbf{P}_{E_2}}     \bigg] \bigg\}         \bigg]  +   \mathrm{exp} \bigg[         \mathrm{log} \bigg\{    \frac{ \mathrm{d}  \textbf{P}^{*}_{E_2}}{ \mathrm{d}  \textbf{P}_{E_1} +  \mathrm{d}  \textbf{P}_{E_2}}     \bigg\}           \bigg]       \bigg\}    \\   +    \frac{1}{6}       \bigg[ 1-  \mathrm{exp} \big[       - 3^{-1} \lambda            \big]           \bigg]  \bigg\{   \mathscr{E} \bigg\{      \mathrm{log} \bigg[   \frac{ \mathrm{d}  \textbf{P}^{*}_{E_1}}{ \mathrm{d}  \textbf{P}_{E_1} +  \mathrm{d}  \textbf{P}_{E_2}}     \bigg] \tau  - \lambda  \underset{1 \leq i \leq \tau}{\sum} \big| \widetilde{\mathscr{U}_i} \big|        \bigg\}       +   \mathscr{E} \bigg\{      \mathrm{log} \bigg[   \frac{ \mathrm{d}  \textbf{P}^{*}_{E_2}}{ \mathrm{d}  \textbf{P}_{E_1} +  \mathrm{d}  \textbf{P}_{E_2}}     \bigg] \tau^{\prime} \\  - \lambda  \underset{1 \leq i \leq \tau}{\sum} \big| \widetilde{\mathscr{V}_i} \big|        \bigg\}      \bigg\}    , 
\end{align*}  }       

\noindent implies,

{\small \begin{align*}
    \underset{n , n^{\prime} \longrightarrow + \infty}{\mathrm{lim}}  \mathscr{E} \bigg\{ \begin{bmatrix}
      M^{\lambda}_n      \\  M^{\lambda}_{n^{\prime}}
   \end{bmatrix} \bigg\}  \equiv \underset{n , n^{\prime} \longrightarrow + \infty}{\mathrm{lim}}  \mathscr{E} \bigg\{ \begin{bmatrix}
     \textit{Martingale along } E_1 \textit{ of the triangular effective random walk}      \\  \textit{Martingale along } E_2 \textit{ of the triangular effective random walk}   
   \end{bmatrix} \bigg\}   \\ \\ = \begin{bmatrix} 1 \\ 1  \end{bmatrix} : = \vec{1} .
\end{align*} }

\noindent \textit{Proof of Lemma 1}. To argue that $  \hat{\mathscr{K}} \big( \lambda \big)     = G^r_{T \subsetneq \textbf{T}} \big( \lambda \big)     \textbf{1}_{\{ R = R_1 : R_1 \text{ } \textit{spans}\text{ }  T \} }           +  G^r_{T^{\prime} \subsetneq      \textbf{T}} \big( \lambda \big)     \textbf{1}_{\{ R = R_2 : R_2 \text{ } \textit{spans} \text{ } T^{\prime} \} }$, which can subsequently be used to obtain the desired representation for the Green's function, we adapt the arguments of \textbf{Lemma} \textit{3.4} in {\color{blue}[4]}. We make use of the observation that,

\begin{align*}
  \underset{t \geq 1}{\sum}     \hat{\mathscr{K}} \big( \lambda , t  \big) =  G^r_{T \subsetneq \textbf{T}} \big( \lambda \big)     \textbf{1}_{\{ R = R_1 : R_1 \text{ } \textit{spans}\text{ }  T \} }           +  G^r_{T^{\prime} \subsetneq      \textbf{T}} \big( \lambda \big)     \textbf{1}_{\{ R = R_2 : R_2 \text{ } \textit{spans} \text{ } T^{\prime} \} }  , 
\end{align*}

\noindent we manipulate,

{\small \begin{align*}
  \hat{\mathscr{K}} \big( \lambda \big)  =  \mathscr{E}  \bigg\{      \textbf{1}_{\{ R = R_1 \} }                \mathrm{exp} \bigg[      \mathrm{d} \bigg[  \frac{\textbf{P}^{*}_{E_1}}{\textbf{P}_{E_1} + \textbf{P}_{E_2}} \bigg]        n      \textbf{1}_{\{    \widetilde{\mathscr{U}} \in \{ 0 , \cdots, R \}  : \widetilde{\mathscr{U}}_n = 0 ,  \underset{1 \leq i \leq n}{\sum} | \widetilde{\mathscr{U}}_i |         = t - n                                                  \}}                      \bigg]                               +   \textbf{1}_{\{ R = R_2 \} }       \mathrm{exp} \bigg[    \mathrm{d} \bigg[ \frac{\textbf{P}^{*}_{E_2}}{\textbf{P}_{E_1} + \textbf{P}_{E_2}}       \bigg] n^{\prime} \\ \times    \textbf{1}_{\{  \widetilde{\mathscr{V}} \in \{ 0 , \cdots, R \}  :   \widetilde{U}_n = 0 ,  \underset{1 \leq i \leq n}{\sum} | \widetilde{\mathscr{V}}_i  |                 =    t - n^{\prime}          \}}                      \bigg] \bigg]        \bigg\}          . 
\end{align*} }

\noindent which has the following decomposition under all times $t$,

\begin{align*}
  \hat{\mathscr{K}} \big( \lambda \big) =  \underset{t \geq 1}{\sum}     \hat{\mathscr{K}} \big( \lambda , t  \big)       .
  \end{align*}

\noindent To this end, write,

{\small \begin{align*}
  \underline{\underset{t \geq 1}{\sum}     \hat{\mathscr{K}} \big( \lambda , t  \big)}     =   \underset{t \geq 1}{\sum}     \mathscr{E}  \bigg\{      \textbf{1}_{\{ R = R_1 \} }                \mathrm{exp} \bigg[      \mathrm{d} \bigg[  \frac{\textbf{P}^{*}_{E_1}}{\textbf{P}_{E_1} + \textbf{P}_{E_2}} \bigg]        n     \textbf{1}_{\{    \widetilde{\mathscr{U}} \in \{ 0 , \cdots, R \}  : \widetilde{\mathscr{U}}_n = 0 ,  \underset{1 \leq i \leq n}{\sum} | \widetilde{\mathscr{U}}_i |         = t - n                                                  \}}                      \bigg]                               +   \textbf{1}_{\{ R = R_2 \} }     \\ \times     \mathrm{exp} \bigg[    \mathrm{d} \bigg[ \frac{\textbf{P}^{*}_{E_2}}{\textbf{P}_{E_1} + \textbf{P}_{E_2}}       \bigg]   n^{\prime} \textbf{1}_{\{  \widetilde{\mathscr{V}} \in \{ 0 , \cdots, R \}  :   \widetilde{U}_n = 0 ,  \underset{1 \leq i \leq n}{\sum} | \widetilde{\mathscr{V}}_i  |                 =    t - n^{\prime}          \}}                      \bigg]    \bigg]        \bigg\}     \\         \\      =    \underset{N \leq t}{\underset{t \geq 1}{\sum}}     \mathscr{E}  \bigg\{      \textbf{1}_{\{ R = R_1 \} }                \mathrm{exp} \bigg[      \mathrm{d} \bigg[  \frac{\textbf{P}^{*}_{E_1}}{\textbf{P}_{E_1} + \textbf{P}_{E_2}} \bigg]         n    \textbf{1}_{\{    \widetilde{\mathscr{U}} \in \{ 0 , \cdots, R \}  : \widetilde{U}_N = 0 ,  N+\underset{1 \leq i \leq N}{\sum} | \widetilde{\mathscr{U}}_i |         = t - n                                                  \}}                      \bigg]                               +   \textbf{1}_{\{ R = R_2 \} }     \\  \times     \mathrm{exp} \bigg[    \mathrm{d} \bigg[ \frac{\textbf{P}^{*}_{E_2}}{\textbf{P}_{E_1} + \textbf{P}_{E_2}}       \bigg]  n^{\prime}  \textbf{1}_{\{  \widetilde{\mathscr{V}} \in \{ 0 , \cdots, R \}  :   \widetilde{U}_N = 0 ,  N+ \underset{1 \leq i \leq N}{\sum} | \widetilde{\mathscr{V}}_i  |                 =    t - n^{\prime}          \}}                      \bigg]  \bigg]        \bigg\}                  \\ \\  =      \underset{1 \leq N \leq + \infty}{\sum}    \mathscr{E}  \bigg\{      \textbf{1}_{\{ R = R_1 \} }                \mathrm{exp} \bigg[      \mathrm{d} \bigg[  \frac{\textbf{P}^{*}_{E_1}}{\textbf{P}_{E_1} + \textbf{P}_{E_2}} \bigg]         n    \textbf{1}_{\{    \widetilde{\mathscr{U}} \in \{ 0 , \cdots, R \}  : \widetilde{U}_N = 0                                 \}}                      \bigg]                               +   \textbf{1}_{\{ R = R_2 \} }      \mathrm{exp} \bigg[    \mathrm{d} \bigg[ \frac{\textbf{P}^{*}_{E_2}}{\textbf{P}_{E_1} + \textbf{P}_{E_2}}       \bigg] n^{\prime} \end{align*}
  
  \begin{align*}   \times    \textbf{1}_{\{  \widetilde{\mathscr{V}} \in \{ 0 , \cdots, R \}  :   \widetilde{U}_N = 0    \}}                      \bigg]  \bigg]        \bigg\}                             \\ \\  =       \underset{1 \leq N \leq + \infty}{\sum}   \bigg[   \underset{1 \leq r \leq + \infty}{\underset{n^{\prime}_1 + \cdots + n^{\prime}_r = N^{\prime} }{\underset{n_1 + \cdots + n_r = N       }{\sum}}}  \bigg[ \underset{1 \leq i \leq r}{\prod}      \bigg[  \mathscr{E}  \bigg\{      \textbf{1}_{\{ R = R_1 \} }                \mathrm{exp} \bigg[      \mathrm{d} \bigg[  \frac{\textbf{P}^{*}_{E_1}}{\textbf{P}_{E_1} + \textbf{P}_{E_2}} \bigg]        n     \textbf{1}_{\{    \widetilde{\mathscr{U}} \in \{ 0 , \cdots, R \}  : \widetilde{U}_{n_i} = 0                                 \}}                      \bigg]                               +   \textbf{1}_{\{ R = R_2 \} }   \\   \times     \mathrm{exp} \bigg[    \mathrm{d} \bigg[ \frac{\textbf{P}^{*}_{E_2}}{\textbf{P}_{E_1} + \textbf{P}_{E_2}}       \bigg]  n^{\prime}   \textbf{1}_{\{  \widetilde{\mathscr{V}} \in \{ 0 , \cdots, R \}  :   \widetilde{U}_{n^{\prime}_i} = 0    \}}                      \bigg]   \bigg]        \bigg\}     \bigg]     \bigg]    \bigg]    \\ \\ =      \underset{1 \leq r \leq + \infty}{\sum}   \bigg[   \underset{1 \leq n^{\prime} \leq + \infty}{\underset{1 \leq n \leq + \infty}{\sum}}      \bigg[  \mathscr{E}  \bigg\{      \textbf{1}_{\{ R = R_1 \} }                \mathrm{exp} \bigg[      \mathrm{d} \bigg[  \frac{\textbf{P}^{*}_{E_1}}{\textbf{P}_{E_1} + \textbf{P}_{E_2}} \bigg]         n     \textbf{1}_{\{    \widetilde{\mathscr{U}} \in \{ 0 , \cdots, R \}  : \widetilde{U}_{n_i} = 0                                 \}}                      \bigg]                               +   \textbf{1}_{\{ R = R_2 \} }     \\  \times      \mathrm{exp} \bigg[    \mathrm{d} \bigg[ \frac{\textbf{P}^{*}_{E_2}}{\textbf{P}_{E_1} + \textbf{P}_{E_2}}       \bigg]   n^{\prime}     \textbf{1}_{\{  \widetilde{\mathscr{V}} \in \{ 0 , \cdots, R \}  :   \widetilde{U}_{n^{\prime}_i} = 0    \}}                      \bigg]  \bigg]        \bigg\}     \bigg]             \bigg]^r   \\ \\ =     \underset{1 \leq r \leq + \infty}{\sum}   \bigg[   {\underset{1 \leq n \leq + \infty}{\sum}}      \bigg[  \mathscr{E}  \bigg\{      \textbf{1}_{\{ R = R_1 \} }                \mathrm{exp} \bigg[      \mathrm{d} \bigg[  \frac{\textbf{P}^{*}_{E_1}}{\textbf{P}_{E_1} + \textbf{P}_{E_2}} \bigg]             \textbf{1}_{\{    \widetilde{\mathscr{U}} \in \{ 0 , \cdots, R \}  : \widetilde{U}_{n_i} = 0                                 \}}                      \bigg]                               +    \underset{1 \leq n^{\prime} \leq + \infty}{\sum}   \textbf{1}_{\{ R = R_2 \} }    \\  \times      \mathrm{exp} \bigg[    \mathrm{d} \bigg[ \frac{\textbf{P}^{*}_{E_2}}{\textbf{P}_{E_1} + \textbf{P}_{E_2}}       \bigg]        \textbf{1}_{\{  \widetilde{\mathscr{V}} \in \{ 0 , \cdots, R \}  :   \widetilde{U}_{n^{\prime}_i} = 0    \}}                      \bigg]   \eta_t \bigg]        \bigg\}     \bigg]             \bigg]^r         \end{align*}

  \begin{align*}   =        \underset{1 \leq r \leq + \infty}{\sum}   \bigg[   {\underset{1 \leq n \leq + \infty}{\sum}}      \bigg[  \mathscr{E}  \bigg\{      \textbf{1}_{\{ R = R_1 \} }                \mathrm{exp} \bigg[      \mathrm{d} \bigg[  \frac{\textbf{P}^{*}_{E_1}}{\textbf{P}_{E_1} + \textbf{P}_{E_2}} \bigg]        n      \textbf{1}_{\{    \widetilde{\mathscr{U}} \in \{ 0 , \cdots, R \}  : \widetilde{U}_{n_i} = 0                                 \}}                      \bigg]    \bigg\}^r                             +   \underset{1 \leq n^{\prime} \leq + \infty}{\sum} \mathscr{E} \bigg\{     \textbf{1}_{\{ R = R_2 \} }    \\ \times      \mathrm{exp} \bigg[    \mathrm{d} \bigg[ \frac{\textbf{P}^{*}_{E_2}}{\textbf{P}_{E_1} + \textbf{P}_{E_2}}       \bigg]    n^{\prime}    \textbf{1}_{\{  \widetilde{\mathscr{V}} \in \{ 0 , \cdots, R \}  :   \widetilde{U}_{n^{\prime}_i} = 0    \}}                      \bigg]  \bigg]        \bigg\}^r            \bigg]  \bigg]  \\ \\ =             \underline{G^r_{T \subsetneq \textbf{T}} \big( \lambda \big)     \textbf{1}_{\{ R = R_1 : R_1 \text{ } \textit{spans}\text{ }  T \} }           +  G^r_{T^{\prime} \subsetneq      \textbf{T}} \big( \lambda \big)     \textbf{1}_{\{ R = R_2 : R_2 \text{ } \textit{spans} \text{ } T^{\prime} \} }}           .
  \end{align*} }

\noindent To obtain the desired representation for the Green's function, introduce,

\begin{align*}
 \vec{1} : =     \begin{bmatrix} \mathrm{log} \big[   \frac{ \mathrm{d}  \textbf{P}^{*}_{E_1}}{ \mathrm{d}  \textbf{P}_{E_1} +  \mathrm{d}  \textbf{P}_{E_2}}     \big] \\ \mathrm{log} \big[   \frac{ \mathrm{d}  \textbf{P}^{*}_{E_2}}{ \mathrm{d}  \textbf{P}_{E_1} +  \mathrm{d}  \textbf{P}_{E_2}}     \big] 
 \end{bmatrix} -   \begin{bmatrix} \lambda \\ \lambda  
 \end{bmatrix} ,
\end{align*}

  \noindent for some $\lambda > 0$. 

\bigskip

\noindent The above vector, given the magnitudes of $\mathrm{d}  \textbf{P}^{*}_{E_1}$ and $\mathrm{d}  \textbf{P}^{*}_{E_2}$ before a normalization with $ \mathrm{d}  \textbf{P}_{E_1} +  \mathrm{d}  \textbf{P}_{E_2}$, necessarily implies that the sets $\widetilde{\mathscr{U}}$, and $\widetilde{\mathscr{V}}$ appearing in,

{\small \begin{align*}
        \underset{1 \leq r \leq + \infty}{\sum} \bigg\{    {\underset{1 \leq n \leq + \infty}{\sum}}      \mathscr{E}  \bigg\{      \textbf{1}_{\{ R = R_1 \} }                \mathrm{exp} \bigg[      \mathrm{d} \bigg[  \frac{\textbf{P}^{*}_{E_1}}{\textbf{P}_{E_1} + \textbf{P}_{E_2}} \bigg]        n      \textbf{1}_{\{    \widetilde{\mathscr{U}} \in \{ 0 , \cdots, R \}  : \widetilde{U}_{n_i} = 0                                 \}}                      \bigg]    \bigg\}^r  \bigg\}         , 
\end{align*} }

\noindent and in,

{\small \begin{align*}
 \underset{1 \leq r \leq + \infty}{\sum} \bigg\{  \underset{1 \leq n^{\prime} \leq + \infty}{\sum}   \mathscr{E} \bigg\{     \textbf{1}_{\{ R = R_2 \} }        \mathrm{exp} \bigg[    \mathrm{d} \bigg[ \frac{\textbf{P}^{*}_{E_2}}{\textbf{P}_{E_1} + \textbf{P}_{E_2}}       \bigg]    n^{\prime}    \textbf{1}_{\{  \widetilde{\mathscr{V}} \in \{ 0 , \cdots, R \}  :   \widetilde{U}_{n^{\prime}_i} = 0    \}}                      \bigg]       \bigg\}^r     \bigg\}    ,
\end{align*} } 

\noindent which are used to express the $r$-th power of the Green's function satisfy,

\[ \left\{\!\begin{array}{ll@{}>{{}}l} 
   \mathrm{log} \bigg[  \mathscr{E} \big[ - \lambda \big| \widetilde{\mathscr{U}} \big|  \big]   \bigg] =  \lambda - \mathrm{log} \bigg[   \frac{ \mathrm{d}  \textbf{P}^{*}_{E_1}}{ \mathrm{d}  \textbf{P}_{E_1} +  \mathrm{d}  \textbf{P}_{E_2}}     \bigg] ,  \\ \\      \mathrm{log} \bigg[  \mathscr{E} \big[ - \lambda \big| \widetilde{\mathscr{V}}  \big| \big]   \bigg] =  \lambda -  \mathrm{log} \bigg[   \frac{ \mathrm{d}  \textbf{P}^{*}_{E_2}}{ \mathrm{d}  \textbf{P}_{E_1} +  \mathrm{d}  \textbf{P}_{E_2}}     \bigg] ,
\end{array}\right. 
\]

\noindent The above set of relations is expressed more compactly with,

\begin{align*}
    \mathrm{log} \bigg[ \mathscr{E} \bigg\{   - \lambda  \begin{bmatrix} \big| \widetilde{\mathscr{U}} \big|  \\ \big| \widetilde{\mathscr{V}} \big|  \end{bmatrix} \bigg\}   \bigg]     = - \vec{1} .
\end{align*}

  \bigskip

  \noindent To argue that the Green's function is composed of contributions which take the form,

{\small \begin{align*}
     \frac{  \mathscr{E} \bigg\{ \mathrm{exp} \bigg[        \bigg\{ \mathrm{log} \bigg[   \frac{ \mathrm{d}  \textbf{P}^{*}_{E_1}}{ \mathrm{d}  \textbf{P}_{E_1} +  \mathrm{d}  \textbf{P}_{E_2}}     \bigg] - \lambda \bigg\}         \bigg]           \bigg\}             \textbf{1}_{ \bigg\{ \widetilde{\mathscr{U}} : \big\{ \widetilde{\mathscr{U}}_n = 0  \big\} , \big\{ \widetilde{\mathscr{U}}_0 = n \big\} \bigg\}  }      }{  \mathscr{E} \bigg\{ \mathrm{exp} \bigg[  \bigg\{ \mathrm{log} \bigg[   \frac{ \mathrm{d}  \textbf{P}^{*}_{E_1}}{ \mathrm{d}  \textbf{P}_{E_1} +  \mathrm{d}  \textbf{P}_{E_2}}     \bigg] - \lambda \bigg\}           \bigg]   \bigg\}   \textbf{1}_{ \bigg\{ \widetilde{\mathscr{U}} :  \big\{ \widetilde{\mathscr{U}}_0 = n \big\} \bigg\}  }                        }        , 
\end{align*}     }

  \noindent and also,

{\small \begin{align*}
         \frac{  \mathscr{E} \bigg\{ \mathrm{exp} \bigg[        \bigg\{ \mathrm{log} \bigg[   \frac{ \mathrm{d}  \textbf{P}^{*}_{E_2}}{ \mathrm{d}  \textbf{P}_{E_1} +  \mathrm{d}  \textbf{P}_{E_2}}     \bigg] - \lambda \bigg\}         \bigg]           \bigg\}             \textbf{1}_{ \bigg\{ \widetilde{\mathscr{V}} : \big\{ \widetilde{\mathscr{V}}_n = 0  \big\} , \big\{ \widetilde{\mathscr{V}}_0 = n \big\} \bigg\}  }      }{  \mathscr{E} \bigg\{ \mathrm{exp} \bigg[  \bigg\{ \mathrm{log} \bigg[   \frac{ \mathrm{d}  \textbf{P}^{*}_{E_2}}{ \mathrm{d}  \textbf{P}_{E_1} +  \mathrm{d}  \textbf{P}_{E_2}}     \bigg] - \lambda \bigg\}           \bigg]   \bigg\}   \textbf{1}_{ \bigg\{ \widetilde{\mathscr{V}} :  \big\{ \widetilde{\mathscr{V}}_0 = n \big\} \bigg\}  }                        }          , 
\end{align*}   }

\noindent for small times $\tau \neq \tau^{\prime} > 0$, it suffices to approximate,

{\small \begin{align*}
       \mathscr{E} \bigg\{ \mathrm{exp} \bigg[        \bigg\{ \mathrm{log} \bigg[   \frac{ \mathrm{d}  \textbf{P}^{*}_{E_1}}{ \mathrm{d}  \textbf{P}_{E_1} +  \mathrm{d}  \textbf{P}_{E_2}}     \bigg] - \lambda \bigg\}         \bigg]           \bigg\}                        , 
\end{align*}     }

\noindent and,

{\small \begin{align*}
     \mathscr{E} \bigg\{ \mathrm{exp} \bigg[        \bigg\{ \mathrm{log} \bigg[   \frac{ \mathrm{d}  \textbf{P}^{*}_{E_2}}{ \mathrm{d}  \textbf{P}_{E_1} +  \mathrm{d}  \textbf{P}_{E_2}}     \bigg] - \lambda \bigg\}         \bigg]           \bigg\}          . 
\end{align*}     }

\noindent Fix the times,

\begin{align*}
   \tau : =   \underset{\textit{times $t$}}{\mathrm{sup}}   \big\{ t :  \textit{the triangular effective random walk hits the origin about the first degree of} \\ \textit{ freedom} \big\}    , \\ \\   \tau^{\prime} : =    \underset{\textit{times $t$}}{\mathrm{sup}}   \big\{ t : \textit{the triangular effective random walk hits the origin about the second degree} \\ \textit{ of freedom}     \big\}     .
\end{align*}

\noindent To bound the first ratio of expected values, write,

{\tiny \begin{align*}
    \frac{  \mathscr{E} \bigg\{ \mathrm{exp} \bigg[        \bigg\{ \mathrm{log} \bigg[   \frac{ \mathrm{d}  \textbf{P}^{*}_{E_1}}{ \mathrm{d}  \textbf{P}_{E_1} +  \mathrm{d}  \textbf{P}_{E_2}}     \bigg] \tau  - \lambda \bigg\}         \bigg]           \bigg\}             \textbf{1}_{ \bigg\{ \widetilde{\mathscr{U}} : \big\{ \widetilde{\mathscr{U}}_n = 0  \big\} , \big\{ \widetilde{\mathscr{U}}_0 = n \big\} \bigg\}  }      }{  \mathscr{E} \bigg\{ \mathrm{exp} \bigg[  \bigg\{ \mathrm{log} \bigg[   \frac{ \mathrm{d}  \textbf{P}^{*}_{E_1}}{ \mathrm{d}  \textbf{P}_{E_1} +  \mathrm{d}  \textbf{P}_{E_2}}     \bigg]  \tau - \lambda \bigg\}           \bigg]   \bigg\}   \textbf{1}_{ \bigg\{ \widetilde{\mathscr{U}} :  \big\{ \widetilde{\mathscr{U}}_0 = n \big\} \bigg\}  }                        }   =        \mathscr{E} \bigg\{ \mathrm{exp} \bigg[        \bigg\{ \mathrm{log} \bigg[   \frac{ \mathrm{d}  \textbf{P}^{*}_{E_1}}{ \mathrm{d}  \textbf{P}_{E_1} +  \mathrm{d}  \textbf{P}_{E_2}}     \bigg] \tau  - \lambda \bigg\}         \bigg]           \bigg\}            \frac{\textbf{1}_{ \bigg\{ {\widetilde{\mathscr{U}} : \big\{ \widetilde{\mathscr{U}}_n = 0  \big\} , \big\{ \widetilde{\mathscr{U}}_0 = n \big\} } \bigg\} }}{   \textbf{1}_{ \bigg\{ {\widetilde{\mathscr{U}} :  \big\{ \widetilde{\mathscr{U}}_0 = n \big\}   } \bigg\} }}               ,
\end{align*}}

\noindent from which one observes, 

{\small
\begin{align*}
   \frac{\textbf{1}_{\big\{   \widetilde{\mathscr{U}} : \big\{ \widetilde{\mathscr{U}_n} = 0 \big\} \cup \big\{ \widetilde{\mathscr{U}_0} = n \big\}         \big\}}}{\textbf{1}_{ \big\{ {\widetilde{\mathscr{U}} :  \big\{ \widetilde{\mathscr{U}}_0 = n \big\}   } \big\} }}  \propto \textbf{P}^*_{E_1} \big[      \widetilde{\mathscr{U}}
   :   \big\{ \widetilde{U}_0 = n \big\} , \big\{ \widetilde{\mathscr{U}_0}  \neq n \big\}                                                  \big]       , 
\end{align*}
}

\noindent can be related to the slab measure through the ratio,

{\small \begin{align*}
   \frac{ \hat{\mathscr{L}_R}  }{\textbf{P}^*_{E_1} \bigg[      \big\{ \widetilde{\mathscr{U}}
   :   \big\{ \widetilde{U}_0 = n \big\} , \big\{ \widetilde{\mathscr{U}_0}  \neq n \big\}              \big\}   ,   \bigg\{  \mathscr{T}_1        > \mathscr{L}_1 - \underset{1 \leq i \leq \gamma_{\mathscr{L}_1}}{\sum} \mathscr{T}_i                \bigg\}                                   \bigg]  }    \\ \\  =   \frac{  \mathscr{E}  \bigg\{        \textbf{1}_{\{ R = R_{\gamma_{\mathscr{L}_1}} \} }       \mathrm{exp} \bigg[    \mathrm{d} \bigg[ \frac{\textbf{P}^{*}_{E_1}}{\textbf{P}_{E_1} + \textbf{P}_{E_2}}       \bigg]    n^{\prime}      \textbf{1}_{\{  \widetilde{\mathscr{V}} \in \{ 0 , \cdots, R \}  :   0 < \widetilde{U}_n <  R ,  \underset{1 \leq i \leq n}{\sum} | \widetilde{\mathscr{V}}_i  |                 =    t - n^{\prime}          \}}                      \bigg]         \bigg\}}{\textbf{P}^*_{E_1} \bigg[      \big\{ \widetilde{\mathscr{U}}
   :   \big\{ \widetilde{U}_0 = n \big\} , \big\{ \widetilde{\mathscr{U}_0}  \neq n \big\}              \big\}   ,   \bigg\{  \mathscr{T}_1        > \mathscr{L}_1 - \underset{1 \leq i \leq \gamma_{\mathscr{L}_1}}{\sum} \mathscr{T}_i                \bigg\}                                   \bigg] }     . 
\end{align*} }

\bigskip

\noindent Similarly, write,

{\tiny \begin{align*}
    \frac{  \mathscr{E} \bigg\{ \mathrm{exp} \bigg[        \bigg\{ \mathrm{log} \bigg[   \frac{ \mathrm{d}  \textbf{P}^{*}_{E_2}}{ \mathrm{d}  \textbf{P}_{E_1} +  \mathrm{d}  \textbf{P}_{E_2}}     \bigg] \tau^{\prime}  - \lambda \bigg\}         \bigg]           \bigg\}             \textbf{1}_{ \bigg\{ \widetilde{\mathscr{V}} : \big\{ \widetilde{\mathscr{V}}_n = 0  \big\} , \big\{ \widetilde{\mathscr{V}}_0 = n \big\} \bigg\}  }      }{  \mathscr{E} \bigg\{ \mathrm{exp} \bigg[  \bigg\{ \mathrm{log} \bigg[   \frac{ \mathrm{d}  \textbf{P}^{*}_{E_2}}{ \mathrm{d}  \textbf{P}_{E_1} +  \mathrm{d}  \textbf{P}_{E_2}}     \bigg] \tau^{\prime}  - \lambda \bigg\}           \bigg]   \bigg\}   \textbf{1}_{ \bigg\{ \widetilde{\mathscr{V}} :  \big\{ \widetilde{\mathscr{V}}_0 = n \big\} \bigg\}  }                        }      =     \mathscr{E} \bigg\{ \mathrm{exp} \bigg[        \bigg\{ \mathrm{log} \bigg[   \frac{ \mathrm{d}  \textbf{P}^{*}_{E_1}}{ \mathrm{d}  \textbf{P}_{E_1} +  \mathrm{d}  \textbf{P}_{E_2}}     \bigg] \tau^{\prime}  - \lambda \bigg\}         \bigg]           \bigg\}            \frac{\textbf{1}_{ \bigg\{ {\widetilde{\mathscr{V}} : \big\{ \widetilde{\mathscr{V}}_n = 0  \big\} , \big\{ \widetilde{\mathscr{V}}_0 = n \big\} } \bigg\} }}{   \textbf{1}_{ \bigg\{ {\widetilde{\mathscr{V}} :  \big\{ \widetilde{\mathscr{V}}_0 = n \big\}   } \bigg\} }}             ,
\end{align*}}

\noindent from which one observes,

{\small
\begin{align*}
  \frac{\textbf{1}_{\big\{   \widetilde{\mathscr{V}} : \big\{ \widetilde{\mathscr{V}_n} = 0 \big\} \cup \big\{ \widetilde{\mathscr{V}_0} = n \big\}         \big\}}}{\textbf{1}_{ \big\{ {\widetilde{\mathscr{V}} :  \big\{ \widetilde{\mathscr{V}}_0 = n \big\}   } \big\} }}  \propto \textbf{P}^*_{E_2} \big[      \widetilde{\mathscr{V}}
   :   \big\{ \widetilde{V}_0 = n \big\} , \big\{ \widetilde{\mathscr{V}_0}  \neq n \big\}                                                  \big]     , 
\end{align*}
}

\noindent can be related to the slab measure through the ratio,

{\small \begin{align*}
       \frac{\hat{\mathscr{L}_R} }{\textbf{P}^*_{E_2} \bigg[      \big\{ \widetilde{\mathscr{V}}
   :   \big\{ \widetilde{V}_0 = n \big\} , \big\{ \widetilde{\mathscr{V}_0}  \neq n^{\prime} \big\}   \big\} , \bigg\{  \mathscr{T}_2        > \mathscr{L}_2 - \underset{1 \leq i \leq \gamma_{\mathscr{L}_2}}{\sum} \mathscr{T}_i                \bigg\}                                                \bigg]  }  \\ \\  =  \frac{  \mathscr{E}  \bigg\{        \textbf{1}_{\{ R = R_{\gamma_{\mathscr{L}_2}} \} }       \mathrm{exp} \bigg[    \mathrm{d} \bigg[ \frac{\textbf{P}^{*}_{E_2}}{\textbf{P}_{E_1} + \textbf{P}_{E_2}}       \bigg]    n^{\prime}      \textbf{1}_{\{  \widetilde{\mathscr{V}} \in \{ 0 , \cdots, R \}  :   0 < \widetilde{U}_n <  R ,  \underset{1 \leq i \leq n}{\sum} | \widetilde{\mathscr{V}}_i  |                 =    t - n^{\prime}          \}}                      \bigg]         \bigg\}  }{\textbf{P}^*_{E_2} \bigg[      \big\{ \widetilde{\mathscr{V}}
   :   \big\{ \widetilde{V}_0 = n \big\} , \big\{ \widetilde{\mathscr{V}_0}  \neq n^{\prime} \big\}   \big\} , \bigg\{  \mathscr{T}_2        > \mathscr{L}_2 - \underset{1 \leq i \leq \gamma_{\mathscr{L}_2}}{\sum} \mathscr{T}_i                \bigg\}                                                \bigg] }     . 
\end{align*} }

\noindent The Green's function can be expressed from the spanning set,

\begin{align*}
  \mathrm{span} \mathscr{B}  
\end{align*}

\noindent associated with the basis,

{\small \begin{align*}
\bigg\{   \mathrm{exp} \bigg[  \bigg\{ \mathrm{log} \bigg[   \frac{ \mathrm{d}  \textbf{P}^{*}_{E_1}}{ \mathrm{d}  \textbf{P}_{E_1} +  \mathrm{d}  \textbf{P}_{E_2}}     \bigg] \tau  - \lambda \bigg\} \bigg]     ,   \mathrm{exp} \bigg[  \bigg\{ \mathrm{log} \bigg[   \frac{ \mathrm{d}  \textbf{P}^{*}_{E_2}}{ \mathrm{d}  \textbf{P}_{E_1} +  \mathrm{d}  \textbf{P}_{E_2}}     \bigg] \tau^{\prime}  - \lambda \bigg\} \bigg]   ,  \mathscr{E} \bigg\{     \textbf{1}_{\big\{ \widetilde{\mathscr{U}} \neq \widetilde{\mathscr{V}} :  \big\{ \widetilde{\mathscr{U}} \big|_{\textit{time} \text{ } 0} =   \widetilde{\mathscr{V}} \big|_{\textit{time} \text{ } 0}   = 0 \big\}    \big\} }             \bigg\}  \\ ,      \mathscr{E} \bigg\{        \textbf{1}_{\big\{ \widetilde{\mathscr{U}} \neq \widetilde{\mathscr{V}} :  \big\{ \widetilde{\mathscr{U}} \big|_{\textit{time} \text{ } 0} =   \widetilde{\mathscr{V}} \big|_{\textit{time} \text{ } 0}   = 0 \big\} , \big\{ \widetilde{\mathscr{U}}_n = \widetilde{\mathscr{V}}_{n^{\prime}} = 0             \big\}   \big\} }                 \bigg\}         \bigg\}   , 
\end{align*} }

\noindent $\mathscr{B}$. Before performing the following computations with martingales, observe that the expected number of excursions of the effective triangular walk,

{\small \begin{align*}
    {\small\mathscr{E} \bigg\{ \mathrm{exp} \bigg[        \bigg\{ \mathrm{log} \bigg[   \frac{ \mathrm{d}  \textbf{P}^{*}_{E_1}}{ \mathrm{d}  \textbf{P}_{E_1} +  \mathrm{d}  \textbf{P}_{E_2}}     \bigg] \tau  - \lambda \bigg\}         \bigg]           \bigg\}        \bigg\}  }          , \\ \\     {\small\mathscr{E} \bigg\{ \mathrm{exp} \bigg[        \bigg\{ \mathrm{log} \bigg[   \frac{ \mathrm{d}  \textbf{P}^{*}_{E_2}}{ \mathrm{d}  \textbf{P}_{E_1} +  \mathrm{d}  \textbf{P}_{E_2}}     \bigg] \tau^{\prime}  - \lambda \bigg\}         \bigg]           \bigg\}        \bigg\}  }            , 
\end{align*} }

\noindent can be related to another probability measure,

\begin{align*}
 \widetilde{\mathscr{P}_0 }   : = \textit{Probabilty measure over the triangular effective random walk for which } {\mathscr{P}}  \big[ \forall \big( x_1, x_2 \big)  \in \big( \mathscr{U} , \mathscr{V} \big) :     \big\{ \mathscr{U} \\ \cap \mathscr{U}_0  \neq \emptyset \big\}  , \big\{ \mathscr{V} \cap \mathscr{V}_0 \neq \emptyset \big\} ,         \big\{    \mathscr{U}_0 = x_1        \big\}   , \big\{  \mathscr{V}_0 = x_2             \big\}      \big]   = \mathscr{E} \bigg\{  - \lambda \begin{bmatrix}
 \big| \mathscr{U}_1 \big|     \\ \big| \mathscr{V}_1 \big| 
 \end{bmatrix} \bigg\}     \widetilde{\mathscr{P}} \big[   \forall \big( \widetilde{x_1}, \widetilde{x_2} \big)  \in \big( \widetilde{\mathscr{U}} , \widetilde{\mathscr{V}} \big) :     \big\{ \widetilde{\mathscr{U}} \\  \cap \widetilde{\mathscr{U}_0}  \neq \emptyset \big\}  , \big\{ \widetilde{\mathscr{V}} \cap \widetilde{\mathscr{V}_0} \neq \emptyset \big\} ,       \big\{    \widetilde{\mathscr{U}_0} = \widetilde{x_1}        \big\}   , \big\{  \widetilde{\mathscr{V}_0} = \widetilde{x_2}             \big\}     \big]  , 
\end{align*}

\noindent is an iid draw from the random sequence

\begin{align*}
      \big\{ \widetilde{\mathscr{P}}_i \big\}_{i \geq 0 }        , 
\end{align*}

\noindent where,

\begin{align*}
   \underset{i >  0}{\prod}  \big\{  \widetilde{\mathscr{P}_0}  \cap   \widetilde{\mathscr{P}_i} \big\}   =     \widetilde{\mathscr{P}_0}  \cap  \bigg\{      \underset{i >  0}{\prod} \widetilde{\mathscr{P}_i}   \bigg\}         . 
\end{align*}

\noindent The sequences of expected values,

{\small \begin{align*}
    \bigg\{ {\small\mathscr{E} \bigg\{ \mathrm{exp} \bigg[        \bigg\{ \mathrm{log} \bigg[   \frac{ \mathrm{d}  \textbf{P}^{*}_{E_1}}{ \mathrm{d}  \textbf{P}_{E_1} +  \mathrm{d}  \textbf{P}_{E_2}}     \bigg] \tau  - \lambda \bigg\}         \bigg]           \bigg\}             \textbf{1}_{ \bigg\{ \widetilde{\mathscr{U}} : \big\{ \widetilde{\mathscr{U}}_n = 0  \big\} , \big\{ \widetilde{\mathscr{U}}_0 = n \big\} \bigg\}  }}    \bigg\}_{n \geq 0 } , \\ \\ \bigg\{   {\small \mathscr{E} \bigg\{ \mathrm{exp} \bigg[        \bigg\{ \mathrm{log} \bigg[   \frac{ \mathrm{d}  \textbf{P}^{*}_{E_2}}{ \mathrm{d}  \textbf{P}_{E_1} +  \mathrm{d}  \textbf{P}_{E_2}}     \bigg] \tau^{\prime}  - \lambda \bigg\}         \bigg]           \bigg\}             \textbf{1}_{ \bigg\{ \widetilde{\mathscr{V}} : \big\{ \widetilde{\mathscr{V}}_{n^{\prime}} = 0  \big\} , \big\{ \widetilde{\mathscr{V}}_0 = n^{\prime} \big\} \bigg\}  }}           \bigg\}_{n^{\prime} \geq 0 } . 
\end{align*}     } 

\noindent can be straightforwardly used to construct the vector,

{\small \begin{align*}
\begin{bmatrix}   \bigg\{ {\small\mathscr{E} \bigg\{ \mathrm{exp} \bigg[        \bigg\{ \mathrm{log} \bigg[   \frac{ \mathrm{d}  \textbf{P}^{*}_{E_1}}{ \mathrm{d}  \textbf{P}_{E_1} +  \mathrm{d}  \textbf{P}_{E_2}}     \bigg] \tau  - \lambda \bigg\}         \bigg]           \bigg\}             \textbf{1}_{ \bigg\{ \widetilde{\mathscr{U}} : \big\{ \widetilde{\mathscr{U}}_n = 0  \big\} , \big\{ \widetilde{\mathscr{U}}_0 = n \big\} \bigg\}  }}    \bigg\}_{n \geq 0 }  \\  \bigg\{   {\small \mathscr{E} \bigg\{ \mathrm{exp} \bigg[        \bigg\{ \mathrm{log} \bigg[   \frac{ \mathrm{d}  \textbf{P}^{*}_{E_2}}{ \mathrm{d}  \textbf{P}_{E_1} +  \mathrm{d}  \textbf{P}_{E_2}}     \bigg] \tau^{\prime}  - \lambda \bigg\}         \bigg]           \bigg\}             \textbf{1}_{ \bigg\{ \widetilde{\mathscr{V}} : \big\{ \widetilde{\mathscr{V}}_{n^{\prime}} = 0  \big\} , \big\{ \widetilde{\mathscr{V}}_0 = n^{\prime} \big\} \bigg\}  }}           \bigg\}_{n^{\prime} \geq 0 }   \end{bmatrix} \end{align*}

\begin{align*} = \begin{bmatrix}  {\small\mathscr{E} \bigg\{ \mathrm{exp} \bigg[        \bigg\{ \mathrm{log} \bigg[   \frac{ \mathrm{d}  \textbf{P}^{*}_{E_1}}{ \mathrm{d}  \textbf{P}_{E_1} +  \mathrm{d}  \textbf{P}_{E_2}}     \bigg] \tau  - \lambda \bigg\}         \bigg]           \bigg\}             \textbf{1}_{ \bigg\{ \widetilde{\mathscr{U}} : \big\{ \widetilde{\mathscr{U}}_n = 0  \big\} , \big\{ \widetilde{\mathscr{U}}_0 = n \big\} \bigg\}  }}     \\     {\small \mathscr{E} \bigg\{ \mathrm{exp} \bigg[        \bigg\{ \mathrm{log} \bigg[   \frac{ \mathrm{d}  \textbf{P}^{*}_{E_2}}{ \mathrm{d}  \textbf{P}_{E_1} +  \mathrm{d}  \textbf{P}_{E_2}}     \bigg] \tau^{\prime}  - \lambda \bigg\}         \bigg]           \bigg\}             \textbf{1}_{ \bigg\{ \widetilde{\mathscr{V}} : \big\{ \widetilde{\mathscr{V}}_{n^{\prime}} = 0  \big\} , \big\{ \widetilde{\mathscr{V}}_0 = n^{\prime} \big\} \bigg\}  }}            \end{bmatrix}_{ \{ n \neq  n^{\prime} \geq 0 \}  }  \\ \\ \overset{(\textit{linearity of } \mathscr{E} [ \cdot ])}{=} \mathscr{E}   \begin{bmatrix}  {\small  \mathrm{exp} \bigg[        \bigg\{ \mathrm{log} \bigg[   \frac{ \mathrm{d}  \textbf{P}^{*}_{E_1}}{ \mathrm{d}  \textbf{P}_{E_1} +  \mathrm{d}  \textbf{P}_{E_2}}     \bigg] \tau  - \lambda \bigg\}         \bigg]                      \textbf{1}_{ \bigg\{ \widetilde{\mathscr{U}} : \big\{ \widetilde{\mathscr{U}}_n = 0  \big\} , \big\{ \widetilde{\mathscr{U}}_0 = n \big\} \bigg\}  }}     \\     {\small   \mathrm{exp} \bigg[        \bigg\{ \mathrm{log} \bigg[   \frac{ \mathrm{d}  \textbf{P}^{*}_{E_2}}{ \mathrm{d}  \textbf{P}_{E_1} +  \mathrm{d}  \textbf{P}_{E_2}}     \bigg] \tau^{\prime}  - \lambda \bigg\}         \bigg]                       \textbf{1}_{ \bigg\{ \widetilde{\mathscr{V}} : \big\{ \widetilde{\mathscr{V}}_{n^{\prime}} = 0  \big\} , \big\{ \widetilde{\mathscr{V}}_0 = n^{\prime} \big\} \bigg\}  }}            \end{bmatrix}_{ \{ n \neq  n^{\prime} \geq 0 \} }  ,
\end{align*}}

\noindent of martingales from the observation that,

{\small \begin{align*}
   \mathscr{E} \bigg\{ \mathrm{exp} \bigg[        \bigg\{ \mathrm{log} \bigg[   \frac{ \mathrm{d}  \textbf{P}^{*}_{E_1}}{ \mathrm{d}  \textbf{P}_{E_1} +  \mathrm{d}  \textbf{P}_{E_2}}     \bigg] \tau  - \lambda \bigg\}         \bigg]    \textbf{1}_{ \bigg\{ \widetilde{\mathscr{U}} : \big\{ \widetilde{\mathscr{U}}_{N^{\prime}} = 0  \big\} , \big\{ \widetilde{\mathscr{U}}_0 = N^{\prime} \big\} \bigg\}  }        \bigg|   \big\{   \textit{Steps of the triangular effective random walk} \\ \textit{between times $0$ and $N^{\prime}-1$ along $E_1$}   \big\}       \bigg\}  \\ \\ =    \mathscr{E} \bigg\{ \mathrm{exp} \bigg[        \bigg\{ \mathrm{log} \bigg[   \frac{ \mathrm{d}  \textbf{P}^{*}_{E_1}}{ \mathrm{d}  \textbf{P}_{E_1} +  \mathrm{d}  \textbf{P}_{E_2}}     \bigg] \tau  - \lambda \bigg\}         \bigg]    \textbf{1}_{ \bigg\{ \widetilde{\mathscr{U}} : \big\{ \widetilde{\mathscr{U}}_{N^{\prime}+1} = 0  \big\} , \big\{ \widetilde{\mathscr{U}}_0 = N^{\prime}+1  \big\} \bigg\}  }           \bigg\}      ,
\end{align*} }

\noindent implies that the Green's function along $E_1$ takes the form,

{\small \begin{align*}
   G_{E_1} \big( \lambda \big) =   \frac{1}{18} \bigg[ 1 +   \mathrm{exp} \big[       - 3^{-1} \lambda            \big]           \bigg]   \mathrm{exp} \bigg[        \bigg\{ \mathrm{log} \bigg[   \frac{ \mathrm{d}  \textbf{P}^{*}_{E_1}}{ \mathrm{d}  \textbf{P}_{E_1} +  \mathrm{d}  \textbf{P}_{E_2}}     \bigg] \bigg\}         \bigg]       +    \frac{1}{6}       \bigg[ 1-  \mathrm{exp} \big[       - 3^{-1} \lambda            \big]           \bigg] \\ \times   \mathscr{E} \bigg\{      \mathrm{log} \bigg[   \frac{ \mathrm{d}  \textbf{P}^{*}_{E_1}}{ \mathrm{d}  \textbf{P}_{E_1} +  \mathrm{d}  \textbf{P}_{E_2}}     \bigg] \tau  - \lambda  \underset{1 \leq i \leq \tau}{\sum} \big| \widetilde{\mathscr{U}_i} \big|        \bigg\}    ,
\end{align*} }

\noindent from the fact that,

{\small \begin{align*}
     \underset{0 \leq k \leq + \infty}{\sum}     \big[ \textit{Degree of each vertex over the triangular lattice} \big]^{-k}  \mathrm{exp} \bigg[       - k \bigg\{ \mathrm{log} \bigg[ \frac{\textbf{P}^*_{E_1}}{\textbf{P}_{E_1} + \textbf{P}_{E_2}}         \bigg] - \lambda    \bigg\}         \bigg]                  \\ \\ =     \underset{0 \leq k \leq + \infty}{\sum}    3^{-k}  \mathrm{exp} \bigg[       - k \bigg\{ \mathrm{log} \bigg[ \frac{\textbf{P}^*_{E_1}}{\textbf{P}_{E_1} + \textbf{P}_{E_2}}         \bigg] - \lambda    \bigg\}         \bigg]                   \\ \\ =   \big\{  1-  \mathrm{exp} \big[       - 3^{-1} \lambda            \big]           \big\}^{-1}                 ,
\end{align*}}

\noindent and also that,

{\small \begin{align*}
   \mathscr{E} \bigg\{ \mathrm{exp} \bigg[        \bigg\{ \mathrm{log} \bigg[   \frac{ \mathrm{d}  \textbf{P}^{*}_{E_2}}{ \mathrm{d}  \textbf{P}_{E_1} +  \mathrm{d}  \textbf{P}_{E_2}}     \bigg] \tau^{\prime}  - \lambda \bigg\}         \bigg]    \textbf{1}_{ \bigg\{ \widetilde{\mathscr{V}} : \big\{ \widetilde{\mathscr{V}}_{N^{\prime\prime}} = 0  \big\} , \big\{ \widetilde{\mathscr{V}}_0 = N^{\prime\prime} \big\} \bigg\}  }        \bigg|   \big\{   \textit{Steps of the triangular effective random walk}        \end{align*}

  \begin{align*}  \textit{between times $0$ and $N^{\prime\prime}-1$ along $E_2$}      \big\}    \bigg\}  \\ \\ =    \mathscr{E} \bigg\{ \mathrm{exp} \bigg[        \bigg\{ \mathrm{log} \bigg[   \frac{ \mathrm{d}  \textbf{P}^{*}_{E_2}}{ \mathrm{d}  \textbf{P}_{E_1} +  \mathrm{d}  \textbf{P}_{E_2}}     \bigg] \tau^{\prime}  - \lambda \bigg\}         \bigg]    \textbf{1}_{ \bigg\{ \widetilde{\mathscr{U}} : \big\{ \widetilde{\mathscr{V}}_{N^{\prime\prime}} = 0  \big\} , \big\{ \widetilde{\mathscr{V}}_0 = N^{\prime\prime}+1  \big\} \bigg\}  }           \bigg\}      ,
\end{align*} }

\noindent from the fact that,

{\small \begin{align*}
    \underset{0 \leq k \leq + \infty}{\sum}     \big[ \textit{Degree of each vertex over the triangular lattice} \big]^{-k}  \mathrm{exp} \bigg[       - k   \bigg\{   \mathrm{log} \bigg[ \frac{\textbf{P}^*_{E_2}}{\textbf{P}_{E_1} + \textbf{P}_{E_2}}         \bigg] - \lambda      \bigg\}   \bigg]              \\ \\ =     \underset{0 \leq k \leq + \infty}{\sum}    3^{-k}  \mathrm{exp} \bigg[       - k \bigg\{ \mathrm{log} \bigg[ \frac{\textbf{P}^*_{E_2}}{\textbf{P}_{E_1} + \textbf{P}_{E_2}}         \bigg] - \lambda    \bigg\}         \bigg]     \\ \\ =   \big\{ 1-  \mathrm{exp} \big[       - 3^{-1} \lambda            \big]           \big\}^{-1}                                                                     ,
\end{align*} }

\noindent implies that the Green's function along $E_2$ takes the form,

{\small \begin{align*}
  G_{E_2} \big( \lambda \big) =   \frac{1}{18}  \bigg[ 1 +  \mathrm{exp} \big[       - 3^{-1} \lambda            \big]           \bigg]   \mathrm{exp} \bigg[         \mathrm{log} \bigg\{    \frac{ \mathrm{d}  \textbf{P}^{*}_{E_2}}{ \mathrm{d}  \textbf{P}_{E_1} +  \mathrm{d}  \textbf{P}_{E_2}}     \bigg\}           \bigg]       +       \frac{1}{6}           \bigg[ 1-  \mathrm{exp} \big[       - 3^{-1} \lambda            \big]           \bigg]  \\ \times  \mathscr{E} \bigg\{      \mathrm{log} \bigg[   \frac{ \mathrm{d}  \textbf{P}^{*}_{E_2}}{ \mathrm{d}  \textbf{P}_{E_1} +  \mathrm{d}  \textbf{P}_{E_2}}     \bigg] \tau^{\prime}  - \lambda  \underset{1 \leq i \leq \tau}{\sum} \big| \widetilde{\mathscr{V}_i} \big|        \bigg\}    ,
\end{align*} }

\noindent As a result the first two terms appearing in the triangular Green's function are proportional to the superposition,

{\small \begin{align*}
 G_{E_1} \big( \lambda \big) \propto  \mathrm{exp} \bigg[      \mathrm{log} \bigg\{    \frac{ \mathrm{d}  \textbf{P}^{*}_{E_1}}{ \mathrm{d}  \textbf{P}_{E_1} +  \mathrm{d}  \textbf{P}_{E_2}}      \bigg\}         \bigg]   +  \mathscr{E} \bigg\{      \mathrm{log} \bigg[   \frac{ \mathrm{d}  \textbf{P}^{*}_{E_1}}{ \mathrm{d}  \textbf{P}_{E_1} +  \mathrm{d}  \textbf{P}_{E_2}}     \bigg] \tau  - \lambda  \underset{1 \leq i \leq \tau}{\sum} \big| \widetilde{\mathscr{U}_i} \big|        \bigg\}   ,
\end{align*} }

\noindent while the remaining two terms appearing in the triangular Green's function are proportional to the superposition,

{\small \begin{align*}
 G_{E_2} \big( \lambda \big) \propto    \mathrm{exp} \bigg[     \mathrm{log} \bigg\{    \frac{ \mathrm{d}  \textbf{P}^{*}_{E_2}}{ \mathrm{d}  \textbf{P}_{E_1} +  \mathrm{d}  \textbf{P}_{E_2}}      \bigg\}         \bigg]   +  \mathscr{E} \bigg\{      \mathrm{log} \bigg[   \frac{ \mathrm{d}  \textbf{P}^{*}_{E_2}}{ \mathrm{d}  \textbf{P}_{E_1} +  \mathrm{d}  \textbf{P}_{E_2}}     \bigg] \tau^{\prime}  - \lambda  \underset{1 \leq i \leq \tau^{\prime}}{\sum} \big| \widetilde{\mathscr{V}_i} \big|        \bigg\}  ,
\end{align*} }

  \noindent With the above representation of the triangular Green's function, to study the behavior of the sequence,
  
  \begin{align*}
  \mathscr{E} \bigg\{ \begin{bmatrix}
      M^{\lambda}_n      \\  M^{\lambda}_{n^{\prime}}
   \end{bmatrix}_{n,n^{\prime} \geq 0 }  \bigg\}  , 
  \end{align*}

  \noindent observe, under taking the limit of infinitely many steps,

 {\small \begin{align*}
   \underset{n , n^{\prime} \longrightarrow + \infty}{\mathrm{lim}}  \mathscr{E} \bigg\{ \begin{bmatrix}
      M^{\lambda}_n      \\  M^{\lambda}_{n^{\prime}}
   \end{bmatrix} \bigg\}  \overset{(\textit{linearity of } \mathscr{E} [ \cdot ])}{=}  \begin{bmatrix}
    \underset{n \longrightarrow + \infty}{\mathrm{lim}}  \mathscr{E}  \big\{   M^{\lambda}_n  \big\}      \\   \underset{n^{\prime} \longrightarrow + \infty}{\mathrm{lim}}  \mathscr{E} \big\{  M^{\lambda}_{n^{\prime}} \big\}  
   \end{bmatrix} =   \begin{bmatrix}
    \underset{n \longrightarrow + \infty}{\mathrm{lim}} \big\{  \mathscr{E}  \big[  M^{\lambda}_n  \big]  \textbf{1}_{ \{ \tau \leq n \}  }  \}  \\   \underset{n^{\prime} \longrightarrow + \infty}{\mathrm{lim}}  \big\{ \mathscr{E} \big[ M^{\lambda}_{n^{\prime}} \big]  \textbf{1}_{ \{ \tau^{\prime } \leq n^{\prime} \} } \big\} 
   \end{bmatrix} \\ \\  =     \begin{bmatrix}
    \underset{n \longrightarrow + \infty}{\mathrm{lim}} \big[  1 -  \mathscr{E}  \big\{  \big[  M^{\lambda}_n  \big]  \textbf{1}_{ \{ \tau >  n  \} }  \big\}  \big]  \\   \underset{n^{\prime} \longrightarrow + \infty}{\mathrm{lim}} \big[ 1 -  \mathscr{E} \big\{  \big[ M^{\lambda}_{n^{\prime}} \big]  \textbf{1}_{ \{ \tau^{\prime } > n^{\prime} \} } \big\}  \big] 
   \end{bmatrix}   \end{align*}

   \begin{align*} =     \underset{n , n^{\prime} \longrightarrow + \infty}{\mathrm{lim}}       \bigg\{  \begin{bmatrix}
  1 \\ 1  \end{bmatrix}  - \begin{bmatrix}
    \mathscr{E} \big\{  \big[  M^{\lambda}_n  \big]  \textbf{1}_{ \{ \tau >  n \} }  \big\}   \\    \mathscr{E} \big\{  \big[ M^{\lambda}_{n^{\prime}} \big]  \textbf{1}_{\{ \tau^{\prime } > n^{\prime} \} } \big\}  
   \end{bmatrix} \bigg\}  \\ \\ =              \begin{bmatrix}
  1 \\ 1  \end{bmatrix}  - \underset{n , n^{\prime} \longrightarrow + \infty}{\mathrm{lim}}   \begin{bmatrix}
    \mathscr{E} \big\{  \big[  M^{\lambda}_n  \big]  \textbf{1}_{\{ \tau >  n \} }  \big\}   \\    \mathscr{E} \big\{  \big[ M^{\lambda}_{n^{\prime}} \big]  \textbf{1}_{\{ \tau^{\prime } > n^{\prime} \} } \big\}  
   \end{bmatrix}     .
  \end{align*} }

  \noindent The above computations sharpen the estimate obtained by Fatou's lemma,

\begin{align*}
  \underset{n , n^{\prime} \longrightarrow + \infty}{\mathrm{lim}}  \mathscr{E} \bigg\{ \begin{bmatrix}
      M^{\lambda}_n      \\  M^{\lambda}_{n^{\prime}}
   \end{bmatrix} \bigg\}  \leq \begin{bmatrix} 1 \\ 1 
   \end{bmatrix}, 
\end{align*}

  \noindent by taking the limit as $n,n^{\prime} \longrightarrow + \infty$ in each component independently.

  \bigskip

\noindent Hence, the desired representation takes the form from combining the above groups of two terms, from which we conclude the argument, as it has been demonstrated that the above limit, as $n \longrightarrow + \infty$, equals $1$ in each component. \boxed{}

\section{Obtaining the desired scaling limit, from a suitable normalization over $\textbf{T}$}

\subsection{Description of objectives}

\noindent We present an adaptation of the truncation procedure defined in the previous section. By making use of the slab measure defined over $\textbf{T}$, we normalize by a probability that paths of the effective triangular walker exceeds a time dependent threshold. By making use of previous remarks, and comparisons, between such thresholds that have been obtained for the kinetic uniform prudent walk over $\textbf{Z}^2$, we formalize the following computations:

\begin{itemize}
    \item[$\bullet$] \textit{(1). Manipulating slab measures that have previously been introduced for prudent walks over the triangular lattice}. Write,

    \begin{align*}
        \hat{\mathscr{L}_R} : =        \mathscr{E}  \bigg\{      \textbf{1}_{\{ R = R_1 \} }                \mathrm{exp}     \bigg[   \mathrm{d} \bigg[  \frac{\textbf{P}^{*}_{E_1}}{\textbf{P}_{E_1} + \textbf{P}_{E_2}} \bigg]     n         \textbf{1}_{\{    \widetilde{\mathscr{U}} \in \{ 0 , \cdots, R \}  : 0 < \widetilde{U}_n <  R ,  \underset{1 \leq i \leq n}{\sum} | \widetilde{\mathscr{U}}_i |         = t - n                                                  \}}                      \bigg]                               +   \textbf{1}_{\{ R = R_2 \} }     \\ \times     \mathrm{exp} \bigg[    \mathrm{d} \bigg[ \frac{\textbf{P}^{*}_{E_2}}{\textbf{P}_{E_1} + \textbf{P}_{E_2}}       \bigg]   n^{\prime}         \textbf{1}_{\{  \widetilde{\mathscr{V}} \in \{ 0 , \cdots, R \}  :   0 < \widetilde{U}_n <  R ,  \underset{1 \leq i \leq n}{\sum} | \widetilde{\mathscr{V}}_i  |                 =    t - n^{\prime}          \}}                      \bigg]         \bigg\}   ,
    \end{align*}

       \noindent corresponding to a slab measure for the normalization $\mathscr{L}$.

    \item[$\bullet$] \textit{(2). Determining the time dependent threshold from the scaling limit normalization $\mathscr{L}$}. Write,

         \begin{align*}
         \textbf{P}^{*}_{\mathscr{T},\alpha} : =   \textbf{P}^{*}  \bigg[ \mathscr{T} > \mathscr{L} - \underset{1 \leq i \leq \alpha}{\sum} \mathscr{T}_i  \bigg]   =   \textbf{P}^{*}  \big[ \mathscr{T} > \mathscr{L} - \big(  \mathscr{T}_1 + \cdots + \mathscr{T}_{\alpha}  \big)  \big]     ,
         \end{align*}   

    \noindent corresponding to the previously defined probability for the triangular effective prudent walk, for the following times, and threshold,

    \begin{align*}
   \mathscr{T}  : = \textit{total time of the excursion of the triangular effective random walk}    , \\ \\  \mathscr{T}_i : = \textit{time that the triangular effective random walk spends at time i of an ex-} \\ \textit{ cursion}    , \\ \\  \underset{1 \leq i \leq \alpha}{\sum} \mathscr{T}_i : = \textit{time that the triangular effective random walk spends at time i of an ex-}  \end{align*}

   \begin{align*} \textit{ cursion, summed over all i } , \\ \\ \alpha : = \textit{threshold for which the inequality }  \bigg\{ \mathscr{T} > \mathscr{L} - \underset{1 \leq i \leq \alpha}{\sum} \mathscr{T}_i  \bigg\} \textit{ is satisfied},        
    \end{align*}

    \noindent respectively. Explicitly, recall that the form of $\alpha$,

    \begin{align*}
    \alpha : = \mathrm{log} \big[ \mathscr{L}^2 \big]    ,
        \end{align*}

\noindent satisfies,

\begin{itemize}
   \item[$\bullet$] \textit{Comparing the magnitude of the scaling limit normalization,} $\mathscr{L}$, \textit{over the triangular lattice to the total time spent at} $\mathrm{log} \big[ \mathscr{L}^2 \big]$ \textit{many excursions}. Denote,

    \begin{align*}
     \mathscr{T} : = \big[  \textit{Time spent along the first degree of freedom}     ,  \textit{Time spent along the second degree of freedom}     \big] \\ = \big[ \mathcal{T} , \widetilde{\mathcal{T}} \big]  ,
    \end{align*}

   \noindent from which one has that,

    \begin{align*}
    \frac{\mathscr{L} - \underset{1 \leq i \leq \mathrm{log} [ \mathscr{L}^2 ]}{\sum} \mathscr{T}_i }{\big( \mathrm{log} \mathscr{L} \big)^4} \longrightarrow 0     ,
    \end{align*}

\noindent as $\mathscr{L} \longrightarrow + \infty$.

   \item[$\bullet$] \textit{Comparing the time spent at $\mathrm{log} \big[ \mathscr{L}^4 \big]$ many excursions of the triangular random walk}. Denote the indicator function,

   \begin{align*}
     \textbf{1}_{\big\{  \mathcal{T} + \widetilde{\mathcal{T}} \geq \mathrm{log} [ \mathscr{L}^4  ]     \big\}}      ,
   \end{align*}

   \noindent corresponding to the occurrence of the event $\big\{  \mathcal{T} + \widetilde{\mathcal{T}} \geq \mathrm{log} [ \mathscr{L}^4  ]   \big\}$, and,

    \begin{align*}
     \mathscr{L} : = \big[ \mathscr{L}_1 , \mathscr{L}_2 \big]   ,
    \end{align*}

   \noindent from which one has that,

{\tiny \begin{align*}
        \bigg\{    \frac{\mathscr{L}_1 - \underset{1 \leq i \leq \mathrm{log} [ \mathscr{L}^2 ]}{\sum} \mathcal{T}_i }{\big( \mathrm{log} \mathscr{L} \big)^2} \longrightarrow 0 ,     \frac{\mathscr{L}_2 - \underset{1 \leq i \leq \mathrm{log} [ \mathscr{L}^2 ]}{\sum} \widetilde{\mathcal{T}_i}  }{\big( \mathrm{log} \mathscr{L} \big)^2} \longrightarrow 0    \bigg\} \Longleftrightarrow     \bigg\{          \frac{\mathscr{L} - \underset{1 \leq i \leq \mathrm{log} [ \mathscr{L}^2 ]}{\sum} \mathscr{T}_i }{\big( \mathrm{log} \mathscr{L} \big)^4} \longrightarrow 0     \bigg\}           , 
\end{align*}}

\noindent as $\mathscr{L} \longrightarrow + \infty$.

\item[$\bullet$] \textit{Comparing the time spend at $\mathrm{log} \big[ \mathscr{L}^6 \big]$ many excursions of the effective triangular random walk}. Denote,

\begin{align*}
    \widetilde{\widetilde{\mathcal{T}}} > 0 , 
\end{align*}

\noindent corresponding to the time spent about the third degree of freedom of the effective triangular random walk, and the indicator function,

\begin{align*}
  \textbf{1}_{\big\{   \mathcal{T} + \widetilde{\mathcal{T}} +    \widetilde{\widetilde{\mathcal{T}}} \geq \mathrm{log} [ \mathscr{L}^6 ]    \big\}}  ,
\end{align*}

\noindent corresponding to the occurrence of the event $\big\{   \mathcal{T} + \widetilde{\mathcal{T}} +    \widetilde{\widetilde{\mathcal{T}}} \geq \mathrm{log} [ \mathscr{L}^6 ]    \big\}$, from which one has that,

    {\small \begin{align*}
    \frac{\mathscr{L}^3  - \underset{1 \leq i \leq \mathrm{log} [ \mathscr{L}^2 ]}{\sum} \big[  \mathcal{T}_i + \widetilde{\mathcal{T}_i} +    \widetilde{\widetilde{\mathcal{T}_i}} \big]  }{\big( \mathrm{log} \mathscr{L} \big)^8} \longrightarrow 0    .
    \end{align*}}

\end{itemize}

    \item[$\bullet$] \textit{(3). Normalizing the first above quantity, corresponding to the slab measure, with the second above quantity, corresponding to the probability that the normalization $\mathscr{L}$}. Write,

   \begin{align*}
        \frac{\hat{\mathscr{L}_R}}{\textbf{P}^{*}_{\mathscr{T},\alpha}}
        ,
        \end{align*}

    \noindent corresponding to the normalization of the slab measure, introduced in the first item, with the probability introduced in the second item.

\item[$\bullet$] \textit{(4). Determining the largest order of the normalized slab measure introduced in the third item above}. Does there exist an up to constant, $C$, where $C$ is strictly positive, that can sharpen the following estimate,

   \begin{align*}
        \frac{\hat{\mathscr{L}_R}}{\textbf{P}^{*}_{\mathscr{T},\alpha}}
           \lesssim t^2       ,
        \end{align*}

\noindent for,

   \begin{align*}
        \frac{\hat{\mathscr{L}_R}}{\textbf{P}^{*}_{\mathscr{T},\alpha}}
        ,
        \end{align*}

\noindent to the inequality, for some strictly positive $C$,

   \begin{align*}
        \frac{\hat{\mathscr{L}_R}}{\textbf{P}^{*}_{\mathscr{T},\alpha}}
        \leq C t^2    ?
        \end{align*}

\end{itemize}

\subsection{Statement of objectives}

\noindent We provide the following statements which will be used to determine the time threshold for $\mathscr{L}$. They include:

\begin{itemize}
    \item[$\bullet$] \textit{(1). The triangular effective random walk does crosses the maximum displaced distance about the first, second, and third, degrees of freedom over $\textbf{T}$}. One has that,

    \begin{align*}
     \frac{1}{\sqrt{3} }    \underset{\mathscr{L} \longrightarrow + \infty}{\mathrm{lim}}            \bigg[           \textbf{P}_{\mathscr{U},\textbf{T}} \big[   \exists i \in \big\{ \mathrm{log} \big[ \mathscr{L}^2 \big], \cdots ,  \gamma_{e_1} \big\}    : \textbf{1}^{\prime} = 1 \textit{ for every i}  \big]        + \textbf{P}_{\mathscr{U},\textbf{T}} \big[   \exists i \in \big\{ \mathrm{log} \big[ \mathscr{L}^2 \big] , \cdots ,  \gamma_{e_2} \big\}   \\  : \textbf{2}^{\prime} = 1   \textit{ for every i} \big]     + \textbf{P}_{\mathscr{U},\textbf{T}} \big[    \exists i \in \big\{ \mathrm{log} \big[ \mathscr{L}^2 \big] , \cdots ,  \gamma_{e_3} \big\}     : \textbf{3}^{\prime} = 1  \textit{ for every i} \big]  \bigg]  ,
    \end{align*}

    \noindent corresponding to the superposition,

\[ \left\{\!\begin{array}{ll@{}>{{}}l} 
 (\textit{Probability that the triangular effective random walk contains more than $\mathrm{log} \big[ \mathscr{L}^2 \big]$  many} \\ \textit{  excursions, conditionally upon the indicator $\textbf{1}^{\prime}$ equaling 1} ) \\  \\ + (\textit{Probability that the triangular effective random walk contains more than $\mathrm{log} \big[ \mathscr{L}^2 \big]$ many} \\ \textit{  excursions,  conditionally upon the indicator $\textbf{2}^{\prime}$ equaling 1} ) \\  \\ +  (\textit{Probability that the triangular effective random walk contains more than $\mathrm{log} \big[ \mathscr{L}^2 \big]$ many} \\ \textit{  excursions, conditionally upon the indicator $\textbf{3}^{\prime}$ equaling 1} ) , 
\end{array}\right. 
\]

    \noindent of probabilities normalized by $\sqrt{3}$, as $\mathscr{L} \longrightarrow + \infty$, for,

    \begin{align*}
        \gamma_{e_1} : = \textit{Path along the first degree of freedom of the triangular effective random walk} , \\ \\   \gamma_{e_2} : = \textit{Path along the second degree of freedom of the triangular effective random walk} , \\ \\   \gamma_{e_3} : = \textit{Path along the third degree of freedom of the triangular effective random walk} .
    \end{align*}

     \item[$\bullet$] \textit{(2). Upper bounding the Radon-Nikodym derivative terms appearing in the definition of the slab measure}. Given a slab measure,
     
\begin{align*}
    \hat{\mathscr{L}_R} : =        \mathscr{E}  \bigg\{      \textbf{1}_{\{ R = R_1 \} }                \mathrm{exp}     \bigg[   \mathrm{d} \bigg[  \frac{\textbf{P}^{*}_{E_1}}{\textbf{P}_{E_1} + \textbf{P}_{E_2}} \bigg]       n       \textbf{1}_{\{    \widetilde{\mathscr{U}} \in \{ 0 , \cdots, R \}  : 0 < \widetilde{U}_n <  R ,  \underset{1 \leq i \leq n}{\sum} | \widetilde{\mathscr{U}}_i |         = t - n                                                  \}}                      \bigg]                               +   \textbf{1}_{\{ R = R_2 \} }      \\ \times    \mathrm{exp} \bigg[    \mathrm{d} \bigg[ \frac{\textbf{P}^{*}_{E_2}}{\textbf{P}_{E_1} + \textbf{P}_{E_2}}       \bigg]   n^{\prime}         \textbf{1}_{\{  \widetilde{\mathscr{V}} \in \{ 0 , \cdots, R \}  :   0 < \widetilde{U}_n <  R ,  \underset{1 \leq i \leq n}{\sum} | \widetilde{\mathscr{V}}_i  |                 =    t - n^{\prime}          \}}                      \bigg]         \bigg\} ,
\end{align*}

     \noindent with the decomposition,

     \begin{align*}
        \hat{\mathscr{L}_R} = \underset{1 \leq j \leq R-1}{\sum} \hat{\mathscr{L}_j} , 
     \end{align*}

     \noindent dependent upon the index of summation $j$, one must also construct an upper bound of the form,

{\small \begin{align*}
    \hat{\mathscr{L}_R} <     \underset{1 \leq j \leq R-1}{\sum} \bigg[   \mathscr{E}  \bigg\{      \textbf{1}_{\{ j = R_1 \} }                \mathrm{exp}     \bigg[   \mathrm{d} \bigg[  \frac{\textbf{P}^{*}_{E_1}}{\textbf{P}_{E_1} + \textbf{P}_{E_2}} \bigg]      n       \textbf{1}_{\{    \widetilde{\mathscr{U}} \in \{ 0 , \cdots, j \}  : 0 < \widetilde{U}_n <  j ,  \underset{1 \leq i \leq n}{\sum} | \widetilde{\mathscr{U}}_i |         = t - n                                                  \}}                      \bigg]                               +   \textbf{1}_{\{ j = R_2 \} }     \\ \times     \mathrm{exp} \bigg[    \mathrm{d} \bigg[ \frac{\textbf{P}^{*}_{E_2}}{\textbf{P}_{E_1} + \textbf{P}_{E_2}}       \bigg]   n^{\prime}          \textbf{1}_{\{  \widetilde{\mathscr{V}} \in \{ 0 , \cdots, j \}  :   0 < \widetilde{U}_n <  j ,  \underset{1 \leq i \leq n}{\sum} | \widetilde{\mathscr{V}}_i  |                 =    t - n^{\prime}          \}}                      \bigg]         \bigg\}  \bigg]              \\ \\ \overset{(*)}{\lesssim}     \underset{1 \leq j \leq R-1}{\sum} \bigg[   \mathscr{E}  \bigg\{      \textbf{1}_{\{ j = R_1 \} }                \mathrm{exp}     \bigg[   \mathrm{d} \bigg[  \frac{\textbf{P}^{**}_{E_1}}{\textbf{P}_{E_1} + \textbf{P}_{E_2}} \bigg]        n     \textbf{1}_{\{    \widetilde{\mathscr{U}} \in \{ 0 , \cdots, j \}  : 0 < \widetilde{U}_n <  j ,  \underset{1 \leq i \leq n}{\sum} | \widetilde{\mathscr{U}}_i |         = t - n                                                  \}}                      \bigg]                               +   \textbf{1}_{\{ j = R_2 \} }     \\ \times     \mathrm{exp} \bigg[    \mathrm{d} \bigg[ \frac{\textbf{P}^{**}_{E_2}}{\textbf{P}_{E_1} + \textbf{P}_{E_2}}       \bigg]   n^{\prime}          \textbf{1}_{\{  \widetilde{\mathscr{V}} \in \{ 0 , \cdots, j \}  :   0 < \widetilde{U}_n <  j ,  \underset{1 \leq i \leq n}{\sum} | \widetilde{\mathscr{V}}_i  |                 =    t - n^{\prime}          \}}                      \bigg]         \bigg\}  \bigg] \end{align*}

    \begin{align*} <       \underset{1 \leq j \leq R-1}{\sum} \bigg[   \mathscr{E}  \bigg\{      \textbf{1}_{\{ j = R_1 \} }                \mathrm{exp}     \bigg[   \mathrm{d} \bigg[  \frac{\textbf{P}^{**}_{E_1}}{\textbf{P}_{E_1} + \textbf{P}_{E_2}} \bigg]  n \textbf{1}_{\{    \widetilde{\mathscr{U}} \in \{ 0 , \cdots, j \}  : 0 < \widetilde{U}_n <  j ,  \underset{1 \leq i \leq n}{\sum} | \widetilde{\mathscr{U}}_i |         = t - n                                                  \}}  \\  \times            \textbf{1}_{\{      \widetilde{\mathscr{U}} :   \widetilde{\mathscr{U}}   \textit{has a decomposition when the effective path reaches a displacement of } \frac{R_1}{2}          \textit{from the reflection } \mathcal{R}_1                   \}}                      \bigg]                         \\       +   \textbf{1}_{\{ j = R_2 \} }        \mathrm{exp} \bigg[    \mathrm{d} \bigg[ \frac{\textbf{P}^{**}_{E_2}}{\textbf{P}_{E_1} + \textbf{P}_{E_2}}       \bigg]      n^{\prime}          \textbf{1}_{\{  \widetilde{\mathscr{V}} \in \{ 0 , \cdots, j \}  :   0 < \widetilde{U}_n <  j ,  \underset{1 \leq i \leq n}{\sum} | \widetilde{\mathscr{V}}_i  |                 =    t - n^{\prime}          \}}       \\ \times    \textbf{1}_{\{   \widetilde{\mathscr{U}} :   \widetilde{\mathscr{U}}   \textit{has a decomposition when the effective path reaches a displacement of } \frac{R_2}{2}        \textit{from the reflection } \mathcal{R}_2          \}}                      \bigg]         \bigg\}  \bigg]                 ,
\end{align*} }

\noindent where,

{\small \begin{align*}
   \big\{ E_1 , \textbf{P}^{**}_{E_1} \big[ \cdot \big] , \textbf{P}^{*}_{E_1} \big[ \cdot \big] :  \textbf{P}^{**}_{E_1} > \textbf{P}^{*}_{E_1} \big\}  \Longleftrightarrow  \bigg\{ E_1, E_2 , \textbf{P}^{**}_{E_1} \big[ \cdot \big] , \textbf{P}^{*}_{E_1} \big[ \cdot \big] , \textbf{P}_{E_1} \big[ \cdot \big] \\ , \textbf{P}_{E_2} \big[ \cdot \big]   : \frac{\textbf{P}^{**}_{E_1}}{\textbf{P}_{E_1} + \textbf{P}_{E_2}}  > \frac{\textbf{P}^{*}_{E_1}}{\textbf{P}_{E_1} + \textbf{P}_{E_2}}  \bigg\} , \\ \\  \big\{ E_2 ,\textbf{P}^{**}_{E_2} \big[ \cdot \big] , \textbf{P}^{*}_{E_2} \big[ \cdot \big] : \textbf{P}^{**}_{E_2} > \textbf{P}^{*}_{E_2} \big\} \Longleftrightarrow  \bigg\{ E_1, E_2 , \textbf{P}^{**}_{E_2} \big[ \cdot \big] , \textbf{P}^{*}_{E_2} \big[ \cdot \big] , \textbf{P}_{E_1} \big[ \cdot \big]   \end{align*}

   \begin{align*} , \textbf{P}_{E_2} \big[ \cdot \big]   :\frac{\textbf{P}^{**}_{E_2}}{\textbf{P}_{E_1} + \textbf{P}_{E_2}}   > \frac{\textbf{P}^{*}_{E_2}}{\textbf{P}_{E_1} + \textbf{P}_{E_2}} \bigg\}  , \\ \\  \mathcal{R}_1 : =  \textit{Reflection of a triangular effective random path about the $R_1$ degree of } \\ \textit{freedom}  , \\ \\            \mathcal{R}_2 : =         \textit{Reflection of a triangular effective random path about the $R_2$ degree of} \\ \textit{ freedom}        . 
\end{align*} } 

     \noindent $\mathcal{R}_1$ and $\mathcal{R}_2$ introduced above, as reflections on triangular effective random paths, additionally satisfy:

        \begin{itemize}
            \item[$\bullet$] \textit{Orthogonality of reflections of triangular effective random paths}. Denote $U_1$ and $V_1$ as,

\begin{align*}
  U_1 \equiv \textit{An excursion of the triangular effective random walk along the first degree of freedom}  , \\ \\    V_1 \equiv   \textit{An excursion of the triangular effective random walk along the second degree of fre-} \\ \textit{edom}  . 
\end{align*}

            \noindent $\mathcal{R}_1$ and $\mathcal{R}_2$ satisfy,

            \begin{align*}
       \mathcal{R}_1 \big( U_1 \cup V_1 \big) + \mathcal{R}_2 \big( U_1 \cup V_1 \big) =  \mathcal{R}_1 \big( U_1 \big) + \mathcal{R}_2 \big( V_1 \big)         .    
           \end{align*}

            \item[$\bullet$] \textit{Reflections of empty paths which have not taken any steps}. $\mathcal{R}_1 \big( \emptyset \big) = \mathcal{R}_2 \big( \emptyset \big) = 0$.

            \item[$\bullet$] \textit{Additivity of the reflections $\mathcal{R}_1$ and $\mathcal{R}_2$}. Denote $U_2$ and $V_2$ as,

           {\small  \begin{align*}
              U_2 \equiv \textit{An excursion of the triangular effective random walk along the first degree of}  \\ \textit{ freedom which is disjoint from $U_1$}   , \\ \\  V_2 \equiv \textit{An excursion of the triangular effective random walk along the second degree }  \\ \textit{ of freedom which is disjoint from $V_1$}   . 
            \end{align*} } 

\noindent $\mathcal{R}_1$ and $\mathcal{R}_2$ satisfy,

{\small \begin{align*}
   \mathcal{R}_1 \big( U_1 \cup U_2 \big) + \mathcal{R}_2 \big( V_1 \cup V_2 \big) =  \mathcal{R}_1 \big( U_1 + U_2 \big) + \mathcal{R}_2 \big( V_1 + V_2 \big) =  \mathcal{R}_1 \big( U_1 \big)  + \mathcal{R}_1 \big(  U_2 \big)  + \mathcal{R}_2 \big( V_1 \big) \\  + \mathcal{R}_2 \big(  V_2 \big)  .
\end{align*}} 
            
        \end{itemize}

\noindent Besides $\mathcal{R}_1$ and $\mathcal{R}_2$, introduce the reflection,

\begin{align*}
  \mathcal{R}_3 : =  \textit{Reflection of a triangular effective random walk} ,
\end{align*}

\noindent satisfies:

\begin{itemize}
\item[$\bullet$] \textit{Time dependence}. In comparison to $\mathcal{R}_1$ and $\mathcal{R}_2$ introduced above, the reflection $\mathcal{R}_3$, as a function of time, exhibits the dependency,

\begin{align*}
  \mathcal{R}_3 \big( t \big)  : =  \textit{Reflection of a triangular effective random walk at time t} ,
\end{align*}

\noindent corresponding to performing a reflection of a triangular effective random path at some $t > 0$.

\item[$\bullet$] \textit{Reflections of empty paths which have not taken any steps}. $\mathcal{R}_3 \big( \emptyset \big) = 0$.

\item[$\bullet$] \textit{The action of the third reflection, $\mathcal{R}_3$, at some prescribed time $t$}. Fix some time $t$. To distinguish the action of the reflection $\mathcal{R}_3$ from the reflections $\mathcal{R}_1$ and $\mathcal{R}_2$ above, we denote,

\begin{align*}
    \mathcal{R}_3   \curvearrowright     \pi_{\textit{effective}}  : =
   \textit{the reflection $\mathcal{R}_3$ is applied to the triangular effective random path $\pi_{\textit{effective}}$}. 
    \end{align*}

\noindent Explicitly, the action of $\mathcal{R}_3$, for some,

\begin{align*}
 \pi_{\textit{Effective}} \in  \mathscr{E}\mathscr{W}_{\textbf{T}}  , 
\end{align*}

\noindent is given by,

{\small \begin{align*}
    \mathcal{R}_3   \curvearrowright     \pi_{\textit{effective}}  : = \left\{\!\begin{array}{ll@{}>{{}}l}          \textit{(1). Partition } \pi_{\textit{effective}} \text{ } \textit{before, and after the time $t$ at which the reflection} \\ \textit{is performed, } \\ \\ \textit{(2). Fix some $t \equiv T>0$. Apply $\mathcal{R}_3$ to the partition $\big( \pi_{\textit{effective}} \big)_{[ 0 , T ]}$ of $\pi_{\textit{effective}}$,} \\ \textit{that is, the effective triangular path which equals $\pi_{\textit{effective}}$ up to time $T$}. \\ \\ \textit{(3). Fix $T_{\textit{final}} > T$ corresponding to the time time step at which the triangular} \\ \textit{effective random walk takes a step. To obtain the effective triangular path after} \\ \textit{ applying $\mathcal{R}_3$, concatenate $\big[ \mathcal{R}_3 \curvearrowright     \big( \pi_{\textit{effective}}   \big)_{[0,T]}   \big]$, which outputs the path } \\ \textit{$\big[ \mathcal{R}_3 \curvearrowright     \big( \pi_{\textit{effective}}   \big)_{[0,T]}   \big] \cup \big( \pi_{\textit{effective}} \big)_{(T,T_{\textit{final}}]}$.}  \\ \\ \textit{(4). By varying the time at which $\mathcal{R}_3$ is applied, one can strictly upper} \\ \textit{  bound $ \mathscr{L}_R$}. 
    \end{array}\right. 
\end{align*}     }

\item[$\bullet$] \textit{Composability of $\mathcal{R}_3$}. Given the definition of the action of $\mathcal{R}_3$ in the previous item above,

\end{itemize}

      \item[$\bullet$] \textit{(3). Expressing the dependency of the up to constant bound provided in (*) upon $R$, $j$, and the total number of steps of a triangular effective path}. The up to constant required for sharpening the up to constants bound provided in $(*)$ above takes the form,

      \begin{align*}
            6  \big( R  - 1 \big)^2     t^2                              . 
      \end{align*}

      \item[$\bullet$] \textit{(4). Indicator functions from the reflections $\mathcal{R}_1, \mathcal{R}_2, \mathcal{R}_3$}. Introduce,

     {\small \begin{align*}
       \textbf{1}_{\mathcal{R}_1} \big[ \mathcal{E}^{\prime} \big]  : =  \textit{Indicator function, corresponding to the first reflection $\mathcal{R}_1$, for which $\mathcal{E}^{\prime}$ belongs to the sample} \\ \textit{ space of the triangular effective walk, $ \mathscr{E}\mathscr{W}_{\textbf{T}} $, after applying $\mathcal{R}_1$}           , \\ \\    \textbf{1}_{\mathcal{R}_2}  \big[ \mathcal{E}^{\prime} \big]  : =   \textit{Indicator function, corresponding to the first reflection $\mathcal{R}_2$, for which $\mathcal{E}^{\prime}$ belongs to the sample} \\ \textit{ space of the triangular effective walk, $ \mathscr{E}\mathscr{W}_{\textbf{T}} $, after applying $\mathcal{R}_2$}              , \\ \\      \textbf{1}_{\mathcal{R}_3} \big[ \mathcal{E}^{\prime} \big]   : =       \textit{Indicator function, corresponding to the first reflection $\mathcal{R}_3$, for which $\mathcal{E}^{\prime}$ belongs to the sample} \\ \textit{ space of the triangular effective walk, $ \mathscr{E}\mathscr{W}_{\textbf{T}} $, after applying $\mathcal{R}_3$}         . 
      \end{align*}  }

\item[$\bullet$] \textit{(5). Upper bounding the slab measure from indicator functions of $\mathcal{R}_1, \mathcal{R}_2$, and $\mathcal{R}_3$}. The desired upper bound takes the form,

\begin{align*}
    3 t^2      ,
\end{align*}

\noindent for some $t>0 $. 

\noindent The up to constants upper bound, 

 \begin{align*}
        \frac{\hat{\mathscr{L}_R}}{\textbf{P}^{*}_{\mathscr{T},\alpha}}
           \lesssim t^3       ,
        \end{align*}

\noindent can be sharpened to the inequality,

 \begin{align*}
        \frac{\hat{\mathscr{L}_R}}{\textbf{P}^{*}_{\mathscr{T},\alpha}}
         \leq C t^3       ,
        \end{align*}

\noindent with the following computations,

{\small \begin{align*}
 \mathscr{L}_R <           \mathscr{E} \bigg\{  \bigg\{      \textbf{1}_{\{ R = R_1 \} }                \mathrm{exp}     \bigg[   \mathrm{d} \bigg[  \frac{\textbf{P}^{*}_{E_1}}{\textbf{P}_{E_1} + \textbf{P}_{E_2}} \bigg]       n      \textbf{1}_{\{    \widetilde{\mathscr{U}} \in \{ 0 , \cdots, R \}  : \widetilde{U}_n = R ,  \underset{1 \leq i \leq n}{\sum} | \widetilde{\mathscr{U}}_i |         = t - n                                                  \}}                      \bigg]                               +   \textbf{1}_{\{ R = R_2 \} }       \mathrm{exp} \bigg[    \mathrm{d} \bigg[ \frac{\textbf{P}^{*}_{E_2}}{\textbf{P}_{E_1} + \textbf{P}_{E_2}}       \bigg]    n^{\prime}       \\ \times    \textbf{1}_{\{  \widetilde{\mathscr{V}} \in \{ 0 , \cdots, R \}  :   \widetilde{U}_n = R ,  \underset{1 \leq i \leq n}{\sum} | \widetilde{\mathscr{V}}_i  |                 =    t - n^{\prime}          \}}                      \bigg]         \bigg\} \textbf{1}_{\mathcal{R}} \big[ \textit{The triangular effective random path crosses a distance $R$} \big]  \bigg\}  \\ \\ <     \mathscr{E} \bigg\{  \bigg\{      \textbf{1}_{\{ R = R_1 \} }                \mathrm{exp}     \bigg[   \mathrm{d} \bigg[  \frac{\textbf{P}^{**}_{E_1}}{\textbf{P}_{E_1} + \textbf{P}_{E_2}} \bigg]          n   \textbf{1}_{\{    \widetilde{\mathscr{U}} \in \{ 0 , \cdots, R \}  : \widetilde{U}_n = R ,  \underset{1 \leq i \leq n}{\sum} | \widetilde{\mathscr{U}}_i |         = t - n                                                  \}}                      \bigg]                               +   \textbf{1}_{\{ R = R_2 \} }       \mathrm{exp} \bigg[    \mathrm{d} \bigg[ \frac{\textbf{P}^{**}_{E_2}}{\textbf{P}_{E_1} + \textbf{P}_{E_2}}       \bigg]    n^{\prime}       \\ \times    \textbf{1}_{\{  \widetilde{\mathscr{V}} \in \{ 0 , \cdots, R \}  :   \widetilde{U}_n = R ,  \underset{1 \leq i \leq n^{\prime}}{\sum} | \widetilde{\mathscr{V}}_i  |                 =    t - n          \}}                      \bigg]         \bigg\} \textbf{1}_{\mathcal{R}} \big[ \textit{The triangular effective random path crosses a distance $R$} \big] \bigg\} \\ \\   <             \mathscr{E} \bigg\{  \bigg\{      \textbf{1}_{\{ R = R_1 \} }                \mathrm{exp}  \bigg[      \bigg[   \mathrm{d}  \bigg[  \frac{\textbf{P}^{**}_{E_1}}{\textbf{P}_{E_1} + \textbf{P}_{E_2}} \bigg]        +    \mathscr{C}_{E_1}      \bigg]  n   \textbf{1}_{\{    \widetilde{\mathscr{U}} \in \{ 0 , \cdots, R \}  : \widetilde{U}_n = R ,  \underset{1 \leq i \leq n}{\sum} | \widetilde{\mathscr{U}}_i |         = t - n                                                  \}}                      \bigg]                               +   \textbf{1}_{\{ R = R_2 \} }   \end{align*}

 \begin{align*} \times      \mathrm{exp} \bigg[  \bigg[   \mathrm{d} \bigg[ \frac{\textbf{P}^{**}_{E_2}}{\textbf{P}_{E_1} + \textbf{P}_{E_2}}       \bigg]    +  \mathscr{C}_{E_2}    \bigg]  n^{\prime}   \textbf{1}_{\{  \widetilde{\mathscr{V}} \in \{ 0 , \cdots, R \}  :   \widetilde{U}_n = R ,  \underset{1 \leq i \leq n}{\sum} | \widetilde{\mathscr{V}}_i  |                 =    t - n^{\prime}          \}}                      \bigg]         \bigg\} \textbf{1}_{\mathcal{R}} \big[ \textit{The triangular } \\  \textit{ effective random path crosses a distance $R$} \big] \bigg\} \\ \\ \overset{\textit{linearity of $\mathscr{E} [ \cdot \big]$}}{= }      \mathscr{E} \bigg\{  \bigg\{      \textbf{1}_{\{ R = R_1 \} }                \mathrm{exp}     \bigg[   \bigg[   \mathrm{d} \bigg[ \frac{\textbf{P}^{**}_{E_1}}{\textbf{P}_{E_1} + \textbf{P}_{E_2}} \bigg]        +    \mathscr{C}_{E_1}      \bigg] n   \textbf{1}_{\{    \widetilde{\mathscr{U}} \in \{ 0 , \cdots, R \}  : \widetilde{U}_n = R ,  \underset{1 \leq i \leq n}{\sum} | \widetilde{\mathscr{U}}_i |         = t - n                                                  \}}                      \bigg] \textbf{1}_{\mathcal{R}_1} \big[ \textit{The triangular ef-} \\ \textit{fective random path crosses a distance $R$} \big]   \bigg\}                \\   +                  \mathscr{E} \bigg\{  \bigg\{       \textbf{1}_{\{ R = R_2 \} }  
        \mathrm{exp} \bigg[  \bigg[   \mathrm{d} \bigg[ \frac{\textbf{P}^{**}_{E_2}}{\textbf{P}_{E_1} + \textbf{P}_{E_2}}       \bigg]    +  \mathscr{C}_{E_2}    \bigg]   n^{\prime}  \textbf{1}_{\{  \widetilde{\mathscr{V}} \in \{ 0 , \cdots, R \}  :   \widetilde{U}_n = R ,  \underset{1 \leq i \leq n}{\sum} | \widetilde{\mathscr{V}}_i  |                 =    t - n^{\prime}          \}}                      \bigg]         \textbf{1}_{\mathcal{R}_2} \big[ \textit{The triangular ef-} \\  \textit{fective random path crosses a distance $R$} \big] \bigg\}               \end{align*}
        
        \begin{align*} <  \mathscr{E} \bigg\{  \bigg\{      \textbf{1}_{\{ R = R_1 \} }                \mathrm{exp}     \bigg[ \bigg\{   \bigg[   \mathrm{d} \bigg[ \frac{\textbf{P}^{**}_{E_1}}{\textbf{P}_{E_1} + \textbf{P}_{E_2}} \bigg]        +    \mathscr{C}_{E_1}      \bigg] n  + \bigg[   \mathrm{d} \bigg[ \frac{\textbf{P}^{**}_{E_1}}{\textbf{P}_{E_1} + \textbf{P}_{E_2}} \bigg]        +    \mathscr{C}_{E_1}      \bigg] \big[ n + 2 \big]    \bigg\} \end{align*}

 \begin{align*}  \times \textbf{1}_{\{    \widetilde{\mathscr{U}} \in \{ 0 , \cdots, R \}  : \widetilde{U}_n = R ,  \underset{1 \leq i \leq n}{\sum} | \widetilde{\mathscr{U}}_i |         = t - n                                                  \}}                      \bigg] \textbf{1}_{\mathcal{R}_1} \big[ \textit{The triangular ef-} \\ \textit{fective random path crosses a distance $R$} \big]   \bigg\}                \\    +                  \mathscr{E} \bigg\{  \bigg\{       \textbf{1}_{\{ R = R_2 \} }  
        \mathrm{exp} \bigg\{  \bigg[  \bigg[   \mathrm{d} \bigg[ \frac{\textbf{P}^{**}_{E_2}}{\textbf{P}_{E_1} + \textbf{P}_{E_2}}       \bigg]    +  \mathscr{C}_{E_2}    \bigg]   n^{\prime}  + \bigg[   \mathrm{d} \bigg[ \frac{\textbf{P}^{**}_{E_2}}{\textbf{P}_{E_1} + \textbf{P}_{E_2}}       \bigg]    +  \mathscr{C}_{E_2}    \bigg]  \big[  n^{\prime} + 2 \big]  \bigg\}  \\ \times  \textbf{1}_{\{  \widetilde{\mathscr{V}} \in \{ 0 , \cdots, R \}  :   \widetilde{U}_n = R ,  \underset{1 \leq i \leq n}{\sum} | \widetilde{\mathscr{V}}_i  |                 =    t - n^{\prime}          \}}                      \bigg]         \textbf{1}_{\mathcal{R}_2} \big[ \textit{The triangular ef-} \\ \textit{fective random path crosses a distance $R$} \big] \bigg\}  \end{align*}
        
        \begin{align*} \overset{\textit{linearity of $\mathscr{E} [ \cdot ]$}}{=}              \mathscr{E} \bigg\{  \bigg\{      \textbf{1}_{\{ R = R_1 \} }                \mathrm{exp}     \bigg[   \bigg[   \mathrm{d} \bigg[ \frac{\textbf{P}^{**}_{E_1}}{\textbf{P}_{E_1} + \textbf{P}_{E_2}} \bigg]        +    \mathscr{C}_{E_1}      \bigg] n   \textbf{1}_{\{    \widetilde{\mathscr{U}} \in \{ 0 , \cdots, R \}  : \widetilde{U}_n = R ,  \underset{1 \leq i \leq n}{\sum} | \widetilde{\mathscr{U}}_i |         = t - n                                                  \}}                      \bigg] \textbf{1}_{\mathcal{R}_1} \big[ \textit{The triangular ef-} \\ \textit{fective random path crosses a distance $R$} \big]   \bigg\}                \\ +        \mathscr{E} \bigg\{  \bigg\{      \textbf{1}_{\{ R = R_1 \} }                \mathrm{exp}     \bigg[   \bigg[   \mathrm{d} \bigg[ \frac{\textbf{P}^{**}_{E_1}}{\textbf{P}_{E_1} + \textbf{P}_{E_2}} \bigg]        +    \mathscr{C}_{E_1}      \bigg] \big[ n + 2 \big]   \textbf{1}_{\{    \widetilde{\mathscr{U}} \in \{ 0 , \cdots, R \}  : \widetilde{U}_n = R ,  \underset{1 \leq i \leq n}{\sum} | \widetilde{\mathscr{U}}_i |         = t - n                                                  \}}                      \bigg] \textbf{1}_{\mathcal{R}_1} \big[ \textit{The triangular ef-} \\ \textit{fective random path crosses a distance $R$} \big]   \bigg\}      \\  \\   +                  \mathscr{E} \bigg\{  \bigg\{       \textbf{1}_{\{ R = R_2 \} }  
        \mathrm{exp} \bigg[  \bigg[   \mathrm{d} \bigg[ \frac{\textbf{P}^{**}_{E_2}}{\textbf{P}_{E_1} + \textbf{P}_{E_2}}       \bigg]    +  \mathscr{C}_{E_2}    \bigg]   n^{\prime}  \textbf{1}_{\{  \widetilde{\mathscr{V}} \in \{ 0 , \cdots, R \}  :   \widetilde{U}_n = R ,  \underset{1 \leq i \leq n}{\sum} | \widetilde{\mathscr{V}}_i  |                 =    t - n^{\prime}          \}}                      \bigg]         \textbf{1}_{\mathcal{R}_2} \big[ \textit{The triangular ef-} \\ \textit{fective random path crosses a distance $R$} \big] \bigg\}    \end{align*}

        \begin{align*} +        \mathscr{E} \bigg\{  \bigg\{       \textbf{1}_{\{ R = R_2 \} }  
        \mathrm{exp} \bigg[  \bigg[   \mathrm{d} \bigg[ \frac{\textbf{P}^{**}_{E_2}}{\textbf{P}_{E_1} + \textbf{P}_{E_2}}       \bigg]    +  \mathscr{C}_{E_2}    \bigg] \big[   n^{\prime} + 2 \big]   \textbf{1}_{\{  \widetilde{\mathscr{V}} \in \{ 0 , \cdots, R \}  :   \widetilde{U}_n = R ,  \underset{1 \leq i \leq n}{\sum} | \widetilde{\mathscr{V}}_i  |                 =    t - n^{\prime}          \}}                      \bigg]         \textbf{1}_{\mathcal{R}_2} \big[ \textit{The triangular ef-} \\ \textit{fective random path crosses a distance $R$} \big] \bigg\}      , \\ \tag{**}
\end{align*} }

\noindent where, as shorthand above, denote,

\begin{align*}
 \textbf{1}_{\mathcal{R}} \big[ \cdot \big] : = \textit{Indicator function where any of the reflections $\mathcal{R} \in \big\{ \mathcal{R}_1, \mathcal{R}_2, \mathcal{R}_3 \big\}$ can be applied to the} \\ \textit{event $\cdot$}   . 
\end{align*}

\noindent Altogether, the above manipulations of the Radon-Nikodym derivative appearing in the truncated slab measure imply,

\begin{align*}
 \big\{   (**) <  3 t^2   \mathscr{L}_R     \big\}   \Longleftrightarrow \big\{  \mathscr{L}_R < 3 t^2   \mathscr{L}_R    \big\} .
\end{align*}

\item[$\bullet$] \textit{Normalizing the ordinary slab measure with the truncated slab measure}. One can strictly upper bound,

\begin{align*}
 \frac{\mathscr{L}_R}{\hat{\mathscr{L}_R}}               , 
\end{align*}

\noindent by making use of the observation in the previous item,

\begin{align*}
 \frac{\mathscr{L}_R}{\hat{\mathscr{L}_R}}        \lesssim      t^2 \textbf{1}_{\{ t \geq R \} } + \textbf{1}_{\{ t < R \} }                  , 
\end{align*}

\noindent for natural $t$, specifically through $\mathcal{R}_1, \mathcal{R}_2$ and $\mathcal{R}_3$. Hence, there exists a strictly positive constant, $C$, for which,

\begin{align*}
 \frac{\mathscr{L}_R}{\hat{\mathscr{L}_R}} <     C \bigg\{   t^2 \textbf{1}_{\{ t \geq R \} } + \textbf{1}_{\{ t< R \} }             \bigg\}     . 
\end{align*}

\end{itemize}

\subsection{Comparison with conditionally defined probabilities in the scaling limit previously obtained by Beffara, Friedli, and Velenik}

\noindent As described in \textit{1.1}, previous computations for the scaling limits of other random walks of interest have been obtained by Beffara, Friedli, Velenik {\color{blue}[4]}. While one adaptations of these approaches could be formed by defining the conditional probabilities, 

  {\small \[    \left\{\!\begin{array}{ll@{}>{{}}l} 
\textbf{P}_{\mathscr{U},\textbf{T}}  \big[               \mathscr{X}_k = m          \big|        \gamma_{[0 , \mathscr{T}_k ]}  \big] =     \textbf{P}_{\mathscr{U},\textbf{T}}   \big[           \eta_{\mathscr{X}_{[\mathscr{T}_k ]}}   = m  \big]             ,\\ \\  \textbf{P}_{\mathscr{U},\textbf{T}}   \big[              \mathscr{Y}_k = m    \big|    \gamma_{[0 , \mathscr{U}_k ]}          \big]   =     \textbf{P}_{\mathscr{U},\textbf{T}}   \big[     \eta_{\mathscr{Y}_{[\mathscr{U}_k ]}} = m      \big]              , \\ \\ \textbf{P}_{\mathscr{U},\textbf{T}}  \big[              \mathscr{Z}_k = m    \big|    \gamma_{[0 , \mathscr{V}_k ]}          \big]   =     \textbf{P}_{\mathscr{U},\textbf{T}}   \big[     \eta_{\mathscr{Z}_{[\mathscr{V}_k ]}} = m      \big]        ,     \end{array}\right. \]  }

\noindent $\forall k \geq 0 , m \geq 1$, where $\eta$ denotes a point process that is dependent upon either $\mathcal{H}_{\mathscr{T}_k}$, or upon $\mathcal{W}_{\mathscr{U}_k}$, each of which respectively correspond to the probabilities of,

{\small \begin{align*}
  \big\{ \mathscr{X}_k = m \big\}   , \\ \\ \big\{ \mathscr{Y}_k = m \big\}  , \\ \\ \big\{ \mathscr{Z}_k = m \big\}  , 
\end{align*} } 

\noindent occurring, one can manipulate the quantities,

 {\small \[    \left\{\!\begin{array}{ll@{}>{{}}l} 
     \frac{\mathscr{L}_R}{\hat{\mathscr{L}_R}}          \bigg|_{\textit{first degree of freedom of the triangular effective random walk}}  : = \textit{Restriction of }   \frac{\mathscr{L}_R}{\hat{\mathscr{L}_R}}  \textit{ to $e_1$} , \\  \\    \frac{\mathscr{L}_R}{\hat{\mathscr{L}_R}}          \bigg|_{\textit{second degree of freedom of the triangular effective random walk}}   : = \textit{Restriction of }   \frac{\mathscr{L}_R}{\hat{\mathscr{L}_R}}  \textit{ to $e_2$}   , \\  \\   \frac{\mathscr{L}_R}{\hat{\mathscr{L}_R}}          \bigg|_{\textit{third degree of freedom of the triangular effective random walk}}  : =  \textit{Restriction of }   \frac{\mathscr{L}_R}{\hat{\mathscr{L}_R}}  \textit{ to $e_3$}    ,     \end{array}\right. \] }

\noindent corresponding to the slab measure normalized with the truncated slab measure. About each degree of freedom of the effective triangular walk, one has that,

{\small \begin{align*}
         \textbf{P}_{\mathscr{U},\textbf{T}}   \bigg[    \text{selecting } \xi_l \text{ } \text{from the iid sequence }  \{ \xi_1 , \cdots , \xi_k \}  \\ \times  \textbf{1}_{ \{ \textit{iid sequence is supported along the first, or second, degrees of freedom $R_1$ or $R_2$}      \}}              \bigg]    = \frac{1}{3} \bigg\{   \frac{1}{2} \bigg\}^{|l|}   . 
\end{align*} }

\noindent Denoting the above probability with $p \big( \beta \big)$, conditionally upon the event that the prudent walker achieves some height $\alpha_i$ conditionally upon the history of the path at time $\mathscr{T}_{\alpha_i}$ occurs. Moreover, the product of probabilities admits the identity,

{\small \begin{align*}
\bigg\{   \underset{1 \leq i \leq m-1}{\prod}  \textbf{P}_{\mathscr{U},\textbf{T}}  \big[               \big\{ X_k = \alpha_i            \big\}  \textbf{1}_{ \{ \textit{iid sequence is supported along the first, or second, degrees of freedom $R_1$ or $R_2$}      \}}              \big]  \bigg\}   \\ \times \textbf{P}_{\mathscr{U},\textbf{T}}  \big[  \big\{ X_k = \beta \big\}   \textbf{1}_{ \{ \textit{iid sequence is supported along the first, or second, degrees of freedom $R_1$ or $R_2$}      \}}      \big] \\ \\ = \bigg\{   \underset{1 \leq i \leq m-1}{\prod}   p \big( \alpha_i \big)         \bigg\}   p \big( \beta \big) , 
 \end{align*} }

\noindent for a suitable transformation $p$ of each $\alpha_i$ and $\beta$. 

\bigskip

\noindent Over $\mathscr{T}_k$ or $\mathscr{U}_k$, the fact that the value of $\mathscr{X}_k$ and $\mathscr{Y}_k$ at the k $th$ step are only dependent upon the vertical, and horizontal, segments of the prudent walker, the first probability conditioned on $\gamma_{[0,\mathscr{T}_k]}$ can be expressed as,

\begin{align*}
\textbf{P}_{\mathscr{U},\textbf{T}}  \big[               
\mathscr{X}_k = m          \big|        \gamma_{[0 , \mathscr{T}_k ]}  \big] =       \underset{\{ \alpha_1 , \cdots , \alpha_n , \beta\}}{\sum} \bigg\{   \underset{1 \leq i \leq m-1}{\prod}  \textbf{P}_{\mathscr{U},\textbf{T}}  \big[  \mathscr{X}_k = \alpha_i \big]  \textbf{P}_{\mathscr{U},\textbf{T}}   \big[ \mathscr{X}_k = \beta \big] \bigg\}   , 
 \end{align*}

 \noindent for times $\big\{  \gamma_{[0,\mathscr{T}_{\alpha_i} ]}\big\}_{i \in \{i_1 , \cdots , i_n , i_{n+1} \}}$. The main task will revolve around evaluating the probability, for each $\alpha_i$, which is related to the probability,

 {\small \begin{align*}
\underset{\textit{increments of the walk along $R_2$}}{\sum } \bigg\{  \underset{\textit{increments of the walk along $R_1$}}{\sum }  \textbf{P}_{\mathscr{U}, \textbf{T}} \big[  \gamma_{R_1}   : \# \big\{ \textit{edges} \in \gamma_{R_1}         \big\} \leq      \mathrm{log} \big[ \mathscr{L}^2 \big]                \big] \\ \times \textbf{P}_{\mathscr{U}, \textbf{T}} \big[              \gamma_{R_2}   : \# \big\{ \textit{edges} \in \gamma_{R_2}         \big\} \leq    \mathrm{log} \big[ \mathscr{L}^2 \big]          \big]  \bigg\}  \\ \\  \lesssim    \underset{R_1, R_2}{\sum} \mathrm{exp} \big[ 2 \mathscr{L}^2 \big] \\ \\ \lesssim    \underset{R_1, R_2}{\sum} \mathrm{exp} \big[  \mathscr{L}^4 \big]   .
 \end{align*} }

\noindent To capture behaviors of the triangular random walk, we argue that the following holds:

\bigskip

\noindent \textbf{Lemma} \textit{2} (\textit{the slab measure, normalized by the truncated slab measure, satisfies the inequality provided at the end of 3.2}). There exists $C>0$ for which,

\begin{align*}
 \frac{\mathscr{L}_R}{\hat{\mathscr{L}_R}} <     C \bigg\{  t^2 \textbf{1}_{\{ t \geq R \} } + \textbf{1}_{\{ t< R \} }             \bigg\}     . 
\end{align*}

\bigskip

\noindent \textit{Proof of Lemma 2}. Apply the steps provided in \textit{3.2}; namely, to argue that a suitable bound for,

\noindent can be used to argue that,

 \begin{align*}
     \frac{1}{\sqrt{3} }               \bigg[           \textbf{P}_{\mathscr{U},\textbf{T}} \big[   \exists i \in \big\{ \mathrm{log} \big[ \mathscr{L}^2 \big], \cdots ,  \gamma_{e_1} \big\}    : \textbf{1}^{\prime} = 1 \textit{ for every i}  \big]        + \textbf{P}_{\mathscr{U},\textbf{T}} \big[   \exists i \in \big\{ \mathrm{log} \big[ \mathscr{L}^2 \big] , \cdots ,  \gamma_{e_2} \big\}   \\  : \textbf{2}^{\prime} = 1   \textit{ for every i} \big]     + \textbf{P}_{\mathscr{U},\textbf{T}} \big[    \exists i \in \big\{ \mathrm{log} \big[ \mathscr{L}^2 \big] , \cdots ,  \gamma_{e_3} \big\}     : \textbf{3}^{\prime} = 1  \textit{ for every i} \big]  \bigg]  , \end{align*}

     \begin{align*} \Updownarrow \end{align*}

      \[    \left\{\!\begin{array}{ll@{}>{{}}l}  \frac{1}{\sqrt{3} }                         \textbf{P}_{\mathscr{U},\textbf{T}} \big[   \exists i \in \big\{ \mathrm{log} \big[ \mathscr{L}^2 \big], \cdots ,  \gamma_{e_1} \big\}    : \textbf{1}^{\prime} = 1 \textit{ for every i}  \big] \longrightarrow 0 \textit{ as } \mathscr{L} \longrightarrow + \infty  ,   \\ \\ \frac{1}{\sqrt{3} }    \textbf{P}_{\mathscr{U},\textbf{T}} \big[   \exists i \in \big\{ \mathrm{log} \big[ \mathscr{L}^2 \big] , \cdots ,  \gamma_{e_2} \big\}    : \textbf{2}^{\prime} = 1   \textit{ for every i} \big] \longrightarrow 0 \textit{ as } \mathscr{L} \longrightarrow + \infty ,   \\ \\  \frac{1}{\sqrt{3} }    \textbf{P}_{\mathscr{U},\textbf{T}} \big[    \exists i \in \big\{ \mathrm{log} \big[ \mathscr{L}^2 \big] , \cdots ,  \gamma_{e_3} \big\}     : \textbf{3}^{\prime} = 1  \textit{ for every i} \big]   \longrightarrow 0 \textit{ as } \mathscr{L} \longrightarrow + \infty ,     \end{array}\right. \]

\noindent as provided in the statement of \textbf{Lemma} \textit{4} below, it suffices to demonstrate that an up to constants estimate of the form,

\begin{align*}
 \frac{\mathscr{L}_R}{\hat{\mathscr{L}_R}}        \lesssim      t^2 \textbf{1}_{\{ t \geq R \} } + \textbf{1}_{\{ t < R \} }                  , 
\end{align*}

\noindent holds. Directly applying the computations provided in \textit{(5)} of \textit{3.2}, the previous section, with the three reflection transformations $\mathcal{R}_1, \mathcal{R}_2$, and $\mathcal{R}_3$ introduced in \textit{(2)} implies the desired estimate, from which we conclude the argument. \boxed{}

\bigskip

\noindent Next, we state how the estimate for normalizing the slab measure, over a distance $R$, with the truncated slab measure implies the constant for which the estimate,

 \begin{align*}
        \frac{\hat{\mathscr{L}_R}}{\textbf{P}^{*}_{\mathscr{T},\alpha}}
           \lesssim t^3       ,
        \end{align*}

\noindent holds.

\bigskip

\noindent \textbf{Lemma} \textit{3} (\textit{under the uniform probability measure over $\textbf{T}$, the last excursion of the effective walk typically has length as most $\mathrm{log} \big[ \mathscr{L}^2 \big]$, up to a constant}). Take $\mathscr{C}_2 > \mathscr{C}_1 > \mathscr{C}$ and $\mathscr{L}$ to be sufficiently large. One has that,

{\small \begin{align*}
\textbf{P}_{\mathscr{U}, \textbf{T}} \bigg[              \gamma_{R_1} ,   \gamma_{R_2}   : \# \big\{ \textit{edges} \in \gamma_{R_1}         \big\} \leq    \mathrm{log} \big[ \mathscr{L}^2 \big]  ,  \# \big\{ \textit{edges} \in \gamma_{R_2}         \big\} \leq    \mathrm{log} \big[ \mathscr{L}^2 \big] \bigg] \end{align*}}

\noindent which is proportional to,

{\small
\begin{align*}
\frac{1}{\sqrt{3} }          \bigg\{                \textbf{P}_{\mathscr{U},\textbf{T}} \big[   \exists i \in  \big\{ \mathrm{log} \big[ \mathscr{L}^2 \big], \cdots ,  \gamma_{e_1} \big\}    : \textbf{1}^{\prime} = 1 \textit{ for every i}  \big]     +    \textbf{P}_{\mathscr{U},\textbf{T}} \big[   \exists i \in \big\{ \mathrm{log} \big[ \mathscr{L}^2 \big] , \cdots ,  \gamma_{e_2} \big\}    : \textbf{2}^{\prime} = 1   \textit{ for every i} \big]   \\  +     \textbf{P}_{\mathscr{U},\textbf{T}} \big[    \exists i \in \big\{ \mathrm{log} \big[ \mathscr{L}^2 \big] , \cdots ,  \gamma_{e_3} \big\}     : \textbf{3}^{\prime} = 1  \textit{ for every i} \big]  \bigg\}         , 
\end{align*} }

\noindent is stochastically dominated by,

{\small \begin{align*}
      \underset{\epsilon = 0, \textit{or } \epsilon =1}{\sum}  \bigg\{  \underset{\tau \leq n }{\sum} \bigg\{ \underset{\tau^{\prime} \leq n^{\prime}}{\sum}   \underset{i \neq j \in \textbf{N}}{\mathrm{sup}} \bigg[     \mathscr{L}_{R_{i-1}} \big[ n , n^{\prime} ,  \epsilon  , R_1 , R_2 \big]        ,    \hat{\mathscr{L}_{R_{i-1}}} \big[ n , n^{\prime} , \epsilon  , R_1 , R_2 \big]   ,     \mathscr{L}_{R_{j-1}} \big[ n , n^{\prime} ,  \epsilon  , R_1 , R_2 \big]      \\    ,      \hat{\mathscr{L}_{R_{j-1}}} \big[ n , n^{\prime} , \epsilon  , R_1 , R_2 \big]           \bigg]        \bigg\}     \bigg\}              , 
\end{align*}}

\noindent and hence,

{\small \begin{align*}
   \textbf{P}_{\mathscr{U}, \textbf{T}} \bigg[              \gamma_{R_1} ,   \gamma_{R_2}   : \# \big\{ \textit{edges} \in \gamma_{R_1}         \big\} \leq    \mathrm{log} \big[ \mathscr{L}^2 \big]  ,  \# \big\{ \textit{edges} \in \gamma_{R_2}         \big\} \leq    \mathrm{log} \big[ \mathscr{L}^2 \big] \bigg]  \\ \\ <    \underset{\epsilon = 0, \textit{or } \epsilon =1}{\sum}  \bigg\{  \underset{\tau \leq n }{\sum} \bigg\{ \underset{\tau^{\prime} \leq n^{\prime}}{\sum}   \underset{i \neq j \in \textbf{N}}{\mathrm{sup}} \bigg[     \mathscr{L}_{R_{i-1}} \big[ n , n^{\prime} ,  \epsilon  , R_1 , R_2 \big]        ,    \hat{\mathscr{L}_{R_{i-1}}} \big[ n , n^{\prime} , \epsilon \\  , R_1 , R_2 \big]   ,     \mathscr{L}_{R_{j-1}} \big[ n , n^{\prime} ,  \epsilon  , R_1 , R_2 \big]          ,      \hat{\mathscr{L}_{R_{j-1}}} \big[ n , n^{\prime} , \epsilon  , R_1 , R_2 \big]           \bigg]        \bigg\}     \bigg\}    \\ \\   <           \mathscr{C}_1 \mathscr{L}^4  \mathrm{exp} \bigg\{            -            \mathscr{C}_2 \mathrm{log} \big[ \mathscr{L}^2   \big]          \bigg\}             \longrightarrow 0           .
\end{align*} }

\noindent

\noindent \textit{Proof of Lemma 3}. We adapt the argument of provided in \textbf{Lemma} \textit{5.3} of {\color{blue}[4]}. Hence it suffices to argue,

\begin{align*}
 \underset{\mathscr{L} \longrightarrow + \infty}{\mathrm{lim}} \bigg\{ \textbf{P}_{\mathscr{U}, \textbf{T}} \big[              \gamma_{R_1}   : \# \big\{ \textit{edges} \in \gamma_{R_1}         \big\} \leq    \mathrm{log} \big[ \mathscr{L}^2 \big]   +  \textbf{P}_{\mathscr{U}, \textbf{T}} \big[              \gamma_{R_2}   : \# \big\{ \textit{edges} \in \gamma_{R_2}         \big\} \leq    \mathrm{log} \big[ \mathscr{L}^2 \big] \big] \bigg\}  \\ \longrightarrow 0 \end{align*} \begin{align*}  \Updownarrow \end{align*}

 \begin{align*} \frac{1}{\sqrt{3} }  \underset{\mathscr{L} \longrightarrow + \infty}{\mathrm{lim}}           \bigg\{                \textbf{P}_{\mathscr{U},\textbf{T}} \big[   \exists i \in  \big\{ \mathrm{log} \big[ \mathscr{L}^2 \big], \cdots ,  \gamma_{e_1} \big\}    : \textbf{1}^{\prime} = 1 \textit{ for every i}  \big]     +    \textbf{P}_{\mathscr{U},\textbf{T}} \big[   \exists i \in \big\{ \mathrm{log} \big[ \mathscr{L}^2 \big] , \\  \cdots ,  \gamma_{e_2} \big\}    : \textbf{2}^{\prime} = 1   \textit{ for every i} \big]     +     \textbf{P}_{\mathscr{U},\textbf{T}} \big[    \exists i \in \big\{ \mathrm{log} \big[ \mathscr{L}^2 \big] , \cdots ,  \gamma_{e_3} \big\}   \\   : \textbf{3}^{\prime} = 1  \textit{ for every i} \big]  \bigg\}   \longrightarrow 0       . 
\end{align*}

\noindent That is, by direct computation one obtains the desired asymptotic behavior from the observation that,

\begin{align*}
  \mathrm{exp} \big[   -   \mathrm{log} \big[ \mathscr{L}^2 \big]     \big]  ,
\end{align*}

\noindent holds, up to the constant, from the rearrangements,

{\tiny \begin{align*}
        \mathscr{E} \bigg\{     \underset{1 \leq i \leq \gamma_{\mathscr{L}}}{\prod}                   \frac{\mathscr{L}_{R_{i-1}}}{\hat{\mathscr{L}_{R_{i-1}}}}          \textbf{1}_{\big\{ \mathscr{L } - \underset{1 \leq  k \leq \gamma_{\mathscr{L}}}{\sum} \mathscr{T}_k \gtrsim \mathrm{log} [ \mathscr{L}^2 ] \big\}}            \bigg\}   \leq   \mathscr{E} \bigg\{ 1 +      \underset{1 \leq i \leq \gamma_{\mathscr{L}}}{\prod}                   \frac{\mathscr{L}_{R_{i-1}}}{\hat{\mathscr{L}_{R_{i-1}}}}          \textbf{1}_{\big\{ \mathscr{L } - \underset{1 \leq  k \leq \gamma_{\mathscr{L}}}{\sum} \mathscr{T}_k \gtrsim \mathrm{log} [ \mathscr{L}^2 ] \big\}}                     \bigg\}   \\ \\ \leq             \mathscr{E} \bigg\{  \bigg[  1  +   \bigg[    \underset{1 \leq j \leq \mathscr{L}_2}{\underset{1 \leq i \leq \mathscr{L}_1}{\prod}}         \widetilde{\mathscr{T}_i }      \mathscr{T}_j     \bigg]       \textbf{1}_{\big\{               i :     \widetilde{\mathscr{T}_i }         \geq \frac{1-i}{2}         \big\}} \textbf{1}_{\big\{          j : {\mathscr{T}_j}         \geq \frac{1-j}{2}         \big\}}      \bigg]                     \textbf{1}_{\big\{ \mathscr{L } - \underset{1 \leq  k \leq \gamma_{\mathscr{L}}}{\sum} \mathscr{T}_k \gtrsim \mathrm{log} [ \mathscr{L}^2 ] \big\}}                  \bigg\}      \\  \\ \leq \mathscr{E} \bigg\{    \bigg[   1  +   \bigg[    \underset{1 \leq j \leq \mathscr{L}_2}{\underset{1 \leq i \leq \mathscr{L}_1}{\prod}}         \widetilde{\mathscr{T}_i }      \mathscr{T}_j            \bigg]       \textbf{1}_{\big\{               i :     \mathrm{sup} \{ \widetilde{\mathscr{T}_i } \}      \gtrsim  \mathrm{log} [ \mathscr{L}^2 ]              \big\}}      \textbf{1}_{\big\{               j :     \mathrm{sup} \{ {\mathscr{T}_j } \}      \gtrsim  \mathrm{log} [ \mathscr{L}^2 ]              \big\}}                                                                    \bigg]   \\ \times              \textbf{1}_{\big\{ \mathscr{L } - \underset{1 \leq  k \leq \gamma_{\mathscr{L}}}{\sum} \mathscr{T}_k \gtrsim \mathrm{log} [ \mathscr{L}^2 ] \big\}}     \bigg\}  \\ \\ \leq    \mathscr{E} \bigg\{         \bigg[   1  +   \bigg[    \underset{\forall i , c_i > 0 }{\underset{1 \leq i \leq \mathscr{L}_1, 1 \leq j \leq \mathscr{L}_2 }{\prod}}      c_i    \widetilde{\mathscr{T}_i }        \mathscr{T}_j   \bigg]       \textbf{1}_{\big\{               i :     \mathrm{sup} \{ \widetilde{\mathscr{T}_i } \}      \gtrsim  \mathrm{log} [ \mathscr{L}^2 ]              \big\}}      \textbf{1}_{\big\{               j :     \mathrm{sup} \{ {\mathscr{T}_j } \}      \gtrsim  \mathrm{log} [ \mathscr{L}^2 ]              \big\}}                                                                    \bigg]           \end{align*}

        \begin{align*}   \times            \textbf{1}_{\big\{ \mathscr{L } - \underset{1 \leq  k \leq \gamma_{\mathscr{L}}}{\sum} \mathscr{T}_k \gtrsim \mathrm{log} [ \mathscr{L}^2 ] \big\}}                                                        \bigg\} \\ \\ \leq     \mathscr{L}^2   \mathscr{E} \bigg\{         \bigg[   1  +   \bigg[    \underset{\forall i , c_i > 0 }{\underset{1 \leq i \leq \mathscr{L}_1 , 1 \leq j \leq \mathscr{L}_2 }{\prod}}      c_i    \widetilde{\mathscr{T}_i }        \mathscr{T}_j   \bigg]       \textbf{1}_{\big\{               i :     \mathrm{sup} \{ \widetilde{\mathscr{T}_i } \}      \gtrsim  \mathrm{log} [ \mathscr{L}^2 ]              \big\}}      \textbf{1}_{\big\{               j :     \mathrm{sup} \{ {\mathscr{T}_j } \}      \gtrsim  \mathrm{log} [ \mathscr{L}^2 ]              \big\}}                                                                    \bigg]         \\ \times            \textbf{1}_{\big\{ \mathscr{L } - \underset{1 \leq  k \leq \gamma_{\mathscr{L}}}{\sum} \mathscr{T}_k \gtrsim \mathrm{log} [ \mathscr{L}^2 ] \big\}}                                                        \bigg\} \\ \\ \leq            \mathscr{L}^2   \mathscr{E} \bigg\{            1  +   \bigg[    \underset{\forall i , c_i > 0 }{\underset{1 \leq i \leq \mathscr{L}_1, 1 \leq j \leq \mathscr{L}_2 }{\prod}}      c_i    \widetilde{\mathscr{T}_i } \mathscr{T}_j          \bigg]       \textbf{1}_{\big\{               i :     \mathrm{sup} \{ \widetilde{\mathscr{T}_i } \}      \gtrsim  \mathrm{log} [ \mathscr{L}^2 ]              \big\}}      \textbf{1}_{\big\{               j :     \mathrm{sup} \{ {\mathscr{T}_j } \}      \gtrsim  \mathrm{log} [ \mathscr{L}^2 ]              \big\}}                                                                  \bigg\}  \\ \times   \underset{\forall i , c_i > 0 }{\underset{1 \leq i \leq \mathscr{L}}{\prod}}   \mathscr{E} \bigg\{         1  +   \bigg[    \underset{\forall i , c_i > 0 }{\underset{1 \leq i \leq \mathscr{L}, 1 \leq j \leq \mathscr{L}_2 }{\prod}}      c_i    \widetilde{\mathscr{T}_i }   \mathscr{T}_j        \bigg]       \textbf{1}_{\big\{               i :     \mathrm{sup} \{ \widetilde{\mathscr{T}_i } \}      \gtrsim  \mathrm{log} [ \mathscr{L}^2 ]              \big\}}      \textbf{1}_{\big\{               j :     \mathrm{sup} \{ {\mathscr{T}_j } \}      \gtrsim  \mathrm{log} [ \mathscr{L}^2 ]              \big\}}                                                                                                                      \bigg\}   \\ \\ \leq              \underset{\forall i, c_i > 0 }{\underset{1 \leq i \leq \mathscr{L}}{\mathrm{sup}}}  c_i              \mathscr{L}^2      \bigg\{                       \mathscr{E} \bigg\{            1  +   \bigg[    \underset{\forall i , c_i > 0 }{\underset{1 \leq i \leq \mathscr{L}_1, 1 \leq j \leq \mathscr{L}_2 }{\prod}}       \widetilde{\mathscr{T}_i }     \mathscr{T}_j      \bigg]       \textbf{1}_{\big\{               i :     \mathrm{sup} \{ \widetilde{\mathscr{T}_i } \}      \gtrsim  \mathrm{log} [ \mathscr{L}^2 ]              \big\}}      \textbf{1}_{\big\{               j :     \mathrm{sup} \{ {\mathscr{T}_j } \}      \gtrsim  \mathrm{log} [ \mathscr{L}^2 ]              \big\}}                                                                  \bigg\}   \\ \times  \underset{\forall i , c_i > 0 }{\underset{1 \leq i \leq \mathscr{L}}{\prod}}     \mathscr{E} \bigg\{         1   +   \bigg[    \underset{\forall i , c_i > 0 }{\underset{1 \leq i \leq \mathscr{L}_1, 1 \leq j \leq \mathscr{L}_2}{\prod}}          \widetilde{\mathscr{T}_i }     \mathscr{T}_j      \bigg]    \textbf{1}_{\big\{               i :     \mathrm{sup} \{ \widetilde{\mathscr{T}_i } \}      \gtrsim  \mathrm{log} [ \mathscr{L}^2 ]              \big\}}      \textbf{1}_{\big\{               j :     \mathrm{sup} \{ {\mathscr{T}_j } \}      \gtrsim  \mathrm{log} [ \mathscr{L}^2 ]              \big\}}                                                                                                                      \bigg\}                                      \bigg\}      .       \tag{$\mathscr{L}^2-\mathrm{Exp}$}   \end{align*}}

\noindent Proceeding, further rearrangements imply,

      {\tiny   \begin{align*}                           (\mathscr{L}^2-\mathrm{Exp}) < \underset{\forall i, c_i > 0 }{\underset{1 \leq i \leq \mathscr{L}}{\mathrm{sup}}} 2  c_i              \mathscr{L}^2      \bigg\{                       \mathscr{E} \bigg\{            1  +   \bigg[    \underset{\forall i , c_i > 0 }{\underset{1 \leq i \leq \mathscr{L}}{\prod}}       \widetilde{\mathscr{T}_i }          \bigg]       \textbf{1}_{\big\{               i :     \mathrm{sup} \{ \widetilde{\mathscr{T}_i } \}      \gtrsim  \mathrm{log} [ \mathscr{L}^2 ]              \big\}}      \textbf{1}_{\big\{               j :     \mathrm{sup} \{ {\mathscr{T}_j } \}      \gtrsim  \mathrm{log} [ \mathscr{L}^2 ]              \big\}}                                                          \bigg\}                                      \bigg\}     \\ \\  <   \mathscr{L}^2     \mathrm{exp} \bigg[   -    \underset{\forall i, c_i > 0 }{\underset{1 \leq i \leq \mathscr{L}}{\mathrm{sup}}}   \frac{c_i             }{100}      \bigg\{          1  +   \bigg[    \underset{\forall i , c_i > 0 }{\underset{1 \leq i \leq \mathscr{L}}{\prod}}       \widetilde{\mathscr{T}_i }          \bigg]       \textbf{1}_{\big\{               i :     \mathrm{sup} \{ \widetilde{\mathscr{T}_i } \}      \gtrsim  \mathrm{log} [ \mathscr{L}^2 ]              \big\}}      \textbf{1}_{\big\{               j :     \mathrm{sup} \{ {\mathscr{T}_j } \}      \gtrsim  \mathrm{log} [ \mathscr{L}^2 ]              \big\}}                                                          \bigg\}                                      \bigg]  \\ \\ = \mathscr{L}^2   \mathrm{exp} \bigg[   -    \underset{\forall i, c_i > 0 }{\underset{1 \leq i \leq \mathscr{L}}{\mathrm{sup}}}   \frac{c_i            }{100}      \bigg\{          1  +   \bigg[    \underset{\forall i , c_i > 0 }{\underset{1 \leq i \leq \mathscr{L}}{\prod}}       \widetilde{\mathscr{T}_i }          \bigg]       \textbf{1}_{\big\{               i :     \mathrm{sup} \{ \widetilde{\mathscr{T}_i } \}      \gtrsim  \mathrm{log} [ \mathscr{L}^2 ]              \big\}}      \textbf{1}_{\big\{               j :     \mathrm{sup} \{ {\mathscr{T}_j } \}      \gtrsim  \mathrm{log} [ \mathscr{L}^2 ]              \big\}}                                                          \bigg\}                            +      \underset{\forall i, c_i > 0 }{\underset{1 \leq i \leq \mathscr{L}}{\mathrm{sup}}}   \frac{c_i            }{100}   -   \underset{\forall i, c_i > 0 }{\underset{1 \leq i \leq \mathscr{L}}{\mathrm{sup}}}   \frac{c_i             }{100}          \bigg] \end{align*}

      \begin{align*}  =  \mathscr{L}^2     \mathrm{exp} \bigg[    -    \bigg[    \underset{\forall i , c_i > 0 }{\underset{1 \leq i \leq \mathscr{L}}{\prod}}       \widetilde{\mathscr{T}_i }          \bigg]       \textbf{1}_{\big\{               i :     \mathrm{sup} \{ \widetilde{\mathscr{T}_i } \}      \gtrsim  \mathrm{log} [ \mathscr{L}^2 ]              \big\}}      \textbf{1}_{\big\{               j :     \mathrm{sup} \{ {\mathscr{T}_j } \}      \gtrsim  \mathrm{log} [ \mathscr{L}^2 ]              \big\}}                                                                              +      \underset{\forall i, c_i > 0 }{\underset{1 \leq i \leq \mathscr{L}}{\mathrm{sup}}}   \frac{c_i            }{100}     \bigg] \\ \\   =    \mathscr{L}^2   \mathrm{exp} \bigg[    -    \bigg[    \underset{\forall i , c_i > 0 }{\underset{1 \leq i \leq \mathscr{L}}{\prod}}       \widetilde{\mathscr{T}_i }          \bigg]         \textbf{1}_{\big\{               i :     \mathrm{sup} \{ \widetilde{\mathscr{T}_i } \}      \gtrsim  \mathrm{log} [ \mathscr{L}^2 ]              \big\}}      \textbf{1}_{\big\{               j :     \mathrm{sup} \{ {\mathscr{T}_j } \}      \gtrsim  \mathrm{log} [ \mathscr{L}^2 ]              \big\}}                    \bigg]  \underset{\forall i, c_i > 0 }{\underset{1 \leq i \leq \mathscr{L}}{\mathrm{sup}}}     \mathrm{exp} \bigg[     \frac{c_i            }{100}              \bigg]                       \\ \\  <    \mathscr{L}^3  \underset{\forall i, c_i > 0 }{\underset{1 \leq i \leq \mathscr{L}}{\mathrm{sup}}}     \mathrm{exp} \bigg[    -    \frac{c_i            }{100}        \widetilde{\mathscr{T}_i }         \textbf{1}_{\big\{               i :     \mathrm{sup} \{ \widetilde{\mathscr{T}_i } \}      \gtrsim  \mathrm{log} [ \mathscr{L}^2 ]              \big\}}      \textbf{1}_{\big\{               j :     \mathrm{sup} \{ {\mathscr{T}_j } \}      \gtrsim  \mathrm{log} [ \mathscr{L}^2 ]              \big\}}                                                                               \bigg]     \\ \\    \lesssim                \mathscr{L}^3  \underset{\forall i, c_i > 0 }{\underset{1 \leq i \leq \mathscr{L}}{\mathrm{sup}}}     \mathrm{exp} \bigg[    -    \frac{c_i            }{100}        \widetilde{\mathscr{T}_i }         \textbf{1}_{\big\{               i :     \mathrm{sup} \{ \widetilde{\mathscr{T}_i } \}      \gtrsim  \mathrm{log} [ \mathscr{L}^2 ]              \big\}}      \textbf{1}_{\big\{               j :     \mathrm{sup} \{ {\mathscr{T}_j } \}      \gtrsim  \mathrm{log} [ \mathscr{L}^2 ]              \big\}}                                                                               \bigg]                                             .
\end{align*} }

\noindent Hence the up to constants estimate,

\begin{align*}
   \mathscr{E} \bigg\{     \underset{1 \leq i \leq \gamma_{\mathscr{L}}}{\prod}                   \frac{\mathscr{L}_{R_{i-1}}}{\hat{\mathscr{L}_{R_{i-1}}}}          \textbf{1}_{\big\{ \mathscr{L } - \underset{1 \leq  k \leq \gamma_{\mathscr{L}}}{\sum} \mathscr{T}_k \gtrsim \mathrm{log} [ \mathscr{L}^2 ] \big\}}            \bigg\}    \lesssim  \mathscr{L}^3 \mathrm{exp} \bigg[    -      \frac{c_i            }{100}                 \alpha \mathrm{log} \big[       \mathscr{L}^2   \big]    \bigg]          , 
\end{align*}

\noindent implies the desired up to constant exponential upper bound, from the observation that,

{\small \begin{align*}
     \frac{  \mathscr{E} \bigg\{      \underset{1 \leq i \leq \gamma_{\mathscr{L}}}{\prod}                   \frac{\mathscr{L}_{R_{i-1}}}{\hat{\mathscr{L}_{R_{i-1}}}}          \textbf{1}_{\big\{ \mathscr{L } - \underset{1 \leq  k \leq \gamma_{\mathscr{L}}}{\sum} \mathscr{T}_k \gtrsim \mathrm{log} [ \mathscr{L}^2 ] \big\}} \bigg\}  }{  \mathscr{E} \bigg\{      \underset{1 \leq i \leq \gamma_{\mathscr{L}}}{\prod}    \frac{\mathscr{L}_{R_{i-1}}}{\hat{\mathscr{L}_{R_{i-1}}}}          \textbf{1}_{\big\{ \underset{1 \leq k \leq \gamma_{\mathscr{L}_1}}{\sum} \mathscr{T}_k = \mathscr{L} \big\}}                              \bigg\}  }     , 
\end{align*} }

\noindent can be used to upper bound the probability,

{\small \begin{align*}
\underset{i \longrightarrow + \infty}{\underset{i^{\prime} \longrightarrow + \infty}{\underset{i^{\prime\prime} \longrightarrow + \infty}{\mathrm{lim}}}}   \textbf{P}_{\mathscr{U}, \textbf{T}}  \bigg[                                \forall \big( i , i^{\prime} , i^{\prime} \big) \in \textbf{N} \times \textbf{N} \times \textbf{N} , \exists  \big[   \gamma_1  ,   \gamma_2 ,  \gamma_3   \big]                 : \big[  \big\{ \textbf{1}^{\prime} = 1  \big\} \textit{ using $\gamma_1$} , \big\{ \textbf{2}^{\prime} = 1 \big\}  \textit{ using $\gamma_2$} , \big\{  \textbf{1}^{\prime} = 1 \textit{ using $\gamma_3$} \big\}   \big]  \\ = \big[ 1 , 1 , 1 \big]       \bigg] ,
\end{align*} }

\noindent for,

{\small \begin{align*}
  \gamma_1 : =  \underset{i \leq e_1 \leq 1}{\bigcup} \gamma_{e_1 }  , \\ \\  \gamma_2 : =  \underset{1 \leq e_2 \leq 2}{\bigcup} \gamma_{e_2 } , \\ \\   \gamma_3 : =  \underset{1 \leq e_3 \leq 3}{\bigcup} \gamma_{e_3 }  
\end{align*} }

\noindent from the following observation. Generally, one would like to compute, with respect to $\mathscr{E} \big[ \cdot \big]$ and a function of the times $\mathscr{T}_1, \cdots, \mathscr{T}_{\mathscr{L}_1}, \widetilde{\mathscr{T}_1}, \cdots, \widetilde{\mathscr{T}_{\mathscr{L}_2}}$, 

\begin{align*}
 \textit{functions} : \big[ \mathscr{T}_1 , \cdots, \mathscr{T}_{\mathscr{L}_1} ,  \widetilde{\mathscr{T}_1} , \cdots , \widetilde{\mathscr{T}_{\mathscr{L}_2}}       \big] \longrightarrow \textbf{T}    , 
\end{align*}

\noindent an expectation of the form,

{\small \begin{align*}
       \mathscr{E} \bigg\{      \underset{1 \leq i \leq \gamma_{\mathscr{L}}}{\prod}                   \frac{\mathscr{L}_{R_{i-1}}}{\hat{\mathscr{L}_{R_{i-1}}}}           \text{ } \textit{functions} \bigg[   \mathscr{T}_1, \cdots, \mathscr{T}_{\mathscr{L}_1}, \widetilde{\mathscr{T}_1}, \cdots, \widetilde{\mathscr{T}_{\mathscr{L}_2}} \bigg]    \bigg\}       , 
\end{align*} } 

\noindent admits the following decomposition,

{\small \begin{align*}
    \mathscr{E} \bigg\{      \underset{1 \leq i \leq \gamma_{\mathscr{L}_1}}{\prod}   \bigg\{      \underset{1 \leq j \leq \gamma_{\mathscr{L}_2}}{\prod}                 \frac{\mathscr{L}_{R_{i-1}}}{\hat{\mathscr{L}_{R_{i-1}}}}   \times      \frac{\mathscr{L}_{R_{j-1}}}{\hat{\mathscr{L}_{R_{j-1}}}}         \bigg\}         \text{ } \textit{functions} \bigg[   \mathscr{T}_1, \cdots, \mathscr{T}_{\mathscr{L}_1}, \widetilde{\mathscr{T}_1}, \cdots, \widetilde{\mathscr{T}_{\mathscr{L}_2}} \bigg]    \bigg\}  \\ \\ =             \mathscr{E} \bigg\{   \underset{1 \leq k \leq \mathscr{L}_1 }{\sum} \bigg\{   \underset{1 \leq i \leq \gamma_{k}}{\prod}   \bigg\{      \underset{1 \leq j \leq \gamma_{\mathscr{L}_2}}{\prod}                 \frac{\mathscr{L}_{R_{i-1}}}{\hat{\mathscr{L}_{R_{i-1}}}}   \times      \frac{\mathscr{L}_{R_{j-1}}}{\hat{\mathscr{L}_{R_{j-1}}}}         \bigg\}        \bigg\}  \text{ } \textit{functions} \bigg[   \mathscr{T}_1, \cdots, \mathscr{T}_{\mathscr{L}_1}, \widetilde{\mathscr{T}_1}, \\ \cdots, \widetilde{\mathscr{T}_{\mathscr{L}_2}} \bigg]    \bigg\} \\ \\ =    \mathscr{E} \bigg\{  \bigg\{ \underset{1 \leq k^{\prime} \leq \mathscr{L}_2 }{\sum}  \bigg\{  \underset{1 \leq k \leq \mathscr{L}_1 }{\sum} \bigg\{   \underset{1 \leq i \leq \gamma_{k}}{\prod}   \bigg\{      \underset{1 \leq j \leq \gamma_{k^{\prime}}}{\prod}                 \frac{\mathscr{L}_{R_{i-1}}}{\hat{\mathscr{L}_{R_{i-1}}}}   \times      \frac{\mathscr{L}_{R_{j-1}}}{\hat{\mathscr{L}_{R_{j-1}}}}         \bigg\}  \bigg\}  \bigg\}       \text{ } \textit{functions} \bigg[   \mathscr{T}_1, \cdots,  \\ \mathscr{T}_{\mathscr{L}_1},  \widetilde{\mathscr{T}_1},  \cdots, \widetilde{\mathscr{T}_{\mathscr{L}_2}} \bigg]    \bigg\}  \\ \\ =   \mathscr{E} \bigg\{  \bigg\{ \underset{1 \leq k^{\prime} \leq \mathscr{L}_2 }{\sum}  \bigg\{  \underset{1 \leq k \leq \mathscr{L}_1 }{\sum} \bigg\{   \underset{1 \leq i \leq \gamma_{k}}{\prod}   \bigg\{      \underset{1 \leq j \leq \gamma_{k^{\prime}}}{\prod}    \mathrm{sup} \bigg[              \frac{\mathscr{L}_{R_{i-1}}}{\hat{\mathscr{L}_{R_{i-1}}}}   \times      \frac{\mathscr{L}_{R_{j-1}}}{\hat{\mathscr{L}_{R_{j-1}}}}   , 1   \bigg]     \bigg\}  \bigg\}  \bigg\}   *_{\mathscr{L}_1 \mathscr{L}_2 }   \bigg\}                   , 
\end{align*} } 

\noindent where,

{\small \begin{align*}
  *_{\mathscr{L}_1 \mathscr{L}_2 }  : =  \underset{1 \leq k^{\prime} \leq \mathscr{L}_2}{\underset{1 \leq k \leq \mathscr{L}_1 }{\sum}}  \mathscr{E} \bigg\{          \underset{1 \leq j \leq \gamma_{k^{\prime}}}{\underset{1 \leq i \leq \gamma_{k}}{\prod} }               \mathrm{sup} \bigg\{        \frac{\mathscr{L}_{R_{i-1}}}{\hat{\mathscr{L}_{R_{i-1}}}}   ,      \frac{\mathscr{L}_{R_{j-1}}}{\hat{\mathscr{L}_{R_{j-1}}}}    , 1    \bigg\}                                             \bigg|  \big\{    \textit{Filtration corresponding to } \mathscr{L}_1 - 1  \big\} \cup  \big\{ \textit{Filtration cor-} \\ \textit{responding to } \mathscr{L}_2 - 1   \big\}             \bigg\}                    . 
\end{align*} }

\noindent is introduced over a scale that is of strictly less than $1$ increment over $\mathscr{L}_1$ and also over $\mathscr{L}_2$. Moreover, from the stopping times $\tau$ and $\tau^{\prime}$,

\begin{align*}
   \tau : =   \underset{\textit{times $t$}}{\mathrm{sup}}   \big\{ t :  \textit{the triangular effective random walk hits the origin about the first degree of} \\ \textit{ freedom} \big\}    , \\ \\   \tau^{\prime} : =    \underset{\textit{times $t$}}{\mathrm{sup}}   \big\{ t : \textit{the triangular effective random walk hits the origin about the second degree} \\ \textit{ of freedom}     \big\}     .
\end{align*}

\noindent previously introduced for obtaining the decomposition,

\begin{align*}
      G^r_{T \subsetneq \textbf{T}} \big( \lambda \big)     \textbf{1}_{\{ R = R_1 : R_1 \text{ } \textit{spans}\text{ }  T \} }           +  G^r_{T^{\prime} \subsetneq      \textbf{T}} \big( \lambda \big)     \textbf{1}_{\{ R = R_2 : R_2 \text{ } \textit{spans} \text{ } T^{\prime} \} }  ,
\end{align*}

\noindent of the Green's function in \textbf{Lemma} \textit{1}, observe,

{\small \begin{align*}
\underset{\epsilon = 0, \textit{or } \epsilon =1}{\sum}  \bigg\{  \underset{\tau \leq n }{\sum} \bigg\{ \underset{\tau^{\prime} \leq n^{\prime}}{\sum}   \underset{i \neq j \in \textbf{N}}{\mathrm{sup}} \bigg[     \mathscr{L}_{R_{i-1}} \big[ n , n^{\prime} ,  \epsilon  , R_1 , R_2 \big]        ,    \hat{\mathscr{L}_{R_{i-1}}} \big[ n , n^{\prime} , \epsilon  , R_1 , R_2 \big]   ,     \mathscr{L}_{R_{j-1}} \big[ n , n^{\prime} ,  \epsilon  , R_1 , R_2 \big]      \\    ,      \hat{\mathscr{L}_{R_{j-1}}} \big[ n , n^{\prime} , \epsilon  , R_1 , R_2 \big]           \bigg]        \bigg\}     \bigg\}   \\ \\ <    \underset{\epsilon = 0, \textit{or } \epsilon =1}{\sum}  \bigg\{  \underset{\tau \leq n }{\sum} \bigg\{ \underset{\tau^{\prime} \leq n^{\prime}}{\sum}   \bigg[     \mathscr{L}_{R_{i-1}} \big[ n , n^{\prime} ,  \epsilon  , R_1 , R_2 \big]        +     \hat{\mathscr{L}_{R_{i-1}}} \big[ n , n^{\prime} , \epsilon  , R_1 , R_2 \big]   +      \mathscr{L}_{R_{j-1}} \big[ n , n^{\prime} ,  \epsilon  , R_1 , R_2 \big]      \\   +       \hat{\mathscr{L}_{R_{j-1}}} \big[ n , n^{\prime} , \epsilon  , R_1 , R_2 \big]           \bigg]        \bigg\}     \bigg\}      \\ \\    <   \underset{1 \leq j \leq \mathscr{L}_2}{\underset{1 \leq i \leq \mathscr{L}_1}{\sum}}    \bigg\{   \underset{\epsilon = 0, \textit{or } \epsilon =1}{\sum}  \bigg\{  \underset{\tau \leq n }{\sum} \bigg\{ \underset{\tau^{\prime} \leq n^{\prime}}{\sum}   \bigg[     \mathscr{L}_{R_{i-1}} \big[ n , n^{\prime} ,  \epsilon  , R_1 , R_2 \big]        +     \hat{\mathscr{L}_{R_{i-1}}} \big[ n , n^{\prime} , \epsilon  , R_1 , R_2 \big]   +      \mathscr{L}_{R_{j-1}} \big[ n , n^{\prime} ,  \epsilon       \\      , R_1 , R_2 \big]          +       \hat{\mathscr{L}_{R_{j-1}}} \big[ n , n^{\prime} , \epsilon  , R_1 , R_2 \big]           \bigg]        \bigg\}     \bigg\}        \bigg\}    \\ \\            \\ \\ =  \underset{1 \leq i \leq \mathscr{L}_1}{\sum}  \bigg\{     \underset{\epsilon = 0, \textit{or } \epsilon =1}{\sum}  \bigg\{  \underset{\tau \leq n }{\sum} \bigg\{  \underset{\tau^{\prime} \leq n^{\prime}}{\sum}       \mathscr{L}_{R_{i-1}} \big[ n , n^{\prime} ,  \epsilon  , R_1 , R_2 \big]     \bigg\} \bigg\}  \bigg\}   +    \underset{1 \leq i \leq \mathscr{L}_1}{\sum}  \bigg\{     \underset{\epsilon = 0, \textit{or } \epsilon =1}{\sum}  \bigg\{  \underset{\tau \leq n }{\sum} \bigg\{  \underset{\tau^{\prime} \leq n^{\prime}}{\sum}       \hat{\mathscr{L}_{R_{i-1}}} \big[ n , \\ n^{\prime} , \epsilon  , R_1 , R_2 \big] \bigg\} \bigg\} \bigg\}     +     \underset{1 \leq j \leq \mathscr{L}_2}{\sum}  \bigg\{     \underset{\epsilon = 0, \textit{or } \epsilon =1}{\sum}  \bigg\{  \underset{\tau \leq n }{\sum} \bigg\{  \underset{\tau^{\prime} \leq n^{\prime}}{\sum}       \mathscr{L}_{R_{j-1}} \big[ n , n^{\prime} ,  \epsilon   , R_1 , R_2 \big]  \bigg\} \bigg\} \bigg\}     \\  +     \underset{1 \leq j \leq \mathscr{L}_2}{\sum}  \bigg\{     \underset{\epsilon = 0, \textit{or } \epsilon =1}{\sum}  \bigg\{  \underset{\tau \leq n }{\sum} \bigg\{  \underset{\tau^{\prime} \leq n^{\prime}}{\sum}       \hat{\mathscr{L}_{R_{j-1}}} \big[ n ,  n^{\prime} , \epsilon  , R_1 , R_2 \big]    \bigg\} \bigg\} \bigg\}    \\ \\    =      \underset{1 \leq i \leq \mathscr{L}_1}{\sum} \bigg[  \bigg\{     \underset{\epsilon = 0, \textit{or } \epsilon =1}{\sum}  \bigg\{  \underset{\tau \leq n }{\sum} \bigg\{  \underset{\tau^{\prime} \leq n^{\prime}}{\sum}       \mathscr{L}_{R_{i-1}} \big[ n , n^{\prime} ,  \epsilon  , R_1 , R_2 \big]     \bigg\} \bigg\}  \bigg\}   +   \bigg\{     \underset{\epsilon = 0, \textit{or } \epsilon =1}{\sum}  \bigg\{  \underset{\tau \leq n }{\sum} \bigg\{  \underset{\tau^{\prime} \leq n^{\prime}}{\sum}       \hat{\mathscr{L}_{R_{i-1}}} \big[ n ,  n^{\prime} , \epsilon    \end{align*}

\begin{align*}    , R_1 , R_2 \big] \bigg\} \bigg\} \bigg\}  \bigg] +     \underset{1 \leq j \leq \mathscr{L}_2}{\sum} \bigg[  \bigg\{     \underset{\epsilon = 0, \textit{or } \epsilon =1}{\sum}  \bigg\{  \underset{\tau \leq n }{\sum} \bigg\{  \underset{\tau^{\prime} \leq n^{\prime}}{\sum}       \mathscr{L}_{R_{j-1}} \big[ n , n^{\prime} ,  \epsilon   , R_1 , R_2 \big]  \bigg\} \bigg\} \bigg\}  \\   +      \bigg\{     \underset{\epsilon = 0, \textit{or } \epsilon =1}{\sum}  \bigg\{  \underset{\tau \leq n }{\sum} \bigg\{  \underset{\tau^{\prime} \leq n^{\prime}}{\sum}       \hat{\mathscr{L}_{R_{j-1}}} \big[ n ,  n^{\prime} , \epsilon  , R_1 , R_2 \big]    \bigg\} \bigg\} \bigg\}     \bigg]                       \end{align*}

\begin{align*} <   \bigg\{  1 + {\underset{1 \leq i \leq \mathscr{L}_1}{\prod}} c_i \widetilde{\mathscr{T}_i} \bigg\}   \hat{\mathscr{L}_{R_{i-1}}}      + \bigg\{  1 +                  {\underset{1 \leq j \leq \mathscr{L}_2}{\prod}}  c_j  \mathscr{T}_j      \bigg\}         \hat{\mathscr{L}_{R_{j-1}}}            \\ \\ <   \frac{\mathscr{C}}{2}  \bigg\{   1 + {\underset{1 \leq i \leq \mathscr{L}_1}{\prod}} c_i \widetilde{\mathscr{T}_i}        +  1 +                  {\underset{1 \leq j \leq \mathscr{L}_2}{\prod}}  c_j  \mathscr{T}_j      \bigg\}     \hat{\mathscr{L}_{R_{i-1}}}     \hat{\mathscr{L}_{R_{j-1}}}          \\ \\  <    \mathscr{C}    \bigg\{ 1 + \underset{1 \leq j \leq \mathscr{L}_2}{\underset{1 \leq i \leq \mathscr{L}_1}{\prod}} c_i c_j \widetilde{\mathscr{T}_i} \mathscr{T}_j  \bigg\}    \hat{\mathscr{L}_{R_{i-1}}}    \hat{\mathscr{L}_{R_{j-1}}}                               \\ \\ <      \mathscr{C}  \mathscr{L}^2   \bigg\{ 1 + \underset{1 \leq j \leq \mathscr{L}_2}{\underset{1 \leq i \leq \mathscr{L}_1}{\prod}} c_i c_j \widetilde{\mathscr{T}_i} \mathscr{T}_j  \bigg\}    \hat{\mathscr{L}_{R_{i-1}}}    \hat{\mathscr{L}_{R_{j-1}}}          \\ \\  <            \mathscr{C}  \mathscr{L}^2   \mathrm{exp}  \bigg\{ -   1 -  \underset{1 \leq j \leq \mathscr{L}_2}{\underset{1 \leq i \leq \mathscr{L}_1}{\prod}} c_i c_j \widetilde{\mathscr{T}_i} \mathscr{T}_j  \bigg\}    \hat{\mathscr{L}_{R_{i-1}}}    \hat{\mathscr{L}_{R_{j-1}}}  \\ \\ <     \mathscr{C}  \mathscr{L}^2   \mathrm{exp}  \bigg\{ \bigg\{ -   1 -  \underset{1 \leq j \leq \mathscr{L}_2}{\underset{1 \leq i \leq \mathscr{L}_1}{\prod}} c_i c_j \widetilde{\mathscr{T}_i} \mathscr{T}_j  \bigg\} \big[ \underset{i,j}{\mathrm{sup}}    \big\{ \hat{\mathscr{L}_{R_{i-1}}} \\ ,    \hat{\mathscr{L}_{R_{j-1}}} \big\} \big]^2  \bigg\}     \\ \\ <        \mathscr{C}  \mathscr{L}^4   \mathrm{exp}  \bigg\{ -   1 -  \underset{1 \leq j \leq \mathscr{L}_2}{\underset{1 \leq i \leq \mathscr{L}_1}{\prod}} c_i c_j \widetilde{\mathscr{T}_i} \mathscr{T}_j  \bigg\}    \\ \\ <      \mathscr{C}_1   \mathscr{L}^4   \mathrm{exp}  \bigg\{ -  \mathscr{C}_2 \bigg\{      1 +   \underset{1 \leq j \leq \mathscr{L}_2}{\underset{1 \leq i \leq \mathscr{L}_1}{\prod}} \widetilde{\mathscr{T}_i} \mathscr{T}_j \bigg\}    \bigg\}       , 
\end{align*} }

\noindent for $\mathscr{C}$ taken to be sufficiently large. Hence, from the above manipulations of the expected value taken over the ratios $\frac{\mathscr{L}_{R_{i-1}}}{\hat{\mathscr{L}_{R_{i-1}}}}$, and $\frac{\mathscr{L}_{R_{j-1}}}{\hat{\mathscr{L}_{R_{j-1}}}}$, for $1 \leq i \leq \mathscr{L}_1$ and $1 \leq j \leq \mathscr{L}_2$ respectively, one has that,

{\small \begin{align*}
      \underset{\epsilon = 0, \textit{or } \epsilon =1}{\sum}  \bigg\{  \underset{\tau \leq n }{\sum} \bigg\{ \underset{\tau^{\prime} \leq n^{\prime}}{\sum}   \underset{i \neq j \in \textbf{N}}{\mathrm{sup}} \bigg[     \mathscr{L}_{R_{i-1}} \big[ n , n^{\prime} ,  \epsilon  , R_1 , R_2 \big]        ,    \hat{\mathscr{L}_{R_{i-1}}} \big[ n , n^{\prime} , \epsilon  , R_1 , R_2 \big]   ,     \mathscr{L}_{R_{j-1}} \big[ n , n^{\prime} ,  \epsilon  , R_1 , R_2 \big]      \\    ,      \hat{\mathscr{L}_{R_{j-1}}} \big[ n , n^{\prime} , \epsilon  , R_1 , R_2 \big]           \bigg]        \bigg\}     \bigg\}  <           \mathscr{C}_1 \mathscr{L}^4  \mathrm{exp} \bigg\{            -            \mathscr{C}_2 \mathrm{log} \big[ \mathscr{L}^2   \big]          \bigg\}                     , 
\end{align*} }

\noindent for $\mathscr{C}_2 > \mathscr{C}_1$ taken sufficiently large, with $\mathscr{C}_2 > \mathscr{C}$ and $\mathscr{C}_1 > \mathscr{C}$. As $\mathscr{L} \longrightarrow + \infty$,

\begin{align*}
\big\{ - \mathrm{log} \big[ \mathscr{L}^2 \big] \longrightarrow 0 \big\}  \Rightarrow \bigg\{ \mathrm{exp} \bigg\{            -            \mathscr{C}_2 \mathrm{log} \big[ \mathscr{L}^2   \big]          \bigg\}     \longrightarrow 0  \bigg\}                , 
\end{align*}

\noindent from which we conclude the argument. \boxed{}

\bigskip

\noindent \textbf{Lemma} \textit{4} (\textit{sending the scaling limit normalization to positive infinity, under a $\sqrt{3}$ normalization}). One has that,

   \[    \left\{\!\begin{array}{ll@{}>{{}}l}  \frac{1}{\sqrt{3} }                         \textbf{P}_{\mathscr{U},\textbf{T}} \big[   \exists i \in \big\{ \mathrm{log} \big[ \mathscr{L}^2 \big], \cdots ,  \gamma_{e_1} \big\}    : \textbf{1}^{\prime} = 1 \textit{ for every i}  \big] \longrightarrow 0 \textit{ as } \mathscr{L} \longrightarrow + \infty  ,   \\ \\ \frac{1}{\sqrt{3} }    \textbf{P}_{\mathscr{U},\textbf{T}} \big[   \exists i \in \big\{ \mathrm{log} \big[ \mathscr{L}^2 \big] , \cdots ,  \gamma_{e_2} \big\}    : \textbf{2}^{\prime} = 1   \textit{ for every i} \big] \longrightarrow 0 \textit{ as } \mathscr{L} \longrightarrow + \infty ,   \\ \\  \frac{1}{\sqrt{3} }    \textbf{P}_{\mathscr{U},\textbf{T}} \big[    \exists i \in \big\{ \mathrm{log} \big[ \mathscr{L}^2 \big] , \cdots ,  \gamma_{e_3} \big\}     : \textbf{3}^{\prime} = 1  \textit{ for every i} \big]  \bigg] \longrightarrow 0 \textit{ as } \mathscr{L} \longrightarrow + \infty ,     \end{array}\right. \]

\noindent holds iff,

    \begin{align*}
     \frac{1}{\sqrt{3} }   \bigg\{  \underset{\mathscr{L} \longrightarrow + \infty}{\mathrm{lim}}            \bigg[           \textbf{P}_{\mathscr{U},\textbf{T}} \big[   \exists i \in \big\{ \mathrm{log} \big[ \mathscr{L}^2 \big], \cdots ,  \gamma_{e_1} \big\}    : \textbf{1}^{\prime} = 1 \textit{ for every i}  \big]        + \textbf{P}_{\mathscr{U},\textbf{T}} \big[   \exists i \in \big\{ \mathrm{log} \big[ \mathscr{L}^2 \big] , \cdots ,  \gamma_{e_2} \big\}   \\  : \textbf{2}^{\prime} = 1   \textit{ for every i} \big]     + \textbf{P}_{\mathscr{U},\textbf{T}} \big[    \exists i \in \big\{ \mathrm{log} \big[ \mathscr{L}^2 \big] , \cdots ,  \gamma_{e_3} \big\}     : \textbf{3}^{\prime} = 1  \textit{ for every i} \big]  \bigg] \bigg\} \\ \longrightarrow 0  ,
    \end{align*}

\noindent holds.

\bigskip

\noindent \textbf{Lemma} \textit{5} (\textit{upper bounding the probability that an excursion of the triangular effective random walk has at most $\mathrm{log} \big[ \mathscr{L}^2 \big]$ many excursions about each degree of freedom}). Denote,

\begin{align*}
  \mathcal{T} : = \mathcal{T}_1 + \cdots + \mathcal{T}_{\gamma_{\mathscr{L}_1}}   =      \textit{Time spent by an excursion of the effective triangular} \\ \textit{ prudent walk along the first degree of freedom}     , \\ \\   \widetilde{\mathcal{T}}  : =  \widetilde{\mathcal{T}_1} + \cdots + \widetilde{\mathcal{T}_{\gamma_{\mathscr{L}_2}}}   =            \textit{Time spent by an excursion of the effective triangular} \\ \textit{ prudent walk along the second degree of freedom}     , \\ \\ \mathscr{L}_1 : =  \textit{First component of the scaling limit normalization $\mathscr{L}$} , \\ \\ \mathscr{L}_2 : = \textit{Second component of the scaling limit normalization $\mathscr{L}$}    .  
\end{align*}

\noindent One has that,

\begin{align*}
         \mathscr{E} \bigg\{     \underset{1 \leq i \leq \gamma_{\mathscr{L}}}{\prod}                   \frac{\mathscr{L}_{R_{i-1}}}{\hat{\mathscr{L}_{R_{i-1}}}}          \textbf{1}_{\big\{  \mathcal{T} + \widetilde{\mathcal{T}} \geq \mathrm{log} [ \mathscr{L}^4  ]     \big\}}          \bigg\}  , \end{align*}

         \noindent implies, 
         \begin{align*} \underset{\mathscr{L}_2 \longrightarrow + \infty}{\mathrm{lim}}     \text{ }        \underset{\mathscr{L}_1 \longrightarrow + \infty}{\mathrm{lim}}                       \textbf{P}_{\mathscr{U}, \textbf{T}} \big[      \forall           \mathscr{L}_1 , \mathscr{L}_2 , \exists \big\{ \alpha_1 \neq \alpha_2 > 0 , EW \sim \mathscr{E} \mathscr{W}_{\textbf{T}} \big\}  :   \# \big\{  \textit{edges} : \textit{edges belong to } EW \big|_1           \big\} \\ +  \# \big\{   \textit{edges}  : \textit{edges belong to } EW \big|_2           \big\}          \geq \mathrm{log} \big[    \mathscr{L}^4   \big]                     \big]                                                 . 
\end{align*}

\noindent \textit{Proof of Lemma 4 and of Lemma 5}. We adapt the arguments provided in the previous result. That is, we make use of identical computations as provided in the previous result, adjusting the constants appearing in the scaling limit constant $\mathscr{L}$, as well as each $\mathscr{T}_i$. Altogether this implies the desired up to constant exponential upper bound, from the fact that,

{\small \begin{align*}
     \frac{  \mathscr{E} \bigg\{      \underset{1 \leq i \leq \gamma_{\mathscr{L}}}{\prod}                   \frac{\mathscr{L}_{R_{i-1}}}{\hat{\mathscr{L}_{R_{i-1}}}}          \textbf{1}_{\big\{ \mathscr{L }_1 - \underset{1 \leq  k \leq \gamma_{\mathscr{L}_1}}{\sum} \mathcal{T}_k + \mathscr{L }_2 - \underset{1 \leq  k \leq \gamma_{\mathscr{L}_2}}{\sum} \widetilde{\mathcal{T}_k}   \gtrsim \mathrm{log} [ \mathscr{L}^2 ] \big\}} \bigg\}  }{  \mathscr{E} \bigg\{      \underset{1 \leq i \leq \gamma_{\mathscr{L}}}{\prod}    \frac{\mathscr{L}_{R_{i-1}}}{\hat{\mathscr{L}_{R_{i-1}}}}          \textbf{1}_{\big\{ \underset{1 \leq k \leq \gamma_{\mathscr{L}_1}}{\sum} \mathcal{T}_k  +  \underset{1 \leq k \leq \gamma_{\mathscr{L}_2}}{\sum}\widetilde{\mathcal{T}_k} = \mathscr{L}_1 + \mathscr{L}_2  \big\}}                              \bigg\}  }   < \mathscr{C}^{\prime}_1 \mathscr{L}^8 \mathrm{exp} \bigg\{            -            \mathscr{C}^{\prime}_2 \mathrm{log} \big[ \mathscr{L}^4   \big]          \bigg\}        , 
\end{align*} } 

\noindent for $\mathscr{C}^{\prime}_2 > \mathscr{C}^{\prime}_1 > \mathscr{C}^{\prime} > 0$ taken sufficiently large, and also that,

{\small \begin{align*}
     \frac{  \mathscr{E} \bigg\{      \underset{1 \leq i \leq \gamma_{\mathscr{L}}}{\prod}                   \frac{\mathscr{L}_{R_{i-1}}}{\hat{\mathscr{L}_{R_{i-1}}}}          \textbf{1}_{\big\{ \mathscr{L }_1 - \underset{1 \leq  k \leq \gamma_{\mathscr{L}_1}}{\sum} \mathcal{T}_k + \mathscr{L }_2 - \underset{1 \leq  k \leq \gamma_{\mathscr{L}_2}}{\sum} \widetilde{\mathcal{T}_k}   \gtrsim \mathrm{log} [ \mathscr{L}^2 ] \big\}} \bigg\}  }{  \mathscr{E} \bigg\{      \underset{1 \leq i \leq \gamma_{\mathscr{L}}}{\prod}    \frac{\mathscr{L}_{R_{i-1}}}{\hat{\mathscr{L}_{R_{i-1}}}}          \textbf{1}_{\big\{ \underset{1 \leq k \leq \gamma_{\mathscr{L}_1}}{\sum} \mathcal{T}_k  +  \underset{1 \leq k \leq \gamma_{\mathscr{L}_2}}{\sum}\widetilde{\mathcal{T}_k} = \mathscr{L}_1 + \mathscr{L}_2  \big\}}                              \bigg\}  }   < \mathscr{C}^{\prime\prime}_1 \mathscr{L}^{10} \mathrm{exp} \bigg\{            -            \mathscr{C}^{\prime\prime}_2 \mathrm{log} \big[ \mathscr{L}^6   \big]          \bigg\}        , 
\end{align*} }

\noindent for $\mathscr{C}^{\prime\prime}_2 > \mathscr{C}^{\prime\prime}_1 > \mathscr{C}^{\prime\prime} > 0$ taken sufficiently large. As $\mathscr{L} \longrightarrow + \infty$,

\begin{align*}
\mathscr{C}^{\prime}_1 \mathscr{L}^8 \mathrm{exp} \bigg\{            -            \mathscr{C}^{\prime}_2 \mathrm{log} \big[ \mathscr{L}^4   \big]          \bigg\}  \longrightarrow 0   , \\ \\  \mathscr{C}^{\prime\prime}_1 \mathscr{L}^{10} \mathrm{exp} \bigg\{            -            \mathscr{C}^{\prime\prime}_2 \mathrm{log} \big[ \mathscr{L}^6   \big]          \bigg\}     \longrightarrow 0   , 
\end{align*}

\noindent from which we conclude the argument. \boxed{}

\bigskip

\noindent We conclude the section by arguing that the desired result with the scaling limit over $\textbf{T}$ holds. The impact of the scaling limit normalization, $\mathscr{L}$ is provided in \textbf{Theorem}.

\bigskip

\noindent \textit{Proof of Theorem}. We demonstrate that the statement in $\textbf{Theorem}$ holds, namely that, with respect to a $\sqrt{\mathscr{L}}$ normalization, the prudent walk over $\textbf{T}$ converges to Brownian motion under diffusive rescaling. First, recall for $\widetilde{\pi_s} \equiv \widetilde{\pi_s} \big( \textbf{T} \big) $ and $\hat{\pi_s} \equiv \hat{\pi_s \big( \textbf{T} \big)}$, to argue that $\textit{(*)}$, one wants to show, under,

\begin{align*}
\textbf{P}_{\textbf{T}} \big[ \cdot \big]  : = \underset{T \in \textbf{T}}{\bigcup } \textbf{P}_{T} \big[ \cdot \big]  , 
\end{align*}

\noindent that,

\begin{align*}
   \bigg\{ \textbf{P}_{\mathscr{U},\textbf{T}}  \bigg[ \bigg\{   \text{ }     \underset{0 \leq s \leq t}{\mathrm{sup}} \bigg| \bigg|         \frac{1}{t} \big[    \hat{\pi_s} - \pi_s    \big]        \bigg|\bigg|_2 \geq \epsilon \bigg\} \text{ } \bigg| \text{ } \mathcal{Q}_1 \text{ }   \bigg]   \approx 1             \bigg\} \Longrightarrow      \bigg\{           \textbf{P}_{\mathscr{U},\textbf{T}}  \big[   \big|       \widetilde{\pi} t     -     \big[   c_{\{ 1,e \} } \textbf{1}_{    \{   e : e \in \mathrm{span} \{ \textit{first degree of freedom}        \}   \}   }  \\   +  c_{\{ 2,e \} } \textbf{1}_{    \{   e : e \in \mathrm{span} \{ \textit{second degree of freedom}        \}   \}   }     \big] t  \big|        \big]                                                                         \bigg\}        \\ \Downarrow \\     \bigg\{               \frac{1}{\sqrt{\mathscr{L}}}    \bigg\{         X_{\lfloor \mathscr{L} t \rfloor }                                       \bigg\}_{\{ t \geq 0\}  }                         \overset{t \longrightarrow + \infty}{\longrightarrow}   \big( \textit{Covariance matrix} \big)^{\frac{1}{2}}  \mathcal{B}^{(a)}  \bigg\}              ,  \tag{*}
\end{align*}

\noindent provided in the statement of  $\textbf{Theorem}$ holds, for,

\begin{align*}
   v_1 : =  \big[ 1 , 0 \big] , \\ \\    v_2 : =  \big[ \frac{1}{2} , \frac{\sqrt{3}}{2} \big]  , \\ \\   v_3 : = \big[ - \frac{1}{2} , \frac{\sqrt{3}}{2} \big]      , \\ \\   v_4 : = ,      \big[  - 1    , 0  \big]      \\ \\      v_5 : =    \big[ - \frac{1}{2} , - \frac{\sqrt{3}}{2}  \big]         , \\ \\    v_6 : =   \big[ \frac{1}{2} , - \frac{\sqrt{3}}{2} \big]       , 
\end{align*}

\noindent observe,

\begin{align*}
   \underset{1 \leq i \leq 6}{\sum} v_i v^{\mathrm{T}}_i = 3 \textbf{I}  , 
\end{align*}

\noindent implies,

\begin{align*}
  \textit{Covariance matrix} = \frac{1}{2} \textbf{I}   .
\end{align*}

\noindent Moreover, given the fact that increments of the Markov chain $X$ over $\textbf{T}$ have \textit{zero} drift, given by the condition,

\begin{align*}
   \mathscr{E} \big\{ X_{k+1} - X_k \big\} = 0  ,
\end{align*}

\noindent for some $k> 0$, in addition to the finite second moment condition,

\begin{align*}
  \mathscr{E} \big\{ \big(  X_{k+1} - X_k    \big) \big( X_{k+1} - X_k \big)^{\mathrm{T}} \big\} < + \infty   .
\end{align*}

\noindent In the forthcoming computations, we. make use of the set of linear combinations spanned by the basis,

\begin{align*}
\underset{1 \leq i \leq 6}{\mathrm{span}} \big\{ v_i \big\}  = \mathrm{span} \bigg\{ \big[ 1 , 0 \big]   , \big[ \frac{1}{2} , \frac{\sqrt{3}}{2} \big]  , \big[ - \frac{1}{2} , \frac{\sqrt{3}}{2} \big]  , \big[ - 1 , 0 \big] , \big[ - \frac{1}{2} , - \frac{\sqrt{3}}{2} \big] ,  \big[ \frac{1}{2} , - \frac{\sqrt{3}}{2} \big]        \bigg\} . 
\end{align*}

\noindent To argue that the following upper bound for an expected value. holds, we make use of the decomposition,

\begin{align*}
   \big\{  X_k \big\}_{k \geq 0} : = \big\{ M_k + R_k \big\}_{k \geq 0 } , 
\end{align*}

\noindent for the Markov chain. From \textbf{Definition} \textit{0}, recall that from the Markov chain associated with the simple random walk over $\textbf{T}$,

\begin{align*}
     \big\{ X_k \big\}_{k \geq 0 }   , 
\end{align*}

\noindent with the transition kernel,

\begin{align*}
  \textbf{P}_{\textbf{T}} \bigg[         X_{k+1} - X_k = a \mathrm{exp} \big[ \frac{ i \pi j}{3}              \big] \bigg] = \frac{1}{6}  , 
\end{align*}

\noindent supported over $\textbf{P}_{\textbf{T}} \big[ \cdot \big]$, under diffusive scaling with mesh size $a > 0$ can be used to obtain Brownian motion,

\begin{align*}
   \mathcal{B}^{(a)}_t : =  a X_{\lfloor   \frac{t}{a^2}  \rfloor }  . 
\end{align*}

\noindent which is supported over $\textbf{T}$. Fix some $C_1$ sufficiently large, as well as , with $C^{\prime}_1 \neq C_2 \neq C_3 > 0$. Hence it suffices to upper bound the expected value,

\begin{align*}
   \mathscr{E}  \bigg\{            \underset{k \leq t}{\mathrm{sup}} \big| \big| M_k \big| \big|^2  \bigg\}  \leq C_1   \mathscr{E} \big\{    \langle  M \rangle_t      \big\}            , 
\end{align*}

\noindent from the rearrangements,

\begin{align*}
       \mathscr{E}  \bigg\{            \underset{k \leq t}{\mathrm{sup}} \big| \big| X_k \big| \big|^2  \bigg\}  \leq 2      \mathscr{E}  \bigg\{            \underset{k \leq t}{\mathrm{sup}} \big| \big| M_k \big| \big|^2  \bigg\}  + 2      \mathscr{E}  \bigg\{            \underset{k \leq t}{\mathrm{sup}} \big| \big| R_k \big| \big|^2  \bigg\}  < 2 C^{\prime}_1 C_2 n + 2 C_3 \leq C_1 n   , 
\end{align*}

\noindent in turn implying,

\begin{align*}
            \frac{1}{t} \langle  M \rangle_t  \overset{\textbf{P}}{\longrightarrow} \textit{Covariance matrix} = \frac{1}{2 } \textbf{I}  .
\end{align*}

\bigskip

\noindent To argue that the sequence,

\begin{align*}
 \bigg\{    \frac{1}{\mathscr{L}} X_{\lfloor t \mathscr{L} \rfloor }   \bigg\}_{t \geq 0 }  \longrightarrow 0   , 
\end{align*}

\noindent can be manipulated to obtain the desired properties of the scaling limit over $\textbf{T}$, observe that it suffices to compute the covariance matrix, and mean, using $\mathcal{B}$ obtained through the above diffusive scaling. That is, with regards to the covariance matrix for mean zero Gaussian processes, write,

{\small \begin{align*}
  \mathrm{Cov} \bigg\{      \big\{ \mathcal{B}^{(a)}_t \big\}_{t \geq 0 } , \big\{ \mathcal{B}^{(a)}_{t^{\prime}}  \big\}_{t^{\prime} \geq 0 }    \bigg\}  =  \mathrm{Cov} \bigg\{      \big\{ \mathcal{B}^{(a)}_t \big\}, \big\{ \mathcal{B}^{(a)}_{t^{\prime}}  \big\}    \bigg\}_{ \{ t \neq t^{\prime} \geq 0  \} }    =   \mathscr{E}  \bigg\{   \big[   \mathcal{B}^{(a)}_t - \mathscr{E} \big\{ \mathcal{B}^{(a)}_t \big\}     \big] \big[  \mathcal{B}^{(a)}_{t^{\prime}}  - \mathscr{E} \big\{ \mathcal{B}^{(a)}_{t^{\prime}} \big\}             \big]       \bigg\}_{ \{ t \neq t^{\prime} \geq 0 \} }                   , 
\end{align*} }

\noindent which is equivalent to, by virtue of,

\begin{align*}
  \mathscr{E} \big\{ \mathcal{B}^{(a)}_{t} \big\}       = 0    , \\ \\  \mathscr{E} \big\{ \mathcal{B}^{(a)}_{t^{\prime}} \big\}       = 0 , 
\end{align*}

\noindent the following expectation,

\begin{align*}
  \mathscr{E} \bigg\{   \big\{ \mathcal{B}^{(a)}_t \big\}_{t \geq 0 } \big\{ \mathcal{B}^{(a)}_{t^{\prime}}  \big\}_{t^{\prime} \geq 0 }        \bigg\} =     \mathscr{E} \bigg\{  \mathcal{B}^{(a)}_t  \mathcal{B}^{(a)}_{t^{\prime}}        \bigg\}_{\{  t \neq t^{\prime} \geq 0\}  } =     a^2 \mathscr{E} \bigg\{    X_{\lfloor  \frac{a}{t^2 }      \rfloor }  X_{\lfloor     \frac{a}{( t^{\prime} )^2}    \rfloor }       \bigg\}_{ \{ t \neq t^{\prime} \geq 0 \} }               , \\   \tag{$a^2-\mathscr{E}$}
\end{align*}

\noindent from which further rearrangements imply,

\begin{align*}
  (a^2-\mathscr{E}) =   a^2 \textbf{I}_2 = a^2 \begin{bmatrix}
        \mathscr{E} \big\{         \big\{ \big( \mathcal{B}^{(1)}_t \big) \big\}^2           \big\}     &   \mathscr{E} \big\{         \big( \mathcal{B}^{(1)}_t \big) \big( \mathcal{B}^{(1)}_{t^{\prime}} \big)        \big\}           \\   \mathscr{E} \big\{ \big( \mathcal{B}^{(1)}_t \big) \big( \mathcal{B}^{(1)}_{t^{\prime}} \big)       \big\}       &       \mathscr{E} \big\{ \big( \mathcal{B}^{(1)}_t \big)^2 \big\}       
  \end{bmatrix} \\ \\ =    \begin{bmatrix}
     a^2 \frac{1}{6} \underset{0 \leq k \leq 5}{\sum}  \mathrm{cos}^2 \big[ \theta_k \big]    &    a^2 \frac{1}{6} \underset{0 \leq k \leq 5}{\sum}  \mathrm{sin}^2 \big[ \theta_k \big]       \\      a^2 \frac{1}{6} \underset{0 \leq k \leq 5}{\sum}  \mathrm{sin}^2 \big[ \theta_k \big]      &       a^2 \frac{1}{6} \underset{0 \leq k \leq 5}{\sum}  \mathrm{cos} \big[ \theta_k \big] \mathrm{sin} \big[ \theta_k \big] \end{bmatrix} \bigg|_{\theta_k = 0 } \\ \\ =  \begin{bmatrix}
   a^2    &  0  \\      0     &   a^2  \end{bmatrix} .
\end{align*}

\noindent Under the action of time reversal, $\mathscr{T}\mathscr{R}$, where,

\begin{align*}
 \mathscr{T}\mathscr{R} :                 \widetilde{\textit{Paths}}          \longrightarrow      \textit{Paths}    , 
\end{align*}

\noindent for,

\begin{align*}
     \widetilde{\textit{Paths}} : =    \underset{\textit{paths}}{\bigcup}          \bigg\{ \textit{paths} , \textit{paths}^{\prime} :  \big\{ \textit{paths } \subsetneq \textbf{T}  \big\} , \bigg\{   \mathrm{Cov} \big\{  \textit{paths} , \textit{paths}^{\prime}     \big\}   =              \mathrm{Cov} \bigg\{      \widetilde{\big\{ \mathcal{B}^{(a)}_t \big\} } , \widetilde{\big\{ \mathcal{B}^{(a)}_{t^{\prime}}  \big\}}    \bigg\}_{ \{ t \neq t^{\prime} \geq 0  \} }                                \bigg\}  \bigg\}      , \\ \\   \textit{Paths}  : =    \underset{\textit{paths}}{\bigcup}          \bigg\{ \textit{paths} , \textit{paths}^{\prime} :  \big\{ \textit{paths } \subsetneq \textbf{T}  \big\} , \bigg\{   \mathrm{Cov} \big\{  \textit{paths} , \textit{paths}^{\prime}     \big\}   =              \mathrm{Cov} \bigg\{     {\big\{ \mathcal{B}^{(a)}_t \big\} } , {\big\{ \mathcal{B}^{(a)}_{t^{\prime}}  \big\}}    \bigg\}_{ \{ t \neq t^{\prime} \geq 0  \} }                                \bigg\}  \bigg\}          , 
\end{align*}

\noindent the covariance,

\begin{align*}
  \widetilde{\mathrm{Cov} \bigg\{      \big\{ \mathcal{B}^{(a)}_t \big\}_{t \geq 0 } , \big\{ \mathcal{B}^{(a)}_{t^{\prime}}  \big\}_{t^{\prime} \geq 0 }    \bigg\}}  =  \mathrm{Cov} \bigg\{      \widetilde{\big\{ \mathcal{B}^{(a)}_t \big\} } , \widetilde{\big\{ \mathcal{B}^{(a)}_{t^{\prime}}  \big\}}    \bigg\}_{ \{ t \neq t^{\prime} \geq 0  \} }                     , 
\end{align*}

\noindent can be expressed as,

\begin{align*}
  \mathscr{E} \bigg\{   \widetilde{\big\{ \mathcal{B}^{(a)}_t \big\}_{t \geq 0 }}  \widetilde{\big\{ \mathcal{B}^{(a)}_{t^{\prime}}  \big\}_{t^{\prime} \geq 0 }}        \bigg\}  = \mathscr{E} \bigg\{    \big\{  \mathcal{B}^{(a)}_{\mathscr{T} \mathscr{R} ( t) } \big\}_{t \geq 0}      \big\{  \mathcal{B}^{(a)}_{\mathscr{T} \mathscr{R} ( t^{\prime}) } \big\}_{t^{\prime} \geq 0}                \bigg\}  . 
\end{align*}

\noindent The original, and time reversed, covariances can be related to one another through the following identity. It reads,

\begin{align*}
   \frac{\mathrm{Cov} \bigg\{      \big\{ \mathcal{B}^{(a)}_t \big\}_{t \geq 0 } , \big\{ \mathcal{B}^{(a)}_{t^{\prime}}  \big\}_{t^{\prime} \geq 0 }    \bigg\} }{\mathrm{Cov} \bigg\{      \widetilde{\big\{ \mathcal{B}^{(a)}_t \big\} } , \widetilde{\big\{ \mathcal{B}^{(a)}_{t^{\prime}}  \big\}}    \bigg\}_{ \{ t \neq t^{\prime} \geq 0  \} }  } =  \frac{ \begin{bmatrix}
     a^2 \frac{1}{6} \underset{0 \leq k \leq 5}{\sum}  \mathrm{cos}^2 \big[ \theta_k \big]    &    a^2 \frac{1}{6} \underset{0 \leq k \leq 5}{\sum}  \mathrm{sin}^2 \big[ \theta_k \big]       \\      a^2 \frac{1}{6} \underset{0 \leq k \leq 5}{\sum}  \mathrm{sin}^2 \big[ \theta_k \big]      &       a^2 \frac{1}{6} \underset{0 \leq k \leq 5}{\sum}  \mathrm{cos} \big[ \theta_k \big] \mathrm{sin} \big[ \theta_k \big] \end{bmatrix} \bigg|_{\theta_k = 0 } }{ \begin{bmatrix}
     a^2 \frac{1}{6} \underset{0 \leq k \leq 5}{\sum}  \mathrm{cos}^2 \big[ \widetilde{\theta_k} \big]    &    a^2 \frac{1}{6} \underset{0 \leq k \leq 5}{\sum}  \mathrm{sin}^2 \big[ \widetilde{\theta_k} \big]       \\      a^2 \frac{1}{6} \underset{0 \leq k \leq 5}{\sum}  \mathrm{sin}^2 \big[ \widetilde{\theta_k} \big]      &       a^2 \frac{1}{6} \underset{0 \leq k \leq 5}{\sum}  \mathrm{cos} \big[ \widetilde{\theta_k} \big] \mathrm{sin} \big[ \widetilde{\theta_k} \big] \end{bmatrix} \bigg|_{\widetilde{\theta_k} = 0 } }  \\ \\ =     1 =   \frac{ \mathscr{E} \bigg\{   \big\{    \mathcal{B}^{(a)}_t   \mathcal{B}^{(a)}_{t^{\prime}}   \big\}   \bigg\}_{ \{ t\neq t^{\prime} \geq 0  \} } } {\mathscr{E} \bigg\{  \big\{     \widetilde{\mathcal{B}^{(a)}_t}   \widetilde{\mathcal{B}^{(a)}_{t^{\prime}}} \big\}     \bigg\}_{ \{ t\neq t^{\prime} \geq 0 \}  }   }              . 
\end{align*}

\noindent Due to the fact that there is \textit{already} isotropy for the walk over $\textbf{T}$ for the first step, in comparison to obtaining isotropy through the diffusive scaling limit, isotropy is already satisfied by the above expression for the covariance. As a result the limit over $\textbf{T}$ is Brownian motion which is isotropic, is hence isotropic. Very similar arguments to those provided above demonstrate that

{\small \begin{align*}
    \bigg\{ \textbf{P}_{\mathscr{U},\textbf{T}}  \bigg[ \bigg\{   \text{ }     \underset{0 \leq s \leq t}{\mathrm{sup}} \bigg| \bigg|         \frac{1}{t} \big[    \hat{\pi_s} - \pi_s    \big]        \bigg|\bigg|_2 \geq \epsilon \bigg\} \text{ } \bigg| \text{ } \mathcal{Q}_1 \text{ }   \bigg]   \approx 1             \bigg\}  , 
\end{align*} } 

\noindent holds, from which we conclude the argument. \boxed{}

\section{Conclusion}

\noindent In this paper we demonstrated how the scaling limit of a walk defined over the triangular lattice can be obtained. Albeit the fact that aspects of the scaling limit share similarities with scaling limits that have previously been obtained over the square lattice, differences emerge. One of the most important differences relates to how the triangular Green's function is related to a system, rather than to a single, Radon-Nikodym derivative. To this end, several computations were performed for the expected number of visits of triangular effective random walks about \textit{each} degree of freedom. Furthermore, the triangular Green's function is also dependent upon a different constant relating to approximating the infinite sums $   \underset{0 \leq k \leq + \infty}{\sum}     \big[ \textit{Degree of each vertex over the triangular lattice} \big]^{-k}  \mathrm{exp} \bigg[       - k \bigg\{ \mathrm{log} \bigg[ \frac{\textbf{P}^*_{E_1}}{\textbf{P}_{E_1} + \textbf{P}_{E_2}}         \bigg] - \lambda    \bigg\}         \bigg]            $, and $   \underset{0 \leq k \leq + \infty}{\sum}     \big[ \textit{Degree of each vertex over the triangular lattice} \big]^{-k}  \mathrm{exp} \bigg[       - k \bigg\{ \mathrm{log} \bigg[ \frac{\textbf{P}^*_{E_2}}{\textbf{P}_{E_1} + \textbf{P}_{E_2}}         \bigg] - \lambda    \bigg\}         \bigg]            $. Intuitively, the natural logarithm of the first and second Radkon-Nikodym derivatives in the system dictate the typical number of expected \textit{complete} visits for the effective triangular random walk.  For the remaining number of steps belonging to an \textit{incomplete} excursion, one computes $ \mathscr{E} \bigg\{      \mathrm{log} \bigg[   \frac{ \mathrm{d}  \textbf{P}^{*}_{E_2}}{ \mathrm{d}  \textbf{P}_{E_1} +  \mathrm{d}  \textbf{P}_{E_2}}     \bigg] \tau  - \lambda  \underset{1 \leq i \leq \tau}{\sum} \big| \widetilde{\mathscr{U}_i} \big|        \bigg\} + \mathscr{E} \bigg\{      \mathrm{log} \bigg[   \frac{ \mathrm{d}  \textbf{P}^{*}_{E_2}}{ \mathrm{d}  \textbf{P}_{E_1} +  \mathrm{d}  \textbf{P}_{E_2}}     \bigg] \tau^{\prime}  - \lambda  \underset{1 \leq i \leq \tau^{\prime}}{\sum} \big| \widetilde{\mathscr{V}_i} \big|        \bigg\} $. With the desired representation of the Green's function, one adapts arguments  obtained from the kinetic scaling limit over the square lattice, as discussed in the previous section.

\section{Declarations}

\subsection{Availability of data and materials}

Not applicable.

\subsection{Competing interests}

Not applicable.

\subsection{Funding}

Not applicable.

\subsection{Authors' contributions}

PR wrote the entire manuscript and performed several rounds of editing.

\subsection{Acknowledgments}

Not applicable.

\section{References}

\noindent [1] Beffara, V., Friedli, S., Velenik, Y. Scaling Limit of the Prudent Walk. \textit{Electron. Comm. Probab.} \textbf{15}: 44-58 (2010).

\bigskip

\noindent [2] Bousquet-Melou, M. Families of prudent self-avoiding walks. \textit{Journal of Combinatorial Theory, Series A} \textbf{117}: 313-344 (2010).

\bigskip

\noindent [3] Bousquet-Melou, M. Lattice paths and random walks. \textit{Course lecture notes}.

\bigskip

\noindent [4] Petrelis, N., Sun, R., Torri, N. Scaling limit of the uniform prudent walk. \textit{Elec. J. Probab.} \textbf{22}: 1-19 (2017).

\bigskip

\end{document}